\RequirePackage{fix-cm}
\documentclass[smallextended,envcountsect]{svjour3}       
\usepackage{geometry}                
\smartqed  
\usepackage{graphicx}
\usepackage{mathtools}
\usepackage{stmaryrd}
\usepackage{amsfonts}
\usepackage{amstext}
\usepackage{amssymb} 
\usepackage{amsbsy}
\renewcommand{\figurename}{Fig.}
\usepackage{subcaption}
\usepackage{float}
\usepackage{multicol}
\usepackage{multirow}
\usepackage[titletoc,toc,title]{appendix}
\usepackage[colorlinks=true]{hyperref} 
\hypersetup{bookmarksnumbered,linkcolor = blue}  
\usepackage{threeparttable}
\newcommand{\heqref}[1]{Eq.~\hyperref[#1]{\eqref{#1}}} 
\newcommand{\heqsref}[2]{Eqs.~\hyperref[#1]{\eqref{#1}} and~\hyperref[#2]{\eqref{#2}}} 
\newcommand{\heqssref}[3]{Eqs.~\hyperref[#1]{\eqref{#1}},~\hyperref[#2]{\eqref{#2}} and~\hyperref[#3]{\eqref{#3}}} 
\newcommand{\hfigref}[1]{\hyperref[#1]{Fig.~\ref{#1}}}
\newcommand{\hfigsref}[2]{\hyperref[#1]{Figs.~\ref{#1}} and~\hyperref[#2]{\ref{#2}}}
\newcommand{\happref}[2]{\hyperref[#1]{Appendix~#2}}
\newcommand{\htablref}[1]{\hyperref[#1]{Table~\ref{#1}}}
\newcommand{\htablsref}[2]{\hyperref[#1]{Tables~\ref{#1}} and~\hyperref[#2]{\ref{#2}}}
\newcommand{\hccite}[1]{\hyperref[#1]{\cite{#1}}}
\newcommand{\secref}[1]{Section~\ref{#1}}

\newcommand{\diff}{\mathrm{d}}
\newcommand{\px}{\partial x}
\newcommand{\dx}{\mathrm{d}x}
\newcommand{\pxi}{\partial \xi}
\newcommand{\dxi}{\mathrm{d}\xi}

\newcommand{\pp}[1]{\partial{#1}}
\DeclareMathOperator{\Omg}{\Omega}
\DeclareMathOperator{\Omge}{\Omega_{e}}
\DeclareMathOperator{\Omgs}{\Omega_{s}}

\DeclareMathOperator{\pOmge}{\partial \Omega_{e}}

\newcommand{\ita}[1]{\textit{#1}}
\newcommand{\bff}[1]{\textbf{#1}}
\newcommand{\bs}[1]{\boldsymbol{#1}}
\newcommand{\undr}[1]{\underline{#1}}
\newcommand{\mytrue}{'\ita{true}' }
\newcommand{\combined}{\ita{combined-mode} }
\newcommand{\physic}{physical-mode }
\newcommand{\jump}[1]{\llbracket #1 \rrbracket}
\newcommand{\aver}[1]{\{ \! \! \{ #1 \} \! \! \} }

\newcommand{\bigaver}[1]{\left \lbrace \! \! \! \left \lbrace#1 \right \rbrace \! \! \! \right \rbrace }
\newcommand{\mymat}[1]{\undr{\bs{#1}  } }
\newcommand{\mycalmat}[1]{\undr{\bs{\mathcal{#1} } } }   
\begin{document}
\title{On the Accuracy and Stability of Various DG Formulations for Diffusion}
\author{Mohammad Alhawwary \and  Z.J. Wang }
\institute{Mohammad Alhawwary (\email{mhawwary@ku.edu}) \and
           Z.J. Wang (\email{zjw@ku.edu}) \at
           Department of Aerospace Engineering, University of Kansas, Lawrence, KS~66045, USA 
}

\date{Received: date / Accepted: date}
\maketitle

\begin{abstract}
In this paper, we study the stability (in terms of the maximum time step) and accuracy (in terms of the wavenumber-diffusion properties) for several popular discontinuous Galerkin (DG) viscous flux formulations. The considered methods include the symmetric interior penalty formulation (\bff{SIPG}), the first and second approaches of Bassi and Rebay (\bff{BR1}, \bff{BR2}), and the local discontinuous Galerkin method (\bff{LDG}). For the purpose of stability, we consider the von Neumann stability analysis method for uniform grids with a periodic boundary condition. In addition, the \ita{combined-mode} analysis approach previously introduced for the wave equation is utilized to analyze the dissipative error. This new approach can be used to study the performance of a particular DG and Runge-Kutta DG (RKDG) scheme for the entire extended wavenumber range. Thus, more insights into the robustness as well as accuracy and efficiency can be obtained. For instance, the \bff{LDG} method provides larger dissipation for high-wavenumber components than the \bff{BR1} and \bff{BR2} approaches for short time simulations in addition to a lower error bound for long time simulations. The \bff{BR1} approach with added dissipation can have desirable properties and stability similar to \bff{BR2}. For \bff{BR2}, the penalty parameter can be adjusted to enhance the performance of the scheme. The results are verified through canonical numerical tests. 
\end{abstract}

\section{Introduction}
The discontinuous Galerkin (DG) method, originally introduced by Reed and Hill~\cite{ReedTriangularMeshMethods1973} to solve the neutron transport equation, is arguably the most popular method in the class of adaptive high-order methods on unstructured grids. This class also includes several other methods such as the spectral difference (SD)~\cite{LiuSpectraldifferencemethod2006}, and the flux reconstruction (FR) or correction procedure via reconstruction (CPR) methods~\cite{HuynhFluxReconstructionApproach2007,WangUnifyingLiftingCollocation2009a,Wangreviewfluxreconstruction2016}. LaSaint and Raviart~\cite{LasaintFiniteElementMethod1974} performed an error analysis for the DG method. It was then further developed for convection-dominated problems and fluid dynamics by many researches, see for example(\cite{CockburnDiscontinuousGalerkinMethods2000,CockburnRungeKuttaDiscontinuous2001,BassiHighOrderAccurateDiscontinuous1997a,BassiHighOrderAccurateDiscontinuous1997,BassiHigherorderaccuratediscontinuous1997,HesthavenNodalDiscontinuousGalerkin2010,ShuHighorderWENO2016}) and the references therein. For elliptic problems, a number of formulations have been developed such as \bff{SIPG}~\cite{DouglasInteriorPenaltyProcedures1976,HartmannSymmetricInteriorPenalty2005}, \bff{BR1}~\cite{BassiHighOrderAccurateDiscontinuous1997a}, \bff{BR2}~\cite{BassiHigherorderaccuratediscontinuous1997}, \bff{LDG}~\cite{CockburnlocaldiscontinuousGalerkin1998} and the compact discontinuous Galerkin (\bff{CDG}) method~\cite{PeraireCompactDiscontinuousGalerkin2008}. In the work of Arnold et al.~\cite{ArnoldUnifiedAnalysisDiscontinuous2002}, the theoretical convergence, stability, and consistency of several  methods have been studied through a unified framework. However, the relative efficiency and dissipation properties of different schemes with respect to the wavenumber were not studied in the literature in great detail, see for example~\cite{GuoSuperconvergencediscontinuousGalerkin2013}. In~\cite{HuynhReconstructionApproachHighOrder2009}, Huynh conducted a Fourier analysis for those formulations in the context of flux reconstruction space discretization schemes in addition to the recovery method of Van Leer and Nomura~\cite{LeerDiscontinuousGalerkinDiffusion2005,vanLeerdiscontinuousGalerkinmethod2007,LeerUnificationDiscontinuousGalerkin2009}. Moreover, Kannan and Wang~\cite{KannanStudyViscousFlux2009} studied several viscous flux formulations for the Navier-Stokes equations using a p-multigrid spectral volume (SV)~\cite{WangSpectralFiniteVolume2002} solver. 

These discretization methods often include a penalty parameter that is utilized to ensure the stability. While the theoretical work of Arnold et al.~\cite{ArnoldUnifiedAnalysisDiscontinuous2002} and other works for the \bff{SIPG}~\cite{Shahbaziexplicitexpressionpenalty2005,EpshteynEstimationpenaltyparameters2007,OwensOptimaltraceinequality2017} method have provided guidance for choosing the penalty parameter to ensure stability, little insights about its effects on the dissipation error and the maximum time step exist. A recent attempt to provide the minimum values necessary for stability of a given energy stable flux reconstruction scheme (ESFR) using either \bff{SIPG} or \bff{BR2} methods is presented in the preprint~\cite{QuaegebeurStabilityEnergyStable2018} following the proofs of~\cite{WilliamsEnergystableflux2013,CastonguayEnergystableflux2013} for the \bff{LDG} method. Additionally, Gassner et al.~\cite{GassnerBR1SchemeStable2018} has studied theoretically the stability of the \bff{BR1} formulation with Gauss-Lobatto DG Spectral Element method (DGSEM) in energy stable split forms and Manzanero et al.~\cite{ManzaneroBassiRebayscheme2018} has showed that this particular version of the \bff{BR1} method is very close to the \bff{SIPG} method in some cases.  In this article, we study the wavenumber-diffusion properties of several methods as well as the effects of the penalty parameter on accuracy, efficiency and maximum time-step for fully-discrete RKDG schemes.

A semi-discrete analysis is performed first for the selected methods, followed by the fully-discrete analysis coupled with RK time integration schemes. In this regard, we investigate the effect of the penalty parameter on the maximum time-step required for stability. In addition, the dissipation properties of all schemes are studied using the~\combined analysis. Simplified closed-form expressions for a number of schemes are also derived which simplifies the implementation of these methods in 1D and establishes the similarities and connections among them. 

The paper is organized as follows.~\secref{sec:num_methods} introduces the basic formulation of the all numerical methods under consideration. We then conduct the semi-discrete analysis in~\secref{sec:sdiscAnalysis}, followed by the fully-discrete analysis for the RKDG schemes in~\secref{sec:fdiscAnalysis}. Numerical  verifications and test cases are presented in Section~\ref{sec:num_results}. Finally, conclusions are summarized in~\secref{sec:conclusions}. 

In this paper, matrices are denoted by either capitalized bold letters or bold math calligraphy letters both with an underscore, vectors are denoted by bold letters, while scalars 
are with plain letters. The columns of a matrix $\mycalmat{A}  \in \mathbb{C}_{m \times n} $ are denoted by $\bs{A}_{i} \in \mathbb{C}^{m}$, $i=1,...,n$. The scalar entries of the matrix are denoted by $A_{ij} \in \mathbb{C}$ such that $\bs{A}_{i}  = [A_{1j}, ..., A_{mj}]$. 

\section{Numerical Methods for Diffusion} \label{sec:num_methods}

In this section we present the basic formulation of the DG method for a one-dimensional linear parabolic diffusion equation of the form%
\begin{align}
   \frac{\partial u}{\partial t} - \gamma \frac{\partial^{2} u}{\partial x^{2}} = 0,  \quad x \in \mathcal{D}, \quad t >0
   \label{eqn:heat_eqn} 
\end{align}%
with periodic boundary conditions, where $u(x,t)$ is the solution in the physical domain, $\mathcal{D}$, and $\gamma$ is a positive constant representing the diffusivity coefficient. By introducing an auxiliary variable $\Theta(x,t)= \partial u/\partial x$, the model $2$nd order~\heqref{eqn:heat_eqn} can be written as a system of $1$st order equations%
\begin{subequations}
\begin{alignat}{2}
\Theta & = \frac{\partial u}{\partial x},  \label{eqn:theta_eqn_1stordersys}\\
\frac{\partial u}{\partial t}  & = \gamma \frac{\partial \Theta}{\partial x}.  
\label{eqn:u_eqn_1stordersys}
\end{alignat}
 \label{eqn:heat_eqn_1stordersys}
\end{subequations}

\subsection{The Discontinuous Galerkin (DG) method}  \label{sec:DG_num}

In the DG framework, the domain $\mathcal{D}$ in one-dimension is discretized into $N_{el}$ of non-overlapping elements, $\Omega_{e} = \left[x_{e-1/2}, x_{e+1/2}\right]$, such that $\mathcal{D}\:=\: \cup_{e=1}^{N_{el}} \Omega_{e}$, and each element has a variable width of $h_{e}$ and a center point $x_{e}$. In addition, DG assumes a reference element with local coordinate $\xi \in [-1,1]$, and defines a linear mapping between a physical and reference element of the form%
\begin{equation}
  \xi =2(x-x_{e})/h_{e}.
  \label{eqn:xi_mapping}
\end{equation}%
On element $\Omega_{e}$, the solution is approximated by a polynomial $u^{e}(x,t)$ of degree $p$ in space, i.e., $u^{e} \in \mathcal{P}^{p}$ which is a finite dimensional space of polynomials of degree at most $p \geq 1$. In addition, DG approximate $\Theta$ by $\Theta^{e}$, a polynomial that belongs to the same solution space of $\mathcal{P}^{p}$ in 1D, whereas for a $d$ dimensional space $\Theta^{e} \in [\mathcal{P}^{p}]^{d}$. The \textit{variational formulation} is obtained by introducing test functions $\upsilon, \psi \in \mathcal{P}^{p}$ and integrating the system of~\heqref{eqn:heat_eqn_1stordersys} over element $\Omega_{e}$%
\begin{subequations}\label{eqn:DG_variational_form}
\begin{alignat}{2}
\int_{\Omge} \Theta^{e} \upsilon \, dx & = \, \int_{\Omge} \frac{\partial u^{e}}{\partial x} \upsilon \: dx, \label{eqn:DG_variational_1steqn}\\
\int_{\Omge} \frac{\partial u^{e}}{\partial t} \psi \, dx & = \, \int_{\Omge} \gamma \frac{\partial \Theta^{e}}{\partial x}  \psi\: dx .
\label{eqn:DG_variational_2ndeqn}
\end{alignat}
\end{subequations}%
Applying integration by parts to the right hand side integrals in~\heqref{eqn:DG_variational_form}, we arrive at the \ita{weak formulation} for an element $\Omge$%
\begin{subequations} \label{eqn:DG_weakform_intbyparts}
\begin{alignat}{2}
\int_{\Omge} \Theta^{e} \upsilon \,  dx  & = \: \left[ u^{e} \upsilon \right]^{x_{e+1/2}}_{x_{e-1/2}} \: - \: \int_{\Omge} \frac{\partial \upsilon}{\partial x} u^{e} \, dx, & \qquad \forall \upsilon \in  \mathcal{P}^{p}, \label{eqn:DG_weak_1steqn_intbyparts} \\
\int_{\Omge} \frac{\partial u^{e}}{\partial t} \psi\: dx & = \gamma \left( -  \int_{\Omge} \: \Theta^{e} \frac{\partial \psi}{\partial x} \: dx \: + \: \left[\Theta^{e} \psi\right]^{x_{e+1/2}}_{x_{e-1/2}}  \right) ,  & \qquad \forall \psi \in  \mathcal{P}^{p}.
\label{eqn:DG_weak_2ndeqn_intbyparts}
\end{alignat}
\end{subequations}%
The interface terms in~\heqref{eqn:DG_weakform_intbyparts}, contain non-unique values for $u^{e}, \Theta^{e}$ at $x_{e-1/2}, x_{e+1/2}$. Hence, we define numerical values for them, i.e., $\hat{u}$ for $u^{e}$ and $\hat{\Theta}$ for $\Theta^{e}$. The DG equations become%
\begin{subequations}\label{eqn:DG_weakform_numflux}
\begin{alignat}{2}
\int_{\Omge} \Theta^{e} \upsilon \:  dx \: &= \: \left[ \hat{u} \upsilon \right]^{x_{e+1/2}}_{x_{e-1/2}} \: - \: \int_{\Omge} \frac{\partial \upsilon}{\partial x} u^{e} \: dx, \label{eqn:DG_1steqn_numflux} \\
\int_{\Omge} \frac{\partial u^{e}}{\partial t} \psi\: dx \: &= \: \gamma \left( -  \int_{\Omge} \: \Theta^{e} \frac{\partial \psi}{\partial x} \: dx \: + \: \left[\hat{\Theta} \psi\right]^{x_{e+1/2}}_{x_{e-1/2}}  \right),
\label{eqn:DG_2ndeqn_numflux}
\end{alignat}
\end{subequations}%
which is usually referred to as the~\textit{flux formulation}. The definition/choice of the numerical functions $\hat{u}, \hat{\Theta}$ at the interfaces has been studied in several papers~\cite{BrezziDiscontinuousGalerkinapproximations2000,ArnoldUnifiedAnalysisDiscontinuous2002,LeerDiscontinuousGalerkinDiffusion2005}. In this report we investigate several well-known DG formulations for the diffusion terms, namely, BR1~\cite{BassiHighOrderAccurateDiscontinuous1997a}, LDG~\cite{CockburnlocaldiscontinuousGalerkin1998}, SIPG~\cite{HartmannSymmetricInteriorPenalty2005}, and BR2~\cite{BassiHigherorderaccuratediscontinuous1997}.  

For the purpose of studying these methods, some preliminary definitions have to be made. In this regard, we define the average operator at an interface $\aver{...}$ and the jump operator across that interface $\jump{...}$ for any quantity $q$ as follows%
\begingroup
\allowdisplaybreaks
\begin{alignat} {2}
\begin{split}
\aver{q}  \, & = \, \frac{1}{2} \, \left(  q^{-} \, + \,  q^{+} \right), \quad \jump{ q }  = \, q^{+} \, \bs{n}^{+} \, + \,  q^{-} \, \bs{n}^{-} , 
\label{eqn:jump_average_operators}
\end{split}
\end{alignat}
\endgroup%
where $\bs{n}^{-} = - \bs{n}^{+}$ are the unit normal vectors and $(.)^{-} \& (.)^{+}$  denote the interface values of any quantity associated with the two elements sharing the interface in which $\bs{n}^{-}$ points outward of the element associated with $(.)^{-}$. We always consider the current element to be $\Omge=\Omg^{-}$ in this article. 

In order to unify the formulation of different DG methods for diffusion, Arnold et al.~\cite{ArnoldUnifiedAnalysisDiscontinuous2002} utilized the \ita{primal form} of the equations. For the first~\heqref{eqn:DG_1steqn_numflux}, let's perform integration by parts once again on the second integral at the RHS so that the result becomes%
\begin{equation}
\int_{\Omge} \Theta^{e} \: \upsilon \:  \dx \: = \: \left[ \left( \hat{u} - u^{e} \right) \: \upsilon \right]^{x_{e+1/2}}_{x_{e-1/2}} \: + \: \int_{\Omge} \frac{\partial u^{e} }{\partial x} \: \upsilon \: \dx.
\label{eqn:DG_strong_1steqn}
\end{equation}%
By selecting $\upsilon \, = \, \pp{\psi} / \pp{x} $, we then substitute the above equation into~\heqref{eqn:DG_2ndeqn_numflux} and rearrange to arrive at the \ita{primal formulation} %
\begin{equation}
\int\limits_{\Omge} \frac{\pp{u^{e}}}{\pp{t}} \psi\, \dx \, = \, \gamma \left( \left[ \hat{\Theta} \, \psi\right]^{x_{e+1/2}}_{x_{e-1/2}} \,  - \, \int\limits_{\Omge}  \frac{\partial u^{e}}{\px} \frac{\pp{\psi}}{\px} \, \dx \, - \left[ ( \hat{u}- u^{e} ) \, \frac{\pp{\psi}}{\px} \right]^{x_{e+1/2}}_{x_{e-1/2}}  \right).
\label{eqn:DG_u_primal_wolift}
\end{equation}%
It is customary in the literature to define a global lifting operator $R^{g}(\jump{u^{e}})$ using~\heqref{eqn:DG_strong_1steqn}, %
\begin{alignat}{2}
 \int\limits_{\Omge} R^{g} ( \jump{u^{e}} ) \,  \upsilon \, \dx  = \int\limits_{\Omge} \left( \Theta^{e} \, - \, \frac{\partial u^{e}}{\px} \right)  \, \upsilon \, \dx
& =  \, \left[ \left( \hat{u} - u^{e} \right) \, \upsilon \right]^{x_{e+1/2}}_{x_{e-1/2}} ,
\label{eqn:R_global_lift_1D_general}
\end{alignat}%
and $\Theta^{e}$ can be obtained by%
\begin{equation}
\Theta^{e}\: = \: \frac{\partial u^{e}}{\px} \, + \, R^{g} ( \jump{u^{e}} ).
\label{eqn:Theta_R_eqn_1D}
\end{equation}%
As a result, the \ita{primal form} reads%
\begin{equation} 
\int\limits_{\Omge} \frac{\partial u^{e}}{\partial t} \psi\, \dx \, \: = \: \gamma \:\left( \: \left[ \hat{\Theta} \: \psi \right]^{x_{e+1/2}}_{x_{e-1/2}}  - \: \int_{\Omge}  \frac{\pp{u^{e}}}{\px} \frac{\pp{\psi}}{\px} \: \dx \: - \int_{\Omge} R^{g}( \jump{u^{e}} ) \: \frac{\pp{\psi}}{\px}  \: \dx  \right).
\label{eqn:DG_u_primal_wlift}
\end{equation}%

For a standard element $\Omega_{s} = [-1,1]$, the solution polynomial $u^{e}$ can be constructed as a weighted sum of some specially chosen local basis functions $\phi \in \mathcal{P}^{p}$ defined on the interval $[-1,1]$%
\begin{equation}
u^{e}(\xi,t) = \sum_{j=0}^{p} U^{e}_{j}(t) \phi_{j}(\xi)
\label{eqn:u_h_xi_phi},
\end{equation}%
where the coefficients $U^{e}_{j}$ are the unknown degrees of freedom (DOFs). In a DG formulation, $\phi$ is chosen to be an orthogonal Legendre polynomial, for both test and basis functions, i.e., $\psi = \phi$.  As a result,~\heqref{eqn:DG_u_primal_wlift} constitutes a $p+1$ system of equations for the unknown DOFs $U^{e}_{j}, j=0,..,p$. The system of equations for $\Omgs$ can be written as%
\begin{equation}
\left[\frac{h_{e}}{\gamma} M_{ll}\right] \frac{\pp{U^{e}_{l}} }{\pp{t}}  \, =  \,   \, \left[ \hat{\Theta} \, \phi_{l} \right]^{1}_{-1} \,  - \, \frac{2}{h_{e}} \sum\limits_{j=0}^{p} S_{lj} U^{e}_{j} \, - \, \int_{\Omgs} R^{g}( \jump{u^{e}} ) \: \frac{ \pp{\phi}}{\pxi}  \: \dxi   , \quad l=0,...,p.
\label{eqn:DG_U_primal_wlift}
\end{equation}%
where $S_{lj}, M_{ll}$ are given by
\begin{equation}
S_{lj} = \int\limits_{-1}^{1} \frac{\partial \phi_{l}}{\pxi} \frac{\partial \phi_{j}}{\pxi } \dxi , \quad M_{ll} = \frac{1}{2} \int\limits_{-1}^{1} \phi_{l} \phi_{l} \dxi. 
\label{eqn:S_M_eqn}
\end{equation}

\subsection{Viscous flux formulations}\label{sec:viscflux_formulations}
\subsubsection{The Symmetric Interior Penalty approach (\bff{SIPG})}%
\label{sec:SIP_num}

In this section we introduce the so-called symmetric interior penalty formulation(\bff{SIPG})~\cite{DouglasInteriorPenaltyProcedures1976,HartmannSymmetricInteriorPenalty2005,Hartmannoptimalorderinterior2008} of the DG method for diffusion. The \bff{SIPG} method defines the numerical fluxes simply as%
\begin{subequations}\label{eqn:SIP_numfluxes}
\begin{equation}
\left.
\begin{aligned} 
\hat{u} \, & = \, \aver{u^{e}}  \\
\hat{\Theta} \, & = \, \bigaver{ \frac{\pp{u^{e}}}{\px} } \, - \,  \alpha^{ip} ( \jump{u^{e}} ) 
\end{aligned} \; \; \right \} \quad \text{at } f, 
\label{eqn:SIP_u_theta_hat}
\end{equation}
where $f$ is an interior face, and since we are interested in interior schemes and periodic boundary conditions this definition suffices. The lift operator $ \alpha \left( \jump{u^{e}} \right)$ at interface $f$ is approximated~\cite{HartmannSymmetricInteriorPenalty2005} by%
\begin{equation}
\alpha^{ip} \left( \llbracket u^{e}\rrbracket \right) = \eta_{f} \frac{C(p)}{\hslash_{f}}  \, \jump{u^{e}},
\label{eqn:SIP_alpha}
\end{equation}
\end{subequations}%
where $\hslash_{f}$ is a carefully selected length-scale, $\eta_{f}$ is a positive number large enough for stability, and $C(p)$ is a function of $p$. In Arnold et al.~\cite{ArnoldUnifiedAnalysisDiscontinuous2002} the lifting parameters were taken to be $\hslash_{f}= h_{f}$ a face length-scale, and $C(p)=1$. Hartmann et al.~\cite{HartmannSymmetricInteriorPenalty2005} suggested $\hslash_{f}= \min(h_{e},h_{f}), \, C(p)=p^{2}$ for the Navier-Stokes equations. Shahbazi~\cite{Shahbaziexplicitexpressionpenalty2005} provided an explicit expression for the penalty term for simplex elements ~\cite{DolejsiAnalysisdiscontinuousGalerkin2005,EpshteynEstimationpenaltyparameters2007}. Moreover, Hillewaert~\cite{HillewaertDevelopmentDiscontinuousGalerkin2013} defined minimum sharp values for $\eta_{f} C(p)$ to guarantee stability for different element types in 2D and 3D and so did Drosson~\cite{Drossonstabilitysymmetricinterior2013}. For an edge in a quadrilateral element, the sharp value for $ \eta_{f} C(p) $  was derived as $(p+1)^{2}$~\cite{HillewaertDevelopmentDiscontinuousGalerkin2013}. In addition, Hillewaert indicated that these parameters affect the conditioning of the implicit system and the minimum of $\eta_{f} C(p)$ results in a suitable conditioning for the system. In this work we choose $C(p)=(p+1)^{2}/2$ and $\hslash_{f}=h_{e}$ which results in a particular \bff{SIPG} method that shares interesting properties with the \bff{BR2} method, as shown later in this article. 

Using the definitions of the numerical fluxes in~\heqref{eqn:SIP_numfluxes}, $ R^{g} ( \jump{u^{e}} )$ can be computed as%
\begin{equation}
 \int\limits_{\Omge} R^{g} ( \jump{u^{e}} ) \,  \frac{\pp{\phi} }{\px} \, \dx  =  \, \left[ \left( \hat{u} - u^{e} \right) \, \frac{\pp{\phi} }{\px} \right]^{x_{e+1/2}}_{x_{e-1/2}} = -\frac{1}{2} \left[  \jump{u^{e}} \frac{\pp{\phi} }{\px} \right]^{x_{e+1/2}}_{x_{e-1/2}} ,
\label{eqn:R_global_lift_1D}
\end{equation}

and by combining the \ita{primal form}~\heqref{eqn:DG_U_primal_wlift} and using~\heqsref{eqn:SIP_numfluxes}{eqn:R_global_lift_1D}, we obtain%
\begingroup
\allowdisplaybreaks
\begin{alignat}{2}
\begin{split}
 \left[\frac{h_{e}}{ \gamma} M_{ll} \right] \frac{\pp U^{e}_{l} }{\pp t} & =  \left[ \bigaver{ \frac{\partial u^{e}}{\px} }_{e+1/2} \, \phi_{l}(1) \, - \, \bigaver{ \frac{\partial u^{e}}{\px} }_{e-1/2} \, \phi_{l}(-1) \right]  \\ & -  \, \eta_{f} \, \frac{C(p)}{h_{e}} \, \Big[ \jump{u^{e}}_{e+1/2} \phi_{l}(1) \, - \, \jump{u^{e}}_{e-1/2} \phi_{l}(-1)  \Big] \\ &  - \, \frac{2}{h_{e}}  \sum\limits_{j=0}^{p} S_{lj} U^{e}_{j}   \\ & + \, \frac{1}{2} \, \left[ \jump{u^{e}}_{e+1/2} \, \frac{\pp \phi_{l}}{\px}(1) \, + \, \jump{u^{e}}_{e-1/2} \, \frac{\pp \phi_{l} }{\px}(-1) \right], \quad l=0,...,p.
\end{split}
\label{eqn:SIP_expandedPrimal}
\end{alignat}
\endgroup%
Note that the \bff{SIPG}~\heqref{eqn:SIP_expandedPrimal} results in a fully-compact scheme, i.e., the elements that are involved in evaluating the time derivative in the above equation are from elements $e-1,\, e, \, e+1$ only. The 1D system can be written in a vector form as%
\begin{equation}
\frac{\partial \mathbf{U}^{e}}{\partial t} = \frac{\gamma}{h^{2}_{e}} \left(\mycalmat{K}^{e-1} \mathbf{U}^{e-1} + \mycalmat{L} \mathbf{U}^{e} +\mycalmat{K}^{e+1} \mathbf{U}^{e+1}   \right) , 
\label{eqn:SIP_vectorform}
\end{equation}%
where $\mathbf{U}^{e}=[U^{e}_{0},...U^{e}_{p}]$ is the vector of unknown DOFs and $\mycalmat{L}, \mycalmat{K}$ are constant matrices. 

\subsubsection{The first approach of Bassi and Rebay (\bff{BR1})}\label{sec:BR1_num}

This method was utilized by Bassi and Rebay~\cite{BassiHighOrderAccurateDiscontinuous1997a} with the RKDG method for the Navier-Stokes equations. In this method, the numerical fluxes at an interface $f$ are defined as%
\begin{subequations}\label{eqn:BR1_numfluxes} 
\begin{equation}
\begin{aligned}
\hat{u} \, & = \, \aver{u^{e}}, \\
\hat{\Theta}  \, & = \, \aver{ \Theta^{e} } \, = \, \bigaver{ \frac{\pp{u^{e}}}{\px} } - \alpha^{br1}(u^{e}), 
\end{aligned}
\label{eqn:BR1_u_theta_hat} 
\end{equation}

where 
\begin{equation}
\alpha^{br1}(u^{e})  = - \aver{ R^{g}(\jump{ u^{e} } ) }.
\label{eqn:BR1_alpha}
\end{equation}
\end{subequations}%
Substituting~\heqsref{eqn:R_global_lift_1D}{eqn:BR1_numfluxes} in the \ita{primal form}~\heqref{eqn:DG_U_primal_wlift}, we get%
\begin{alignat}{2}
\begin{split}
\left[\frac{h_{e}}{ \gamma} M_{ll} \right]  \frac{\pp U^{e}_{l} }{\pp t}  & =  \left[ \bigaver{ \frac{\partial u^{e}}{\px} }_{e+1/2}  \, \phi_{l}(1) \, - \, \bigaver{  \frac{\partial u^{e}}{\px} }_{e-1/2}  \, \phi_{l}(-1) \right] \\ 
& + \bigg[ \aver{ R^{g} }_{e+1/2} \, \phi_{l}(1) \, - \, \aver{ R^{g} }_{e-1/2} \, \phi_{l}(-1) \bigg]  \\ &  - \,\frac{2}{h_{e}} \sum\limits_{j=0}^{p} S_{lj} U^{e}_{j}   \\ & + \, \frac{1}{2} \, \left[ \jump{ u^{e}}_{e+1/2} \, \frac{\pp \phi_{l}}{\px}(1) \, + \, \jump{ u^{e} }_{e-1/2} \, \frac{\pp \phi_{l}}{\px}(-1) \right], \quad l=0,...,p.
\end{split}
\label{eqn:BR1_system_primal_1D}
\end{alignat}%
A further simplification of the update equation is obtained by computing the averages of $R^{g}$ at the element interfaces using~\heqref{eqn:R_global_lift_1D_general} and following ~\heqref{eqn:R_global_lift_1D}%
\begin{alignat}{2}
\begin{split}
\aver{ R^{g}  }_{e+1/2} & = -\frac{C(p)}{h} \jump{u^{e}}_{e+1/2} + \frac{C\!1(p)}{h} \left( \jump{u^{e}}_{e+3/2}  + \jump{u^{e}}_{e-1/2} \right)   \\
\aver{ R^{g} }_{e-1/2} & =- \frac{C(p)}{h} \jump{u^{e}}_{e-1/2} + \frac{C\!1(p)}{h}  \left( \jump{u^{e}}_{e+1/2}  + \jump{u^{e}}_{e-3/2}  \right) 
\end{split}
\label{eqn:BR1_Rg_simplified}
\end{alignat}%
where 
\begin{equation}
C\!1(p) = \frac{(-1)^{p+1}}{4} (p+1). 
\end{equation}%
Unfortunately, this method is weekly unstable and non-compact~\cite{BassiHigherorderaccuratediscontinuous1997,BrezziDiscontinuousGalerkinapproximations2000,ArnoldUnifiedAnalysisDiscontinuous2002,HuynhReconstructionApproachHighOrder2009}. This can be inferred by inspecting~\heqref{eqn:BR1_Rg_simplified}. In~\heqref{eqn:BR1_Rg_simplified} the computation of $R^{g}$ depends on not only the immediate neighbors of $\Omge$, namely $e-1,e+1$ but on their neighbors as well, i.e., $e-2,e+2$. In~\cite{GassnerBR1SchemeStable2018,ManzaneroBassiRebayscheme2018}, it has been shown that the \bff{BR1} method for the DGSEM with Gauss Lobatto points may possess the compactness and stability properties as the \bff{SIPG} method under certain conditions. %
%
\subsubsection{A modified \bff{BR1} method}\label{sec:BR1_stabilized}

If a penalty term is added to $\hat{\Theta}$ in~\heqref{eqn:BR1_u_theta_hat} of the standard \bff{BR1} method similar to the penalty term in the case of the \bff{SIPG} method, the scheme becomes identical to the \bff{IP} approach of Kannan and Wang~\cite{KannanStudyViscousFlux2009}. Using this penalty term the numerical flux becomes%
\begin{equation}
\hat{\Theta} =\bigaver{ \frac{\pp{u^{e}}}{\px} } + \aver{ R^{g} } -   \eta_{f}^{br1} \frac{C\!2(p)}{h}  \jump{u^{e}} =   \bigaver{ \frac{\pp{u^{e}}}{\px} } - \alpha^{br1}_{s}(\jump{u^{e}}) .
\label{eqn:BR1_Thetahat_penalty}
\end{equation}
where $\alpha_{s}^{br1}$ is the new stabilization term. By choosing $C\!2(p) = C(p)$ then from~\heqref{eqn:BR1_Rg_simplified} we can see that the terms involving $\jump{u^{e}}_{f}$ at the interface under consideration can be combined to simplify the formulation. This allows for easy comparison with the \bff{SIPG} and  \bff{BR2} which includes a similar jump term, as shown next. The new stabilized term $\alpha_{s}^{br1}(\jump{u^{e}})$ can be written as%
\begin{alignat}{2}
\begin{split}
\alpha_{s}^{br1}\big|_{e+1/2} & = \tilde{\eta}_{f}^{br1} \frac{C(p)}{h} \jump{u^{e}}_{e+1/2} - \frac{C\!1(p)}{h} \left( \jump{u^{e}}_{e+3/2}  + \jump{u^{e}}_{e-1/2} \right)   \\
\alpha_{s}^{br1}\big|_{e-1/2} & =  \tilde{\eta}_{f}^{br1} \frac{C(p)}{h} \jump{u^{e}}_{e-1/2} - \frac{C\!1(p)}{h}  \left( \jump{u^{e}}_{e+1/2}  + \jump{u^{e}}_{e-3/2}  \right) 
\end{split}
\label{eqn:BR1_stabilized_Rg_simplified}
\end{alignat}%
where $\tilde{\eta}_{f}^{br1} = 1 + \eta_{f}^{br1}$. The only difference between this modified \bff{BR1} scheme and the original one is the additional penalty parameter $\eta_{f}^{br1}$. In the rest of the paper, the name \bff{BR1} implies the stabilized version of the \bff{BR1} approach while for $\eta_{f}^{br1}=0$ we recover the standard \bff{BR1} method.  

\subsubsection{The second approach of Bassi and Rebay (\bff{BR2})}%
\label{sec:BR2_num}

The \bff{BR2} scheme originally introduced by Bassi and Rebay~\cite{BassiHigherorderaccuratediscontinuous1997,BassiNumericalevaluationtwo2002} as a modification for their \bff{BR1} scheme~\cite{BassiHighOrderAccurateDiscontinuous1997a}. The basic idea of the \bff{BR2} scheme is to define a local lift operator $r^{br2}_{f}\big( \jump{u^{e}} \big)$ for each face $f$ as a polynomial of degree $p$ instead of the global lift operator $R^{g}\big(\jump{u^{e}}\big)$. In~\heqref{eqn:DG_U_primal_wlift}, $r^{br2}_{f}\big( \jump{u^{e}} \big)$ is utilized only for the surface integral while $R^{g}$ in~\heqref{eqn:R_global_lift_1D} is utilized for the volume integral. This way the \bff{BR2} is very similar to the \bff{SIPG} method but with different stabilization term $\alpha^{br2} (u^{e} )$. The numerical fluxes are defined as%
\begin{subequations}\label{eqn:BR2_numfluxes} 
\begin{equation}
\begin{aligned}
\hat{u} \, & = \, \aver{u^{e}}, \\
\hat{\Theta} \, & = \,\bigaver{ \frac{\pp{u^{e}}}{\px} } \, - \,  \alpha^{br2}( u^{e} ) ,
\end{aligned}
\label{eqn:BR2_u_theta_hat}
\end{equation}
where%
\begin{equation} 
\alpha^{br2} ( u^{e} ) \, = \, - \eta^{br2}_{f}  \aver{ r^{br2}_{f} \big( \jump{u^{e}} \big) } , \quad \eta^{br2}_{f} \geq N_{f},
\label{eqn:BR2_alpha}
\end{equation}
\end{subequations}%
and $N_{f}$ is the number of faces of $\Omge$. In 1D this is just $N_{f}=1$. The local lift operator $r_{f}=r^{br2}_{f} \big( \jump{u^{e}} \big)$ at an interface $f$ is given by%
\begin{equation*}
\int\limits_{\Omge^{\pm}} r_{f}^{\pm} \upsilon \, \dx  =  \int\limits_{f \, \in \, \pOmge}  ( \hat{u}- u^{\pm} )  \bs{n}_{f}^{\pm} \,  \upsilon \, \diff s \,  =    - \frac{1}{2} \int\limits_{f \, \in \, \pOmge}  \jump{u^{e}} \upsilon \, \diff s,
\label{eqn:BR2_locaLiftIntegral}
\end{equation*}%
and on $\Omgs$, the polynomial definition of $r_{f}$ is%
\begin{equation} 
  r_{f} \,  = \, \sum\limits_{i=0}^{p} \, \tilde{r}_{f,i} \, \upsilon_{i}  ( \xi ). 
  \label{eqn:BR2_r_polynomial}
  \end{equation}%
If we take $\upsilon =\phi \bs{n}$, we get the normal local lift coefficient $r_{n}$ required for the surface integral in~\heqref{eqn:DG_U_primal_wlift} by performing the following integration%
\begin{equation}
\int\limits_{\Omge^{\pm}} r_{f}^{\pm} \, \phi^{\pm} \bs{n}^{\pm}  \, \dx = \int\limits_{\Omge^{\pm}} r_{n}^{\pm} \, \phi^{\pm}  \, \dx  = -\frac{1}{2}  \Big[ \jump{u^{e}} \, \phi^{\pm} \bs{n}^{\pm}  \Big]_{f} . 
\end{equation}%
%
Accordingly, in 1D, the \bff{SIPG} and \bff{BR2} DG methods become closely related to the extent that they reduce to the same update equation for the DOFs~\heqref{eqn:SIP_expandedPrimal}. It turns out that the term involving $\eta_{f} C(p) / \hslash_{f}$ in the surface integral of the numerical flux, results in a simple form for the \bff{BR2} method that is equivalent to the \bff{SIPG}. This was previously pointed out by Huynh~\cite{HuynhReconstructionApproachHighOrder2009} and Quaegebeur et al.~\cite{QuaegebeurStabilityEnergyStable2018} has mathematically proved the same result for the FR-DG schemes. In~\heqref{eqn:SIP_expandedPrimal}, the \bff{BR2} scheme is obtained by taking the stabilization parameters as%
\begin{equation} 
C(p) \: = \: \frac{(p+1)^{2}}{2}, \quad \eta^{br2}_{f}=\eta_{f}=1 , \quad \hslash_{f} = h_{e}.
\label{eqn:BR2_SIPG_equiv_Cp}
\end{equation} %
Therefore, in the rest of the paper we focus on the \bff{BR2} method. For $p=0$, these parameters result in an inconsistent \bff{BR2} scheme as reported by Huynh~\cite{HuynhReconstructionApproachHighOrder2009}. This can be remedied if $\eta_{f}=2$ for this case only. %
%
\subsubsection{The Local Discontinuous Galerkin method (\bff{LDG})}\label{sec:LDG_num}%
%
The \bff{LDG} was introduced by Cockburn and Shu~\cite{CockburnlocaldiscontinuousGalerkin1998} as a generalization for the \bff{BR1} scheme. In this formulation the numerical fluxes in 1D are given by%
\begin{subequations}\label{eqn:LDG_numfluxes}
\begin{alignat} {2}
\hat{u} & = \aver{ u^{e} } -\frac{\beta_{f}}{2} \jump{u^{e}} =  \left(1-\beta_{f} \right) u^{-} \, + \,  \beta_{f} \, u^{+} ,  \label{eqn:LDG_uhat}\\  
\hat{\Theta} & = \aver{ \Theta^{e} } + \frac{\beta_{f}}{2} \jump{ \Theta^{e} }  = \beta_{f} \, \Theta^{-} \, + \,  \left(1-\beta_{f} \right)  \Theta^{+}  \, - \, \frac{\eta_{f}^{ldg}}{h_{f}}  \jump{u^{e}},
\label{eqn:LDG_thetahat}
\end{alignat}
  \end{subequations}%
where $\beta_{f} \in \{0,1\}$ is a face switch and $\eta^{ldg}_{f}$ is a numerical stabilization (penalty) parameter. Substituting~\heqref{eqn:LDG_numfluxes} into~\heqref{eqn:DG_U_primal_wlift}, we obtain the following definition for $\Theta^{e}$%
\begin{alignat}{2}
 \Theta^{e} = \frac{\pp{ u^{e} } }{\px} + R^{g} \big( \beta \jump{u^{e}} \big), 
 \label{eqn:LDG_Theta_R_beta}
\end{alignat}%
where, for  $\Omge=\Omg^{-}$,%
\begin{alignat}{2}
\int\limits_{\Omge}  R^{g} \big( \beta \jump{u^{e}} \big) \upsilon \, \dx \, & = \oint\limits_{\pOmge} ( \hat{u}- u^{-} ) \, \bs{n}^{-} \upsilon^{-} \, \diff s = - \oint\limits_{\pOmge} \beta_{f} \jump{u^{e}} \,  \upsilon^{-} \, \diff s, 
\end{alignat}%
and by selecting $\upsilon =\phi \bs{n}$ we get the normal $R^{g}_{n} \big( \beta \jump{u^{e}} \big)$%
\begin{equation}
\int\limits_{\Omge}  R^{g} \, \phi \, \bs{n}^{-} \, \dx \, = \int\limits_{\Omge}  R^{g}_{n} \,  \phi \, \dx  \,  =   - \sum\limits_{f \in \pOmge} \beta_{f} \big( \jump{u^{e}}  \phi^{-} \bs{n}^{-} \big)_{f} \, .
\end{equation}%
The above equation implies that the integral has support only at one face in 1D since one can choose either $\left(\beta_{e+1/2} = 1, \beta_{e-1/2}=0\right)$ or $\left(\beta_{e+1/2} = 0, \beta_{e-1/2}=1\right)$ and this makes it compact in 1D. Using these forms and selecting $\left(\beta_{e+1/2}=1, \, \beta_{e-1/2}=0\right)$, the primal form of the \bff{LDG} can be written as%
\begingroup
\allowdisplaybreaks
\begin{alignat}{2}
\begin{split}
\left[ \frac{h_{e}}{ \gamma} M_{ll} \right]  \frac{\pp U^{e}_{l} }{\pp t} \, & =  \left[   \frac{\pp{u^{e}} }{\px}\Big|_{e+1/2} \, \phi_{l}(1) -  \frac{\pp{u^{e-1}} }{\px}\Big|_{e-1/2} \, \phi_{l}(-1) \right]  \\ & - \frac{2 C(p) + \eta^{ldg}_{f}}{h_{e}}   \Big[   \jump{u^{e}}_{e+1/2} \, \phi_{l}(1) \,  -  \, \jump{u^{e}}_{e-1/2} \, \phi_{l}(-1) \Big] \\ 
 & - \, \frac{2}{h_{e}}  \sum\limits_{j=0}^{p} S_{lj} U^{e}_{j}  + \,  \jump{u^{e}}_{e+1/2} \, \frac{\pp{\phi_{l}}}{\px}(1), \quad l=0,...,p , 
\end{split}
\label{eqn:LDG_primal_1D}
\end{alignat}
\endgroup%
with $C(p) = (p+1)^{2}/2$ as in the \bff{SIPG} and \bff{BR2} approaches. %
%
\subsection{Time integration methods\label{sec:num_method_RK}}
%
For time integration we utilize the Runge-Kutta method which is applied to the ordinary differential equation (ODE) resulting from the DG space discretization of~\heqref{eqn:heat_eqn}. This ODE takes the form%
\begin{equation}
\frac{d \bs{u} }{d t}= \gamma \mycalmat{A} \bs{u},
\label{eqn:ODE_heateqn}
\end{equation}%
where $\bs{u}$ is the vector of all unknown global element-wise DOFs, and $ \mycalmat{A}$ is the space discretization operator. Applying a RK scheme to~\heqref{eqn:ODE_heateqn} results in an update formula for $\bs{u}$ at $t=t+\Delta t$ of the following form%
\begin{equation}
\bs{u}(t+\Delta t) = \left[ \mycalmat{I} + \sum_{m=1}^{s} \frac{\left(\gamma \Delta t \: \mycalmat{A} \right)^{m}}{m!} \right] \bs{u}(t)  = \mathcal{P}(\gamma \Delta t  \mycalmat{A}) \bs{u}(t) ,
\label{eqn:RK_update_form_DG}
\end{equation}%
where $\mycalmat{I}$ is the identity matrix and $s=2$ for second-order RK2~\cite{GottliebStrongStabilityPreservingHighOrder2001}, $s=3$ for third-order RK3~\cite{GottliebStrongStabilityPreservingHighOrder2001}, and $s=4$ for the fourth-order classical RK4~\cite{ButcherNumericalAnalysisOrdinary1987}, and $\mathcal{P}$ is a polynomial of degree $s$. %
%
\section{Semi-discrete Fourier Analysis } \label{sec:sdiscAnalysis}%

The linear parabolic diffusion equation defined in~\heqref{eqn:heat_eqn} with an initial wave solution
\begin{equation}
 u(x,0) = u_{o}(x) = e^{ikx} ,
 \label{eqn:linheat_IC}
\end{equation}%
admits a wave solution of the form %
\begin{equation}
u(x,t) = e^{i k x -\omega t},
\label{eqn:lin_diffusion_exact_sol}
\end{equation}%
where $k$ is the spatial wavenumber, and $\omega$ denotes the frequency that admits the exact dissipation relation $\omega=\gamma k^{2}$. In a temporal Fourier analysis~\cite{VichnevetskyFourierAnalysisNumerical1982} a prescribed wavenumber $k$ is assumed for the initial condition~\heqref{eqn:linheat_IC} and different spatial and temporal schemes are applied to~\heqref{eqn:heat_eqn} in order to study their diffusion properties based on the numerical frequency $\tilde{\omega}$. 

The semi-discrete analysis utilizes only the spatial discretization for~\heqref{eqn:heat_eqn} so that dissipation properties can be studied. The method used in this section for DG is similar to the one previously presented in~\cite{HuAnalysisDiscontinuousGalerkin1999,MouraLineardispersiondiffusion2015,AlhawwaryFourierAnalysisEvaluation2018} for the linear advection equation, and Guo et. al~\cite{GuoSuperconvergencediscontinuousGalerkin2013} for the linear heat equation among others. For DG schemes, applying the spatial discretization to ~\heqref{eqn:heat_eqn} results in a system of semi-discrete equations of the form~\heqref{eqn:SIP_vectorform}, in which, the element-wise DOFs $U^{e}_{j}$ are then related to the wave solution~\heqref{eqn:lin_diffusion_exact_sol} through projection onto the solution basis. Therefore, the element-wise DOFs $U^{e}_{l}$ are computed as%
\begin{equation}
U_{l}^{e}(0) = \frac{\int_{\Omge} u^{e}(x,0) \phi(x) dx }{ \int_{\Omge} \phi(x) \phi(x) dx }= \frac{\int_{\Omgs} u(x_{e}+\xi h/2,0) \phi_{l}(\xi) d\xi }{L_{ll}}.
\label{eqn:sdisc_initsolproj}
\end{equation}%
For the initial wave form ~\heqref{eqn:linheat_IC}, these DOFs can be written as%
\begin{equation}
U_{l}^{e}(0) = \hat{\mu}_{l} \: e^{i k x_{e} },
\label{eqn:lin_diffusion_sol_coef_blochwave}
\end{equation}%
where%
\begin{equation}
\hat{\mu}_{l} = \frac{\int_{-1}^{1}  e^{ik \left(\xi h/2\right) } \phi_{l}(\xi) d\xi }{L_{ll}}.
\label{eqn:mu_projection}
\end{equation}%
It is easy to see that the exact DG solution can be expressed as %
\begin{equation}
U_{l}^{e}(t) =  \mu_{l} e^{i k x_{e} - \omega t}.
\label{eqn:sdisc_waveform}
\end{equation}%
where $\mu_{l}=\hat{\mu}_{l}$ for the exact solution. Seeking a similar solution to~\heqref{eqn:lin_diffusion_sol_coef_blochwave} and substituting into~\heqref{eqn:SIP_vectorform} we get%
\begin{equation}
\frac{\partial \mathbf{U}^{e}}{\partial t} = \frac{\gamma}{h^{2}} \left[ \mycalmat{K}^{e-1} e^{-Ikh} + \mycalmat{L}+ \mycalmat{K}^{e+1} e^{Ikh} \right] \mathbf{U}^{e},
\label{eqn:sdisc_updateeqn}
\end{equation}%
and by differentiating this equation, we get the semi-discrete relation%
\begin{equation}
\lambda\, \bs{\mu} = \mycalmat{A} \, \bs{\mu}, \qquad \lambda = -\left( \frac{h^2}{\gamma} \tilde{\omega} \right).
\label{eqn:semidisc_rel}
\end{equation}%
Assuming a real wavenumber $k$ and a numerical frequency $\tilde{\omega}$, the semi-discrete system~\heqref{eqn:semidisc_rel} constitutes an eigenvalue problem with  $p+1$ eigenpairs, $(\lambda_{j},\bs{\mu}_{j})$, $j=0,..,p$, for the matrix $\mycalmat{A}$. It turned out that all the eigenvalues of the matrix  $\mycalmat{A}$ are real for all schemes under consideration.  The general solution~\heqref{eqn:sdisc_waveform} can be written as a linear expansion in the eigenvectors space%
\begin{equation}
\mathbf{U}^{e}(t) = \sum_{j=0}^{p} \vartheta_{j} \bs{\mu}_{j} \: e^{i k x_{e} - \tilde{\omega}_{j} t } .
\label{eqn:sdisc_gen_sol_eigenvector_space}
\end{equation}%
The expansion coefficients $\vartheta_{j}$ are obtained by satisfying the initial condition~\heqref{eqn:heat_eqn} for $t=0$. As a result, these coefficients are given by%
\begin{equation}
\hat{\bs{\mu}}= \sum_{j=0}^{p} \vartheta_{j} \bs{\mu}_{j}, \quad \text{or} \quad \bs{\vartheta} = \mymat{\mathcal{M}}^{-1}  \hat{\bs{\mu}},
\label{eqn:eigenmode_weights}
\end{equation}%
where $\mymat{\mathcal{M}}=\left[\bs{\mu}_{0},...,\bs{\mu}_{p}\right]$ is the matrix of eigenvectors, and $\bs{\vartheta}=\left[\vartheta_{0},...,\vartheta_{p}\right]^{T}$.  The numerical solution can now be computed using~\heqref{eqn:u_h_xi_phi} with the DOFs given by~\heqref{eqn:sdisc_gen_sol_eigenvector_space} to yield%
\begin{equation}
u^{e}(\xi,t) = \sum\limits_{l=0}^{p} \sum\limits_{j=0}^{p}  \vartheta_{j} \mathcal{M}_{l,j} \: e^{i k x_{e} - \tilde{\omega}_{j} t }  \phi_{l}(\xi),
\label{eqn:sdisc_fullsolform}
\end{equation}%
while the exact solution can be written as%
\begin{equation}
u^{e}_{ex}(\xi,t) = \sum\limits_{l=0}^{p} \hat{\mu}_{l} \: e^{i  k x_{e} - \omega t }  \phi_{l}(\xi).
\label{eqn:exact_sdisc_fullsolform}
\end{equation}%
The numerical solution~\heqref{eqn:sdisc_fullsolform} is essentially a linear combination of $p+1$ eigenmodes, each has its own numerical diffusion behavior. In order to assess the numerical dissipation of DG schemes based on a certain eigenmode, we define the non-dimensional wavenumber as%
\begin{equation}
K = +\sqrt{\frac{\omega}{\gamma}}  \left(\frac{h}{p+1}\right) =   \frac{kh}{p+1}, \quad 0 \leq kh \leq (p+1) \pi.
\end{equation}%
Similarly, by defining the modified/numerical wavenumber $k_{m} = +\sqrt{\tilde{\omega}/\gamma}$ we can obtain a similar expression for the non-dimensional $K_{m}$. Note that $h/(p+1)$ is a measure of the smallest length-scale that can be captured by a DG scheme~\cite{MouraLineardispersiondiffusion2015,AlhawwaryFourierAnalysisEvaluation2018}. The numerical scheme introduces a modified wavenumber $k_{m}$ which induces the numerical dissipation through $\operatorname{\mathcal{R}e}(k_{m})$, whereas the numerical dispersion for all schemes studied in this work is zero, since $\mathcal{I}m(k_{m}) = 0$. 

In the literature the ability of polynomial-based high-order methods to admit an eigensolution structure has been interpreted in different ways~\cite{HuAnalysisDiscontinuousGalerkin1999,VandenAbeeleDispersiondissipationproperties2007,VincentInsightsNeumannanalysis2011,MouraLineardispersiondiffusion2015}. However, they all agree on that there is one physically relevant (physical) mode that characterizes the behavior of the scheme as a whole while the other modes are considered parasite. Recently, the present authors~\cite{AlhawwaryFourierAnalysisEvaluation2018} proposed the \combined analysis whereby the \mytrue behavior of DG and RKDG schemes can be studied for the entire range of wavenumbers. In~\cite{AlhawwaryFourierAnalysisEvaluation2018} it was verified that the \physic of a DG and RKDG scheme serves as a good approximation for the \mytrue behavior  in the low wavenumber range $K \lesssim \pi/2$ for the linear advection equation. For the linear diffusion equation, the physical-mode idea cannot be extended readily and previous studies have focused on a subset of the extended wavenumber range ($0 \leq kh \leq (p+1)\pi$) as in~\cite{HuynhReconstructionApproachHighOrder2009,KannanStudyViscousFlux2009}. As a result, the~\combined\cite{AlhawwaryFourierAnalysisEvaluation2018} analysis is imperative for a complete and more informative analysis. %
%
\subsection{Eigenmode analysis}\label{sec:sdisc_eigenmode_analysis}%
\begin{figure}[H]
\begin{subfigure}[h]{0.5\textwidth}
    \includegraphics[width=0.975\textwidth]{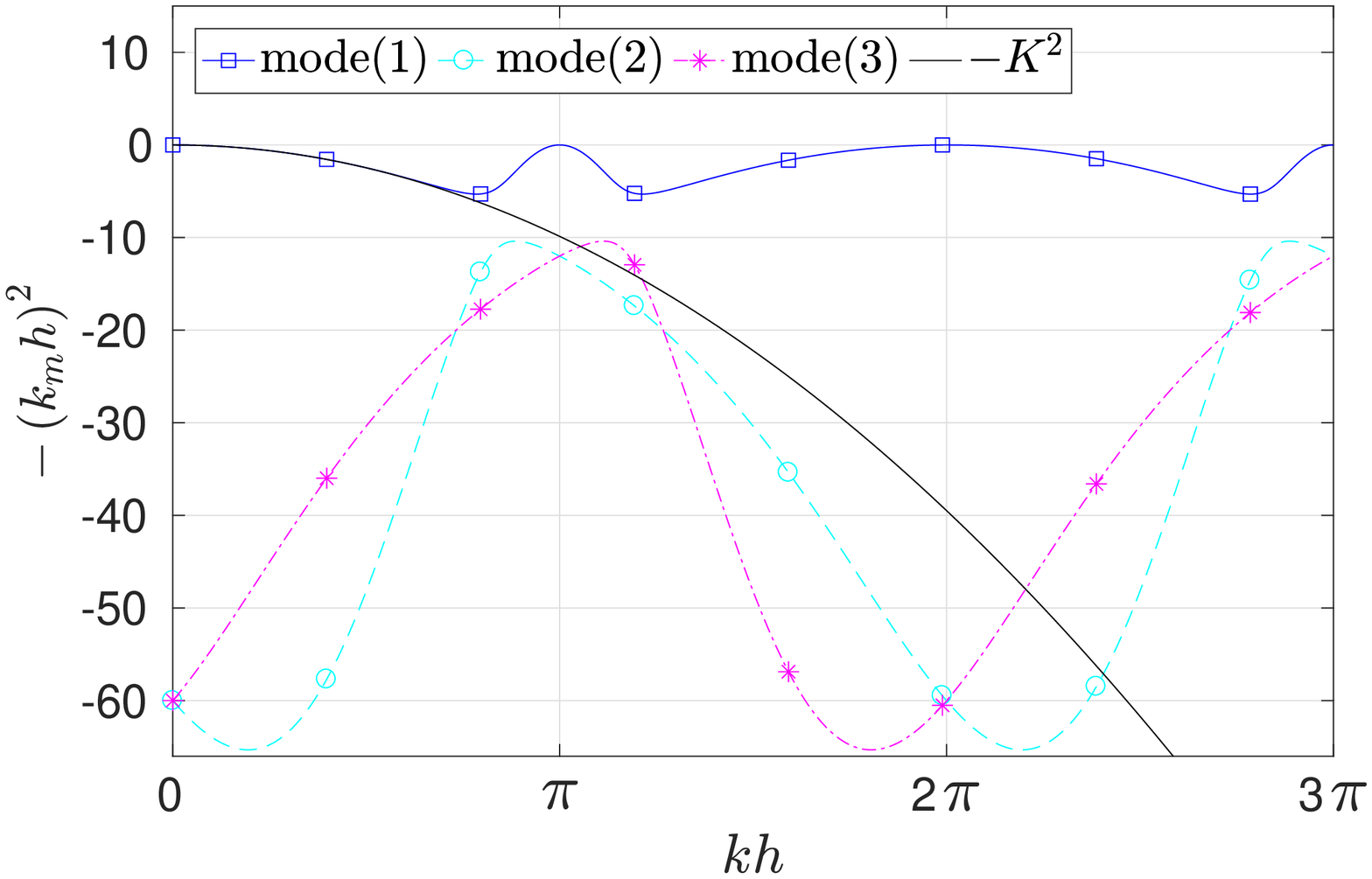} 
    \caption{Numerical dissipation of BR1p$2$-$\eta0$} 
    \label{fig:BR1p2eta0_wd}
    \end{subfigure}
      \, 
   \begin{subfigure}[h]{0.5\textwidth}
    \includegraphics[width=0.975\textwidth]{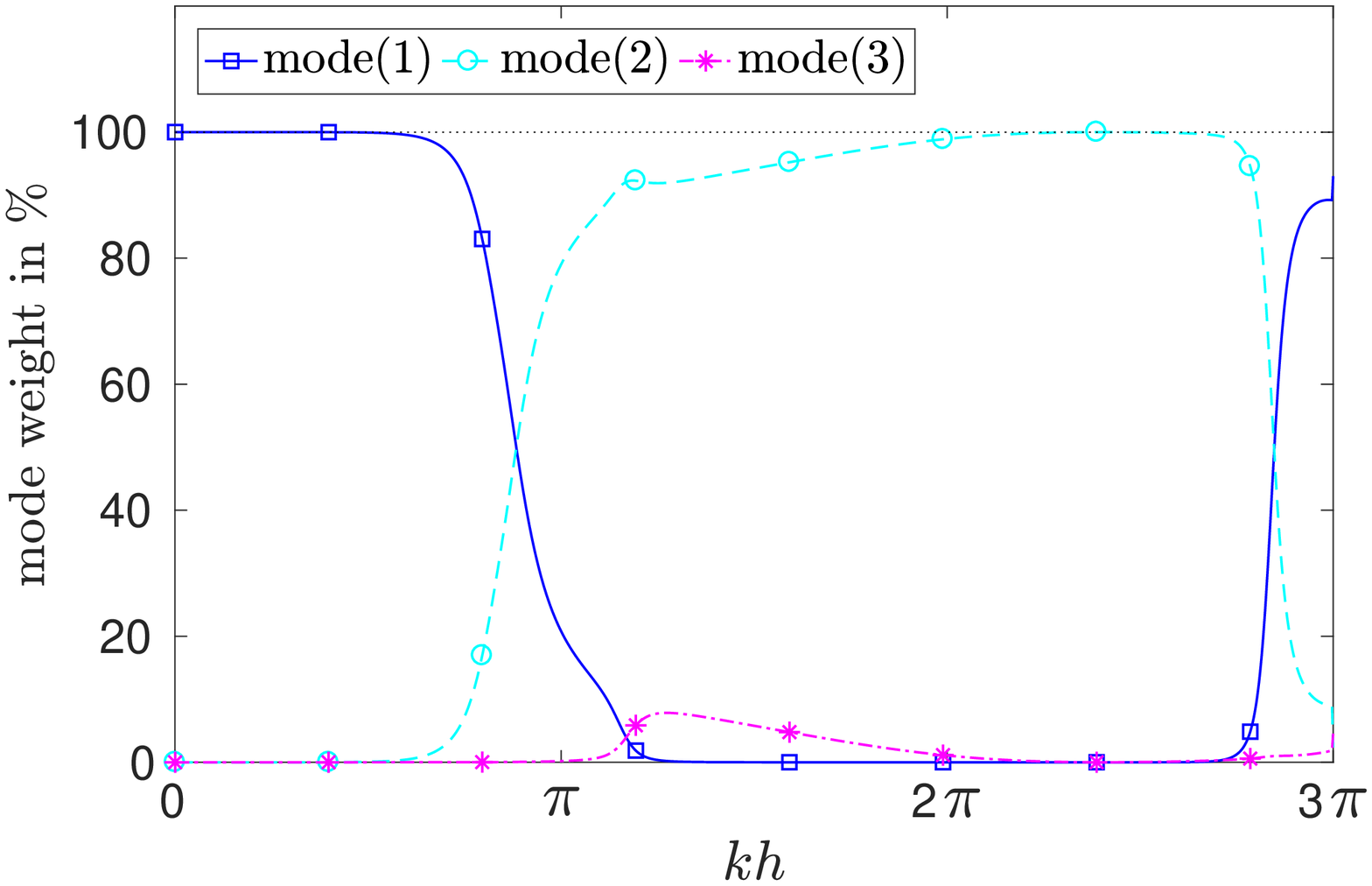}  
    \caption{Relative energy distribution, BR1p$2$-$\eta0$} 
    \label{fig:BR1p2eta0_weights}
    \end{subfigure}  \\ \\ \\
    \begin{subfigure}[h]{0.5\textwidth}
    \includegraphics[width=0.975\textwidth]{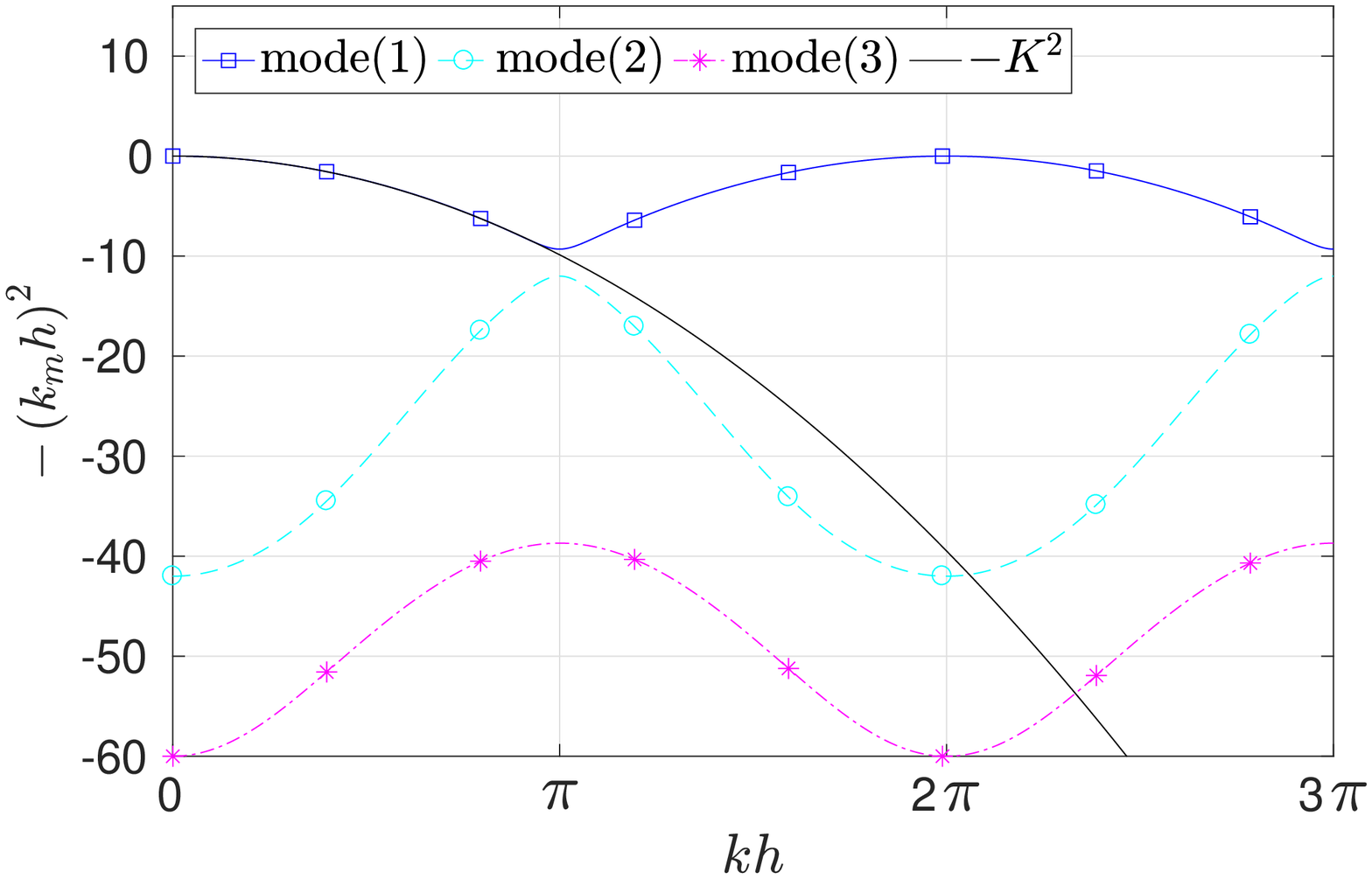} 
    \caption{Numerical dissipation of BR$2$p$2$-$\eta1$ } 
    \label{fig:BR2p2eta0_wd}
    \end{subfigure}
    \, 
   \begin{subfigure}[h]{0.5\textwidth}
    \includegraphics[width=0.975\textwidth]{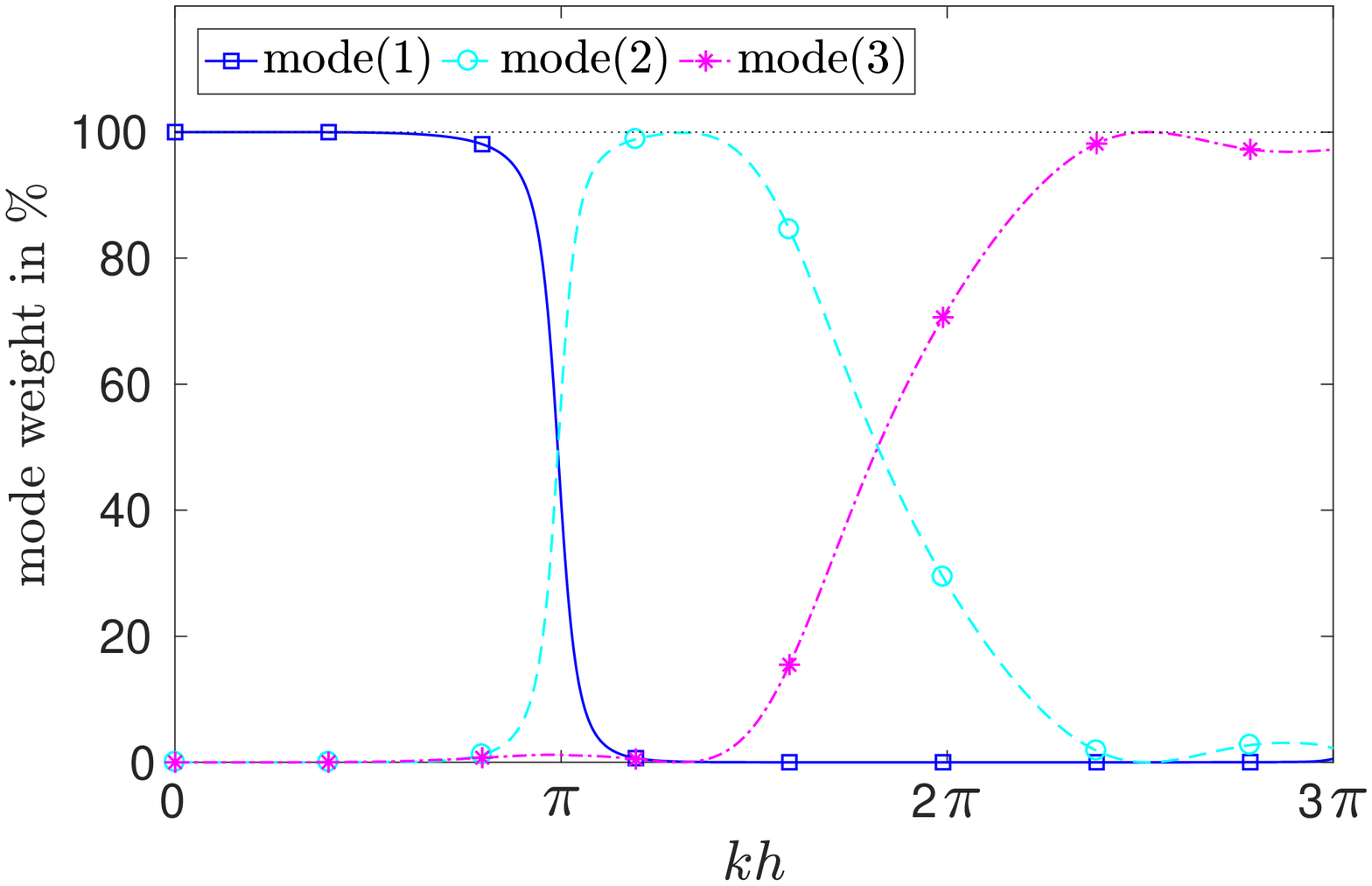}  
    \caption{Relative energy distribution, BR$2$p$2$-$\eta1$} 
    \label{fig:BR2p2eta0_weights}
    \end{subfigure}  \\ \\ \\
    \begin{subfigure}[h]{0.5\textwidth}
    \includegraphics[width=0.975\textwidth]{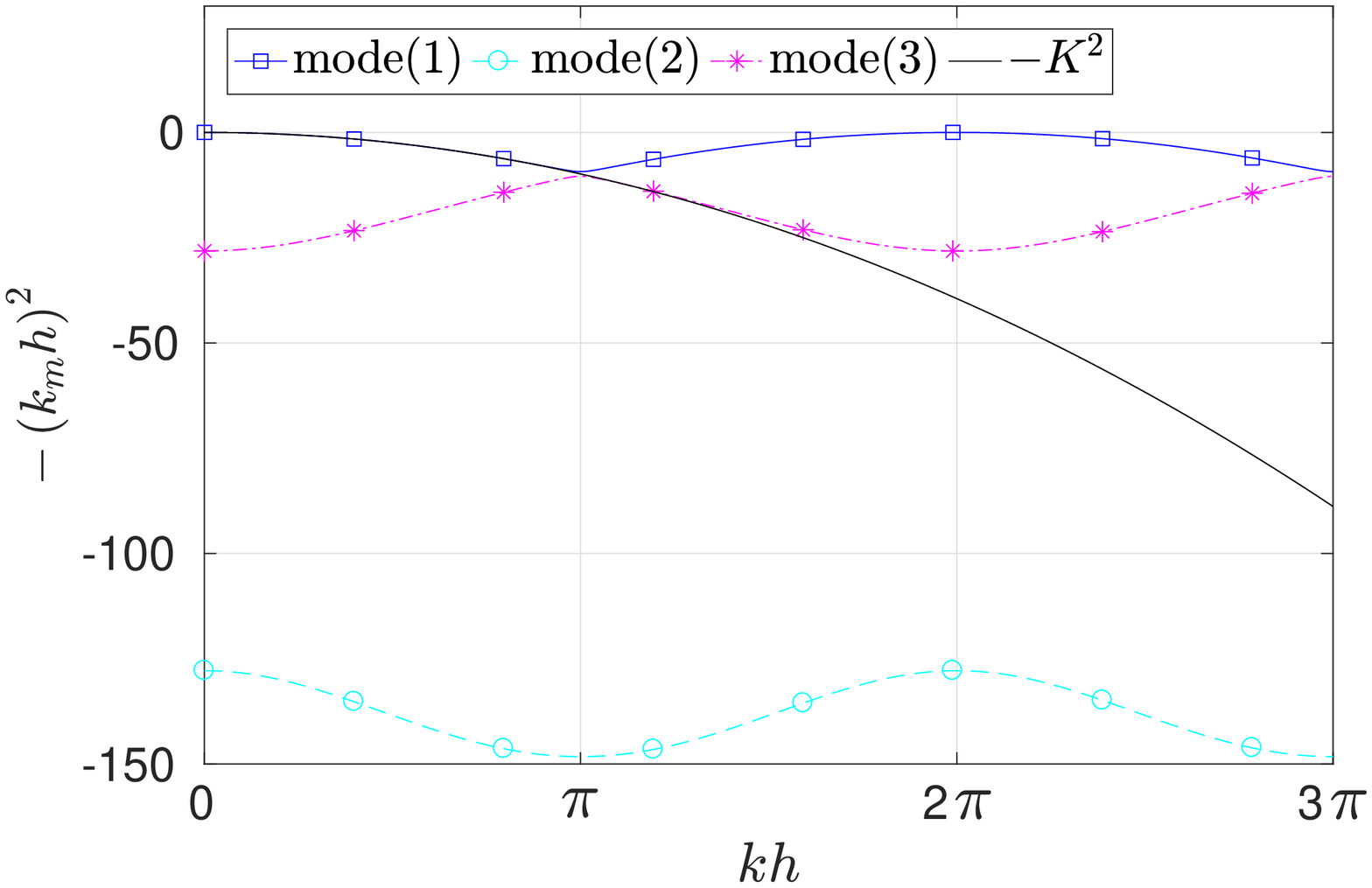} 
    \caption{Numerical dissipation of LDGp$2$-$\eta0$} 
    \label{fig:LDGp2eta0_wd}
    \end{subfigure}
    \, 
   \begin{subfigure}[h]{0.5\textwidth}
    \includegraphics[width=0.975\textwidth]{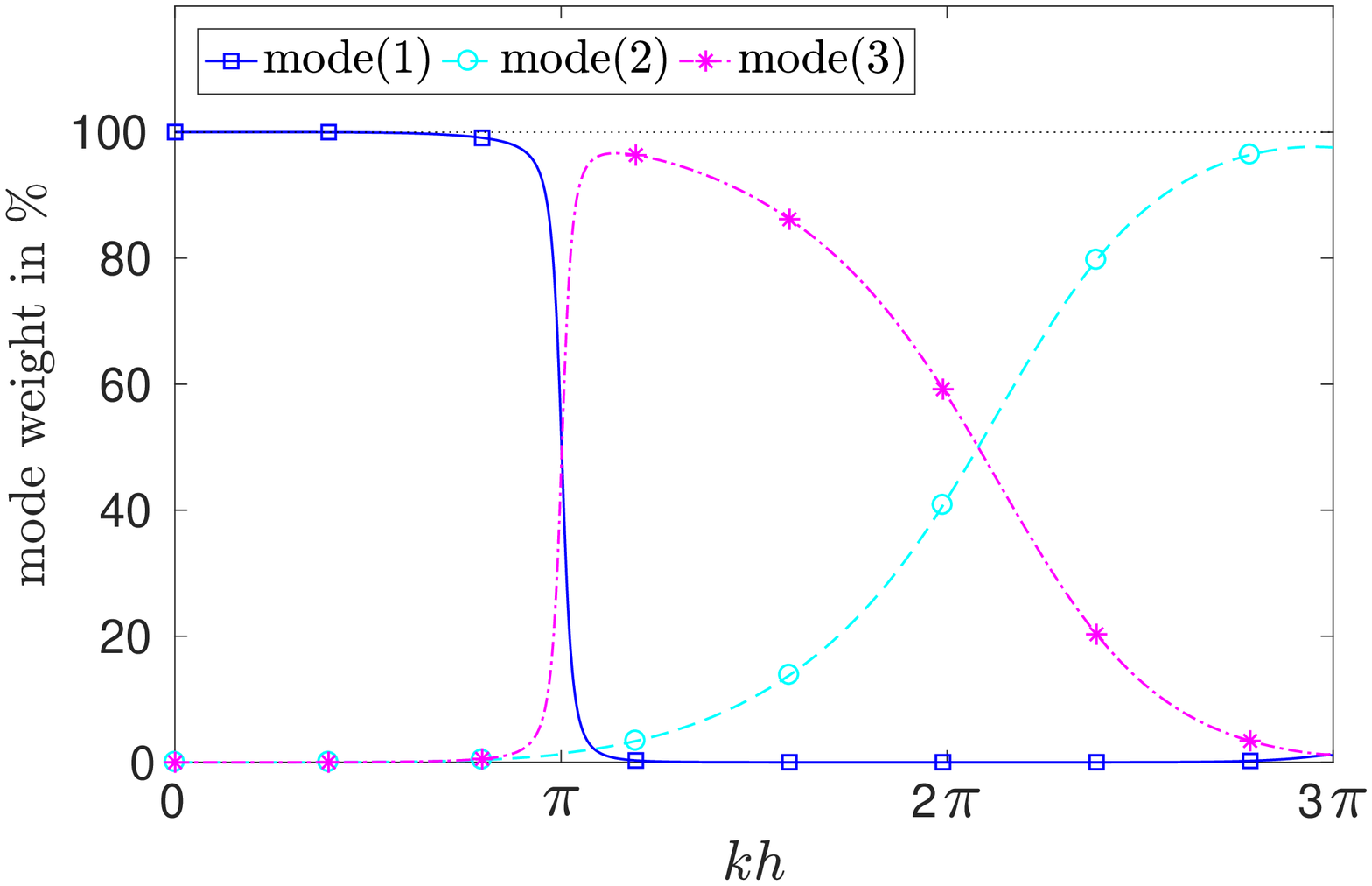}  
    \caption{Relative energy distribution, LDGp$2$-$\eta0$} 
    \label{fig:LDGp2eta0_weights}
    \end{subfigure}
	\caption{Semi-discrete dissipation and eigenmode relative energy curves for DGp$2$ schemes with \bff{BR1}, \bff{BR2}, and \bff{LDG} formulations. This plot contains all three eigenmodes of the scheme with mode(1) as the \physic.}
    \label{fig:sdisc_physical_mode_allschemes}
\end{figure}

%
In this section, we present the dissipation analysis based on the individual modes. By defining a relative energy measure~\cite{AsthanaHighOrderFluxReconstruction2015} among eigenmodes one can gain more insight about the influence of each mode along the extended wavenumber range. This relative energy measure is defined based on the weights $\vartheta$ of the normalized eigenvectors
\begin{equation}
\Gamma_{l} =  |\vartheta_{l}|^{2}/ \sum\limits_{j=0}^{p} |\vartheta_{j}|^{2}, \quad l=0,...,p,
\end{equation}%
where $\Gamma_{i}$ is the relative energy share of eigenmode $l$. ~\hfigref{fig:sdisc_physical_mode_allschemes} presents the eigenmode dissipation curves as well as relative energy distribution curves for the p$2$ spatial DG schemes with \bff{BR1}, \bff{BR2}, and \bff{LDG} formulations. In this figure, we have chosen to study the three formulations in their standard forms with minimum $\eta$. Note that in the rest of this paper we drop the subscript $f$ from $\eta_{f}$ since we are working with uniform grids and hence $\eta_{f}=\eta$, constant for all faces. In addition, if a particular order p and $\eta$ parameter are specified we use the following notation as an example \bff{BR1}p$2$-$\eta0$, \bff{BR2}p$2$-$\eta1$, and \bff{LDG}p$2$-$\eta0$ for the p$2$ spatial DG schemes. From~\hfigref{fig:sdisc_physical_mode_allschemes}, we can see that mode(1), in the left column of sub-figures, which closely agrees with the classical physical-mode definition~\cite{HuAnalysisDiscontinuousGalerkin1999,HuynhReconstructionApproachHighOrder2009,KannanStudyViscousFlux2009}, only approximates the exact dissipation up to $kh \approx \pi$. For multi-degrees of freedom high-order methods, the wavenumber range of the initial Fourier-mode can be extended up to $kh \leq (p+1)\pi$ due to their local nDOFs of ($p+1$) in each element. 

Therefore, it is expected that these schemes should approximate the exact dissipation for a reasonable subset of the extended wavenumber range. Unfortunately, through~\hfigref{fig:sdisc_physical_mode_allschemes} this is not true and this problem was not carefully addressed in the literature. If one examines the relative energy dissipation of each eigenmode in the right column of subfigures~\hfigref{fig:sdisc_physical_mode_allschemes}, it can be seen clearly that mode(1) has the highest energy share up to $kh\approx \pi$ and then the highest energy is contained in another mode which in turn loses it after some range to the third mode. This means that the coupling between modes is very important and they all cooperate to approximate the exact dissipation relation for the entire range of wavenumbers. 

For the \bff{BR1}p$2$-$\eta0$ scheme, we can see that it is the only scheme in~\hfigref{fig:sdisc_physical_mode_allschemes} that has a mode with exactly zero dissipation at $kh = 3\pi$ and this raises several instability and robustness issues~\cite{BassiHigherorderaccuratediscontinuous1997,BrezziDiscontinuousGalerkinapproximations2000,ArnoldUnifiedAnalysisDiscontinuous2002}. Huynh~\cite{HuynhReconstructionApproachHighOrder2009} has reported that this formulation results in a singular matrix for steady problems.  A possible remedy for this problem is to add a penalty term to the \bff{BR1} method as in~\secref{sec:BR1_stabilized}. By increasing $\eta^{br1}$, the relative energy distribution between eigenmodes changes and the more dissipative modes (mode(2), mode(3)~\hfigref{fig:BR1p2eta0_wd}) at $kh=\pi$ start to have higher energies and hence enhances the performance of the scheme. Through~\combined analysis (shown next) of the \bff{BR1} approach, we are able to see that it has a very small decaying rate at the Nyquist frequency $kh=(p+1)\pi$ and can actually be nearly non-decaying for very long time in the simulation. This behavior is somehow similar to the DG scheme with central fluxes for linear advection~\cite{AlhawwaryFourierAnalysisEvaluation2018}. %
 Moreover,~\hfigref{fig:sdisc_physical_mode_allschemes} shows that for $\pi \leq kh \lesssim 2\pi$, \bff{LDG}p$2$-$\eta0$ has a dominant mode that better approximates the exact dissipation relation in this range than \bff{BR2}p$2$-$\eta1$. Additionally, the eigenmode of \bff{LDG}p$2$-$\eta0$ that has the highest energy for $kh \gtrsim 2\pi$ is more dissipative than that of the \bff{BR2}p$2$-$\eta1$. These observation suggest that \bff{LDG} schemes could be more accurate and robust for moderate to high wavenumbers and this is confirmed through the \ita{combined-mode} analysis in the subsequent sections. %
\begin{table}[H]%
\caption{ Minimum penalty $\eta$ for a stable DG scheme with \bff{BR2}. } 
\centering 
\small
\begin{tabular}{c|c|c|c|c|c|c|c|c}
\hline
 $p$ & $1$ & $2$ & $3$ & $4$ & $5$ & $6$ & $7$ & $8$ \\
 \hline
 $\eta$ & $0.50$ & $0.67$ & $0.75$ & $0.80$ & $0.83$ & $0.86$ & $0.88$ & $0.89$ \\
\hline
\end{tabular}
\label{table:minEta_stable_sdiscBR2}
\end{table}%
\begin{table}[H]%
\caption{ Minimum penalty $\eta$ for a stable \bff{LDG} scheme. } 
\centering 
\small
\begin{tabular}{c|c|c|c|c|c|c|c|c}
\hline
 $p$ & $1$ & $2$ & $3$ & $4$ & $5$ & $6$ & $7$ & $8$ \\
 \hline
 $\eta$ & $-3$ & $-5$ & $-7$ & $-9$ & $-11$ & $-13$ & $-15$ & $-17$ \\
\hline
\end{tabular}
\label{table:minEta_stable_sdiscLDG}
\end{table}%
In order to to study the von Neumann stability of a certain scheme, the semi-discrete $k_{m}$ of all eigenmodes of a DG diffusion scheme can be utilized. For a non-growing eigenmode solution, $-(k_{m} h)^{2} \leq 0$ and for strictly decaying mode $-(k_{m}h)^{2} < 0$. In studying the stability of a certain scheme we assume a prescribed wavenumber range $ kh \in [0, (p+1) \pi]$ and for a certain value of $\eta $ we check the sign of $-(k_{m} h)^{2}$ for all eigenmodes. If the sign of $-(k_{m} h)^{2}$ is positive this means the scheme has a growing mode and hence is unstable. Since the penalty term using lifting coefficients is necessary for the stability of a \bff{BR2} scheme~\cite{ArnoldUnifiedAnalysisDiscontinuous2002} we choose a range of $\eta^{br2}  \in (0,3] $ while for \bff{LDG} schemes we tried both negative and positive values for $\eta^{ldg}$. It turned out that the \bff{LDG} scheme is indeed stable for negative $\eta$ values unlike the \bff{BR2} which has a minimum positive value.

\htablsref{table:minEta_stable_sdiscBR2}{table:minEta_stable_sdiscLDG} presents the stability analysis results for \bff{BR2} and \bff{LDG}, respectively. From these tables, we can conclude that the minimum $\eta$ for the stabilization of \bff{BR2} schemes has the following relation%
\begin{subequations}\label{eqn:minEta_LDGBR2}
\begin{equation}
\min\left(\eta^{br2} \right) = \frac{p}{p+1},  \quad \lim_{p\rightarrow\infty} \min\left(\eta^{br2} \right) =1.0 , \quad p\geq 1,
\label{eqn:minEta_BR2}
\end{equation}%
while for the \bff{LDG}%
\begin{equation}
\min\left( \eta^{ldg} \right) = - \left( 2p+1 \right),  \quad p\geq 1.
\label{eqn:minEta_LDG}
\end{equation}
\end{subequations}%
On the other hand, the stabilized \bff{BR1} schemes has a minimum $\eta^{br1}=0$ for all orders. The optimal choice for $\eta$ obviously depends on a number of factors concerning accuracy, stability, and robustness of a diffusion scheme. This is investigated in more detail in the rest of this paper. %
%
\subsection{Combined-mode analysis}\label{sec:sdisc_combinedMode}
%
The~\combined analysis was first introduced in~\cite{AlhawwaryFourierAnalysisEvaluation2018} in order to study the \mytrue behavior of DG and RKDG schemes. In this approach all eigenmodes are considered in computing the response of the semi-discrete system for an initial wave mode. We extend the~\combined analysis to DG schemes for diffusion using the linear heat equation. In this regard, it is convenient to define two non-dimensional time scales as%
\begin{alignat}{2}
\tau  = \frac{\gamma t}{h^{2}} , \qquad \qquad \tau_{p} = (p+1)^{2} \tau,
\end{alignat}%
where $\tau, \tau_{p}$ are the non-dimensional diffusion time scales based on the cell width $h$ and the smallest length-scale of a DG scheme $h/(p+1)$, respectively. The time scale $\tau$ represents the well-known diffusion time-scale for single degree of freedom methods through an area $h^{2}$. The second time-scale $\tau_{p}$ represents the time required for the diffusion of a wave per DOF of a DG scheme through an area $(h/p+1))^{2}$. 

The energy based on the $L_{2}$ norm of a complex functions at any time $\tau_{p}$ is given by%
\begin{alignat}{2}
\begin{split}
E^{e}(k,\tau_{p}) & = \sqrt{ \frac{\int\limits_{-1}^{1} \bigl|u^{e}(\xi,\tau_{p}) \bigr|^{2} \, d\xi } {\int\limits_{-1}^{1} d\xi} } = \sqrt{\frac{1}{2}  \int\limits_{-1}^{1} \bigl| u^{e}(\xi,\tau_{p}) \bigr|^{2} \,  d\xi }  = \sqrt{\frac{1}{2}  \int\limits_{-1}^{1} \Biggl|  \sum\limits_{l=0}^{p}  U_{l} \phi_{l}(\xi)    \Biggr|^{2}  d\xi} \\ & = \sqrt{\frac{1}{2}  \int\limits_{-1}^{1} \Biggl| \sum\limits_{l=0}^{p} \sum\limits_{j=0}^{p}  \vartheta_{j} \mathcal{M}_{l,j} \: e^{- K^{2}_{m,j}  \tau_{p} } \:  \phi_{l}(\xi) \Biggr|^{2} \,  d\xi}  ,
\end{split}
\label{eqn:sdisc_E}
\end{alignat}%
while the exact energy is defined by%
\begin{alignat}{2}
E_{ex}^{e}(k,\tau_{p}) & = \sqrt{\frac{1}{2}  \int\limits_{-1}^{1} \bigl| u_{ex}^{e}(\xi,\tau_{p}) \bigr|^{2} \,  d\xi } = \sqrt{\frac{1}{2}  \int\limits_{-1}^{1} \Biggl| \sum\limits_{l=0}^{p} \hat{\mu}_{l} \: e^{- K^{2} \tau_{p} } \: \phi_{l}(\xi) \Biggr|^{2} \,  d\xi } ,
\label{eqn:sdisc_Eex}
\end{alignat}
where we have eliminated $e^{ikx_{e}}$ since $|e^{ikx_{e}}|=1$. This definition of the energy is equivalent to an averaged amplitude over $\Omgs$. Hence, the~\mytrue diffusion factor can be written as%
\begin{alignat}{2}
G_{true}(k,\tau_{p})  = \frac{E^{e}(k,\tau_{p})}{E^{e}_{ex}(k,0)}
\label{eqn:sdisc_G}
\end{alignat}%
while the \physic diffusion factor is given by%
\begin{equation}
G_{phys}(k,\tau_{p}) = e^{- \tilde{\omega} \, t} = e^{-K_{m}^{2}  \, \tau_{p}}.
\label{eqn:sdisc_Gphys}
\end{equation}%
In addition, an absolute diffusion error can be defined as%
\begin{equation}
\bigl| \Delta G(\tau_{p})  \bigr| =  \bigl| G_{exact} - G_{numerical} \bigr|  = \Biggl| \frac{E_{ex}(k,\tau_{p}) - E(k,\tau_{p})}{E_{ex}(k,0)} \Biggr|
\label{eqn:true_Gerr}
\end{equation}%
\begin{figure}[H]
\centering
    \includegraphics[width=0.65\textwidth]{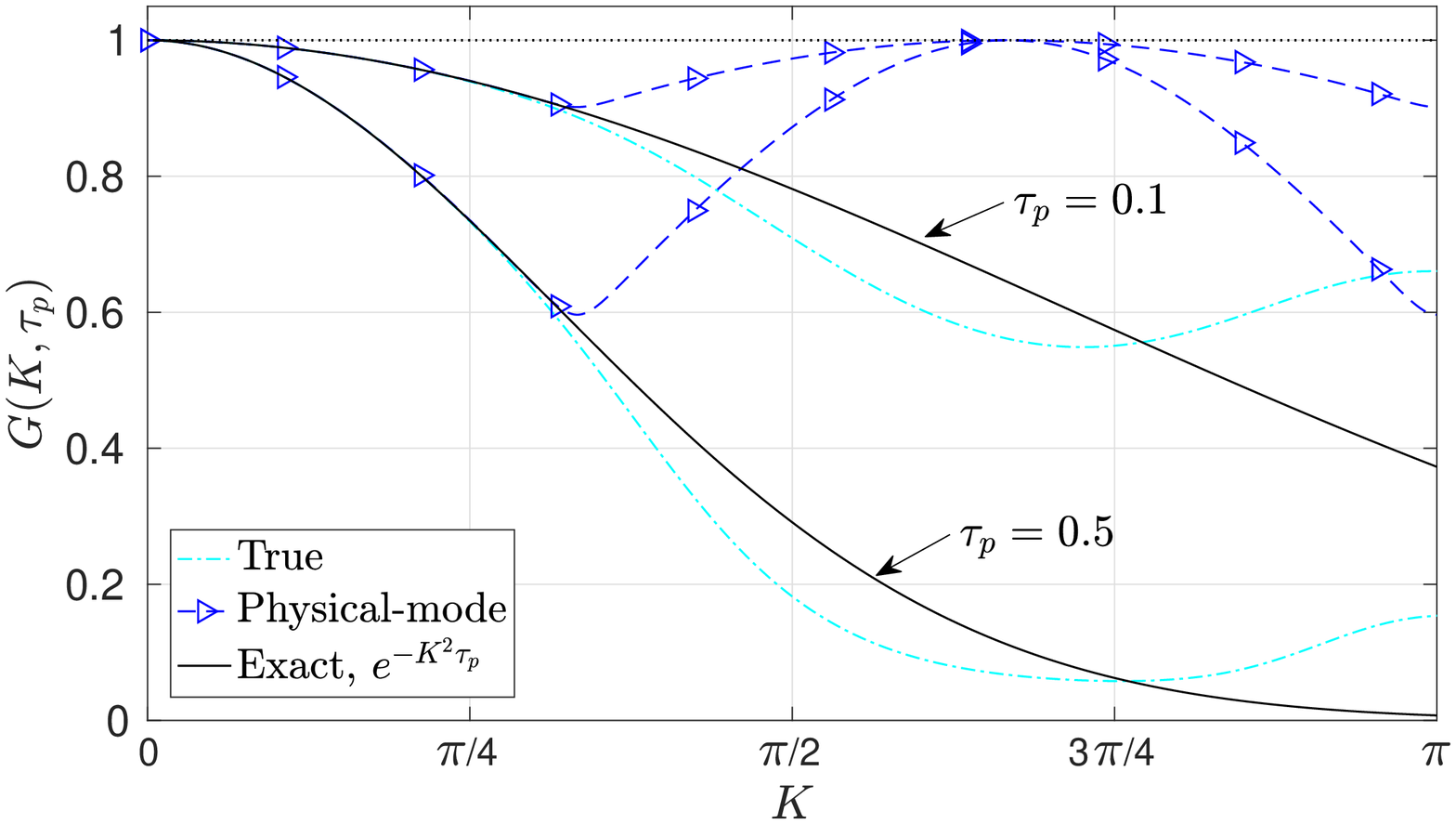} 
	\caption{Semi-discrete~\combined analysis of the BR$2$p$2$-$\eta1$. Comparison of~\mytrue and~\physic behaviors.}
    \label{fig:BR2p2eta1_G}
\end{figure}%
\hfigref{fig:BR2p2eta1_G} displays the results of the semi-discrete~\combined analysis in comparison to the standard Fourier analysis based on the~\physic. From this figure it can be seen clearly that the~\physic approximates the exact dissipation and is close to the~\mytrue behavior for only a small range of wavenumbers $K \lesssim \pi/3$. In addition, the~\mytrue behavior of DG schemes for diffusion is reasonably accurate in comparison with the exact one for a wide range of wavenumbers as expected. It has also less dissipation for high wavenumbers than the exact dissipation, which while still stable it could trigger some aliasing and stability issues for nonlinear problems. It is always desirable to have more dissipation for high wavenumbers than the exact dissipation in under-resolved large eddy simulations (LES). This indicates the effectiveness of the~\combined analysis for studying the~\mytrue behavior of DG schemes for diffusion and assessing their robustness and stability. %
%
\subsubsection{Diffusion characteristics, and the effects of the penalty parameter and polynomial orders}
\label{sec:sdisc_accuracy_eta}
%
In this section, we distinguish between short and long time diffusion behaviors. The classical single eigenmode analysis suggests that the individual eigenvalue dissipation characteristics with respect to the exact dissipation remains the same for long time, i.e., long time is just an accumulation of more errors.  In contrast, through \combined analysis one can see that the combination of all modes results in a very different behavior for short and long times.

\begin{figure}[H]
\begin{subfigure}[h]{0.5\textwidth}
    \includegraphics[width=0.965\textwidth]{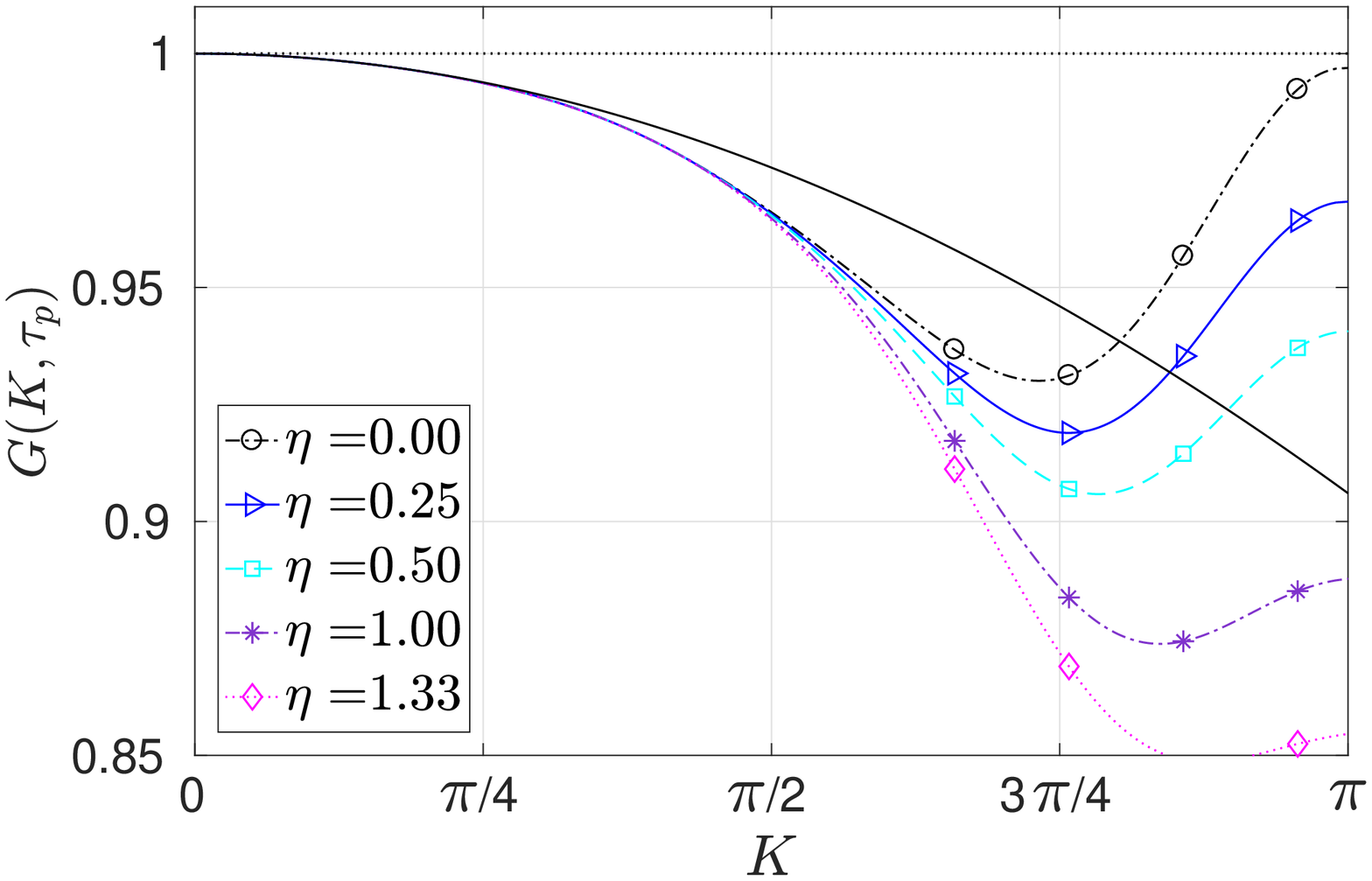} 
    \caption{Diffusion factor for BR$1$p$2$} 
    \label{fig:BR1p2_compareEta_sdisc_G_shortTime}
    \end{subfigure}
    \, \,
   \begin{subfigure}[h]{0.5\textwidth}
    \includegraphics[width=0.965\textwidth]{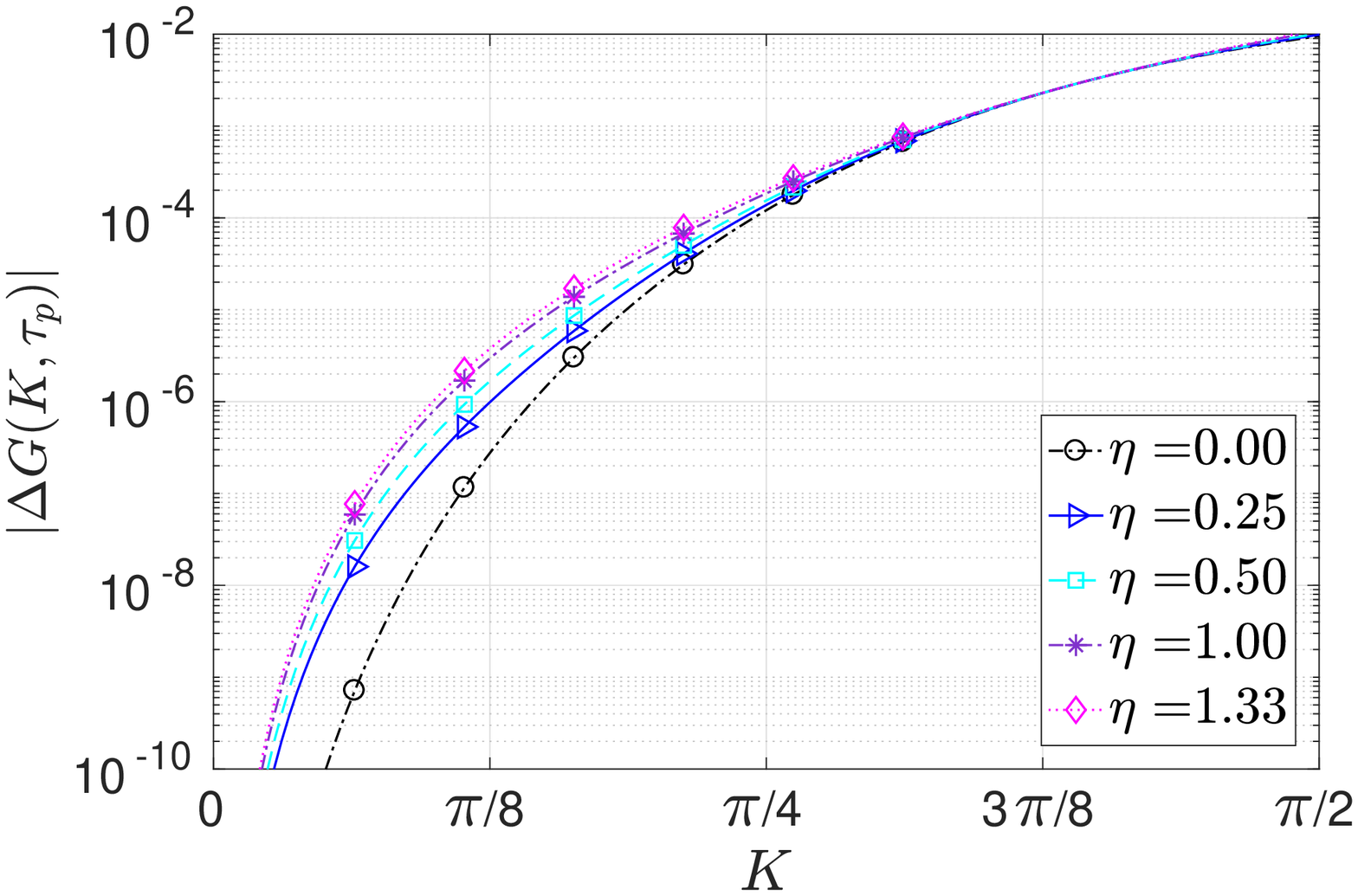}  
    \caption{Diffusion error for BR$1$p$2$} 
    \label{fig:BR1p2_compareEta_sdisc_Gerr_shortTime}
    \end{subfigure}   \\ \\ \\
    \begin{subfigure}[h]{0.5\textwidth}
    \includegraphics[width=0.965\textwidth]{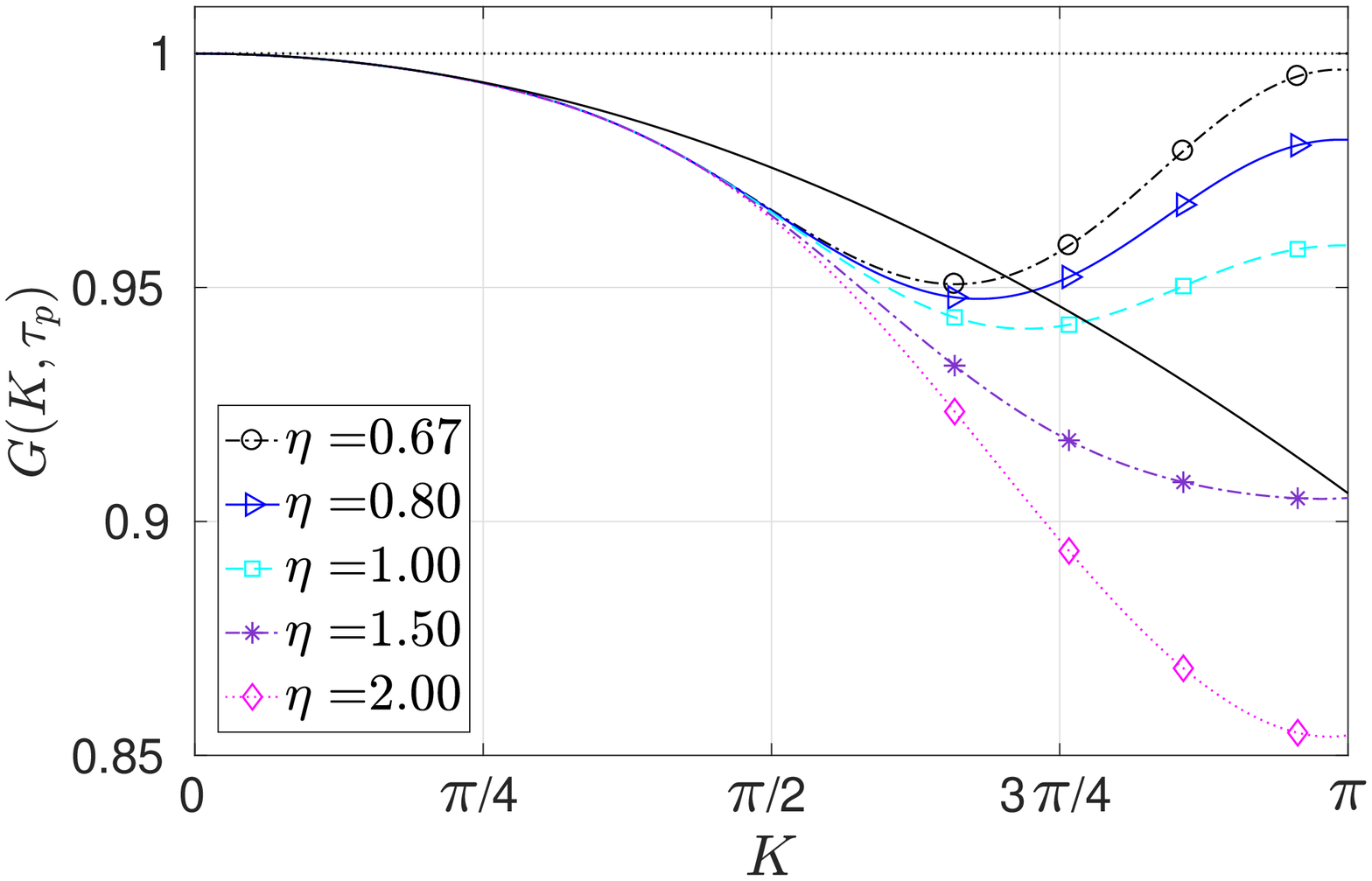} 
    \caption{Diffusion factor for BR$2$p$2$} 
    \label{fig:BR2p2_compareEta_sdisc_G_shortTime}
    \end{subfigure}
    \, \,
   \begin{subfigure}[h]{0.5\textwidth}
    \includegraphics[width=0.965\textwidth]{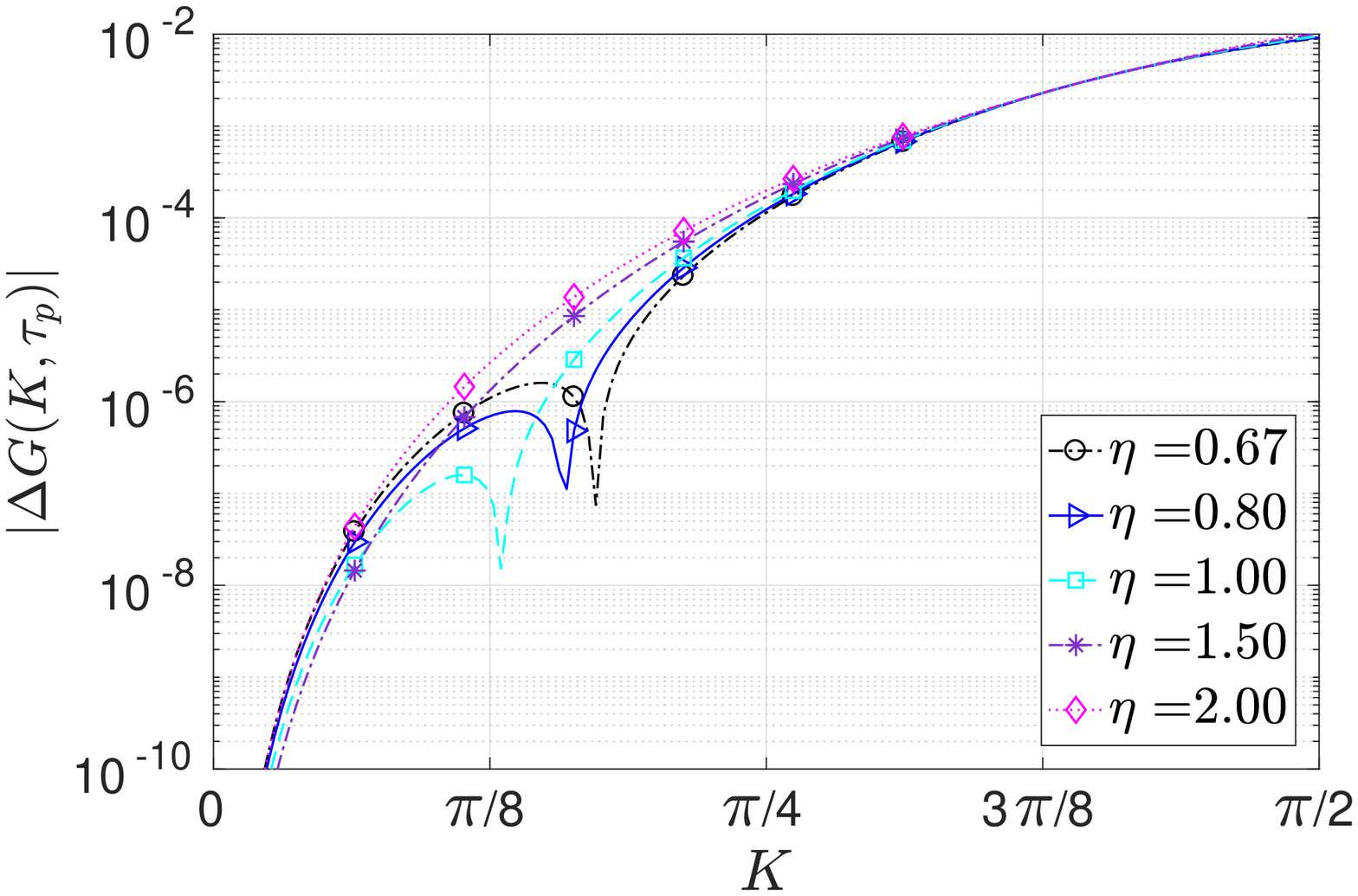}  
    \caption{Diffusion error for BR$2$p$2$} 
    \label{fig:BR2p2_compareEta_sdisc_Gerr_shortTime}
    \end{subfigure} \\ \\ \\
    \begin{subfigure}[h]{0.5\textwidth}
    \includegraphics[width=0.965\textwidth]{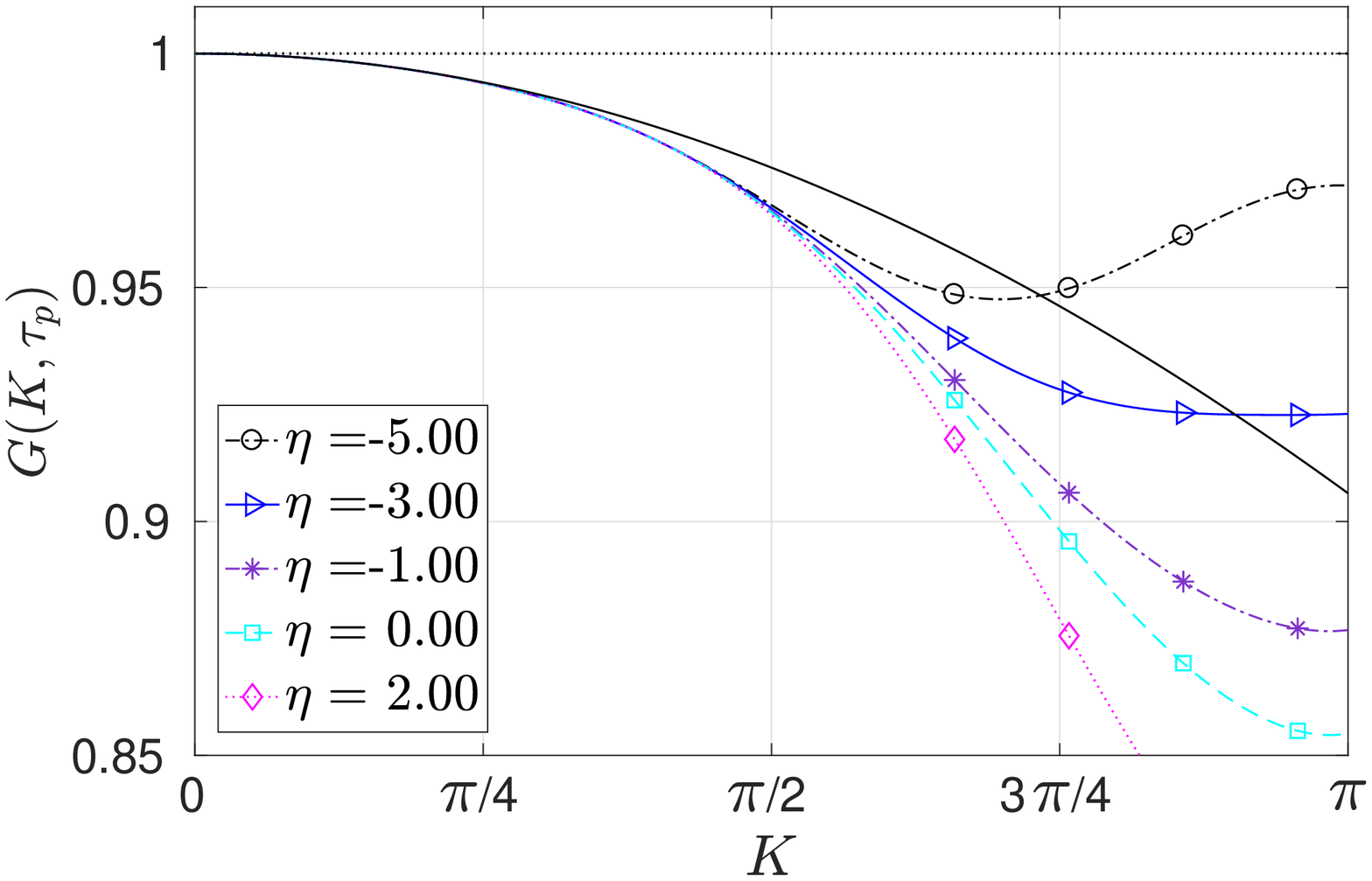} 
    \caption{Diffusion factor for LDGp$2$} 
    \label{fig:LDGp2_compareEta_sdisc_G_shortTime}
    \end{subfigure}
    \, \,
   \begin{subfigure}[h]{0.5\textwidth}
    \includegraphics[width=0.965\textwidth]{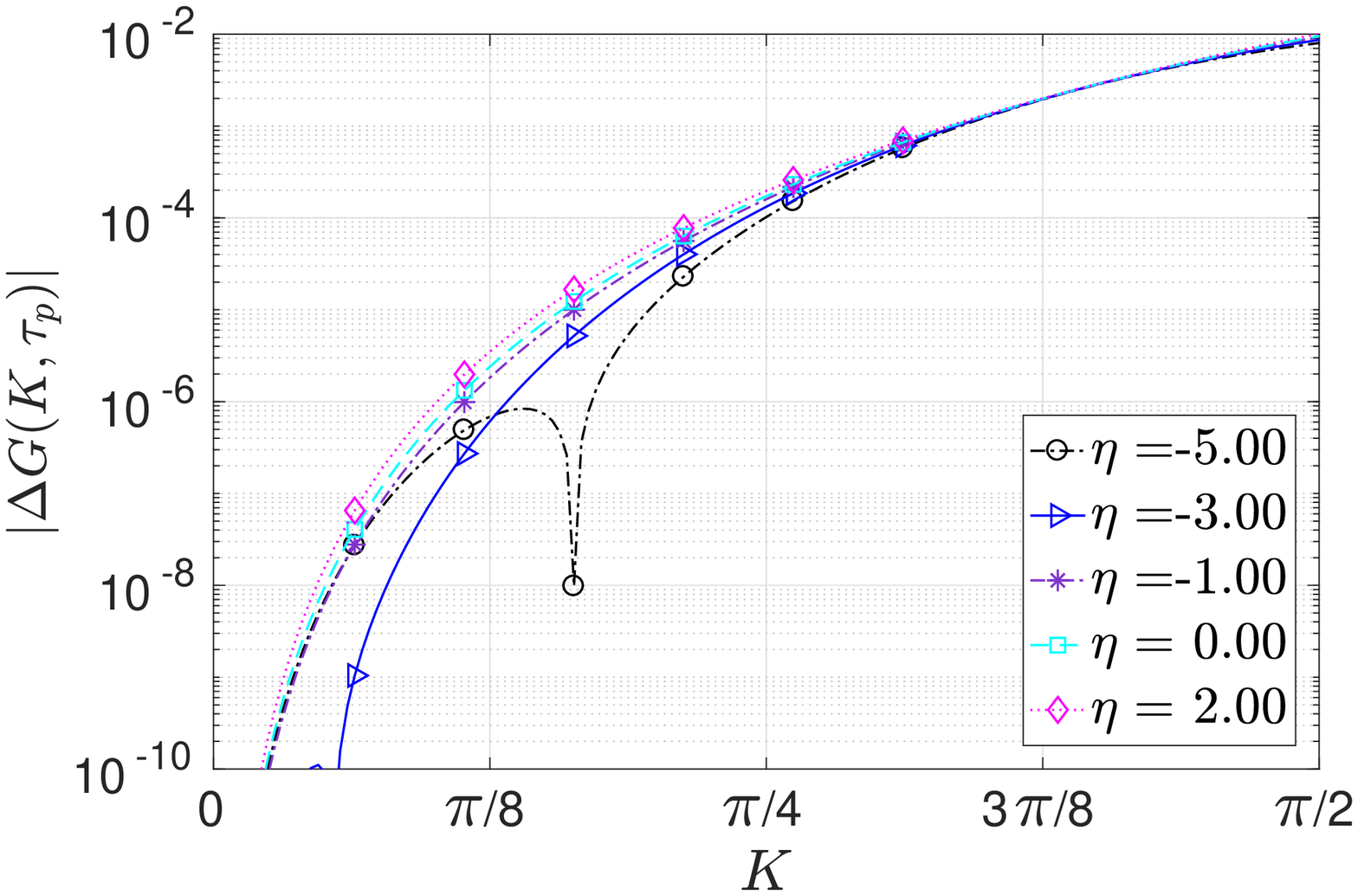}  
    \caption{Diffusion error for LDGp$2$} 
    \label{fig:LDGp2_compareEta_sdisc_Gerr_shortTime}
    \end{subfigure}
	\caption{Effect of the penalty parameter on the semi-discrete~\mytrue behavior of diffusion schemes for a short time $\tau_{p}=0.01$. In the left column of figures, the solid line without symbols represents the exact diffusion factor $G_{ex}=e^{-K^{2} \tau_{p}}$.}
    \label{fig:compareEta_sdisc_shortTime}
\end{figure}

\begin{figure}[H]
\begin{subfigure}[h]{0.5\textwidth}
    \includegraphics[width=0.965\textwidth]{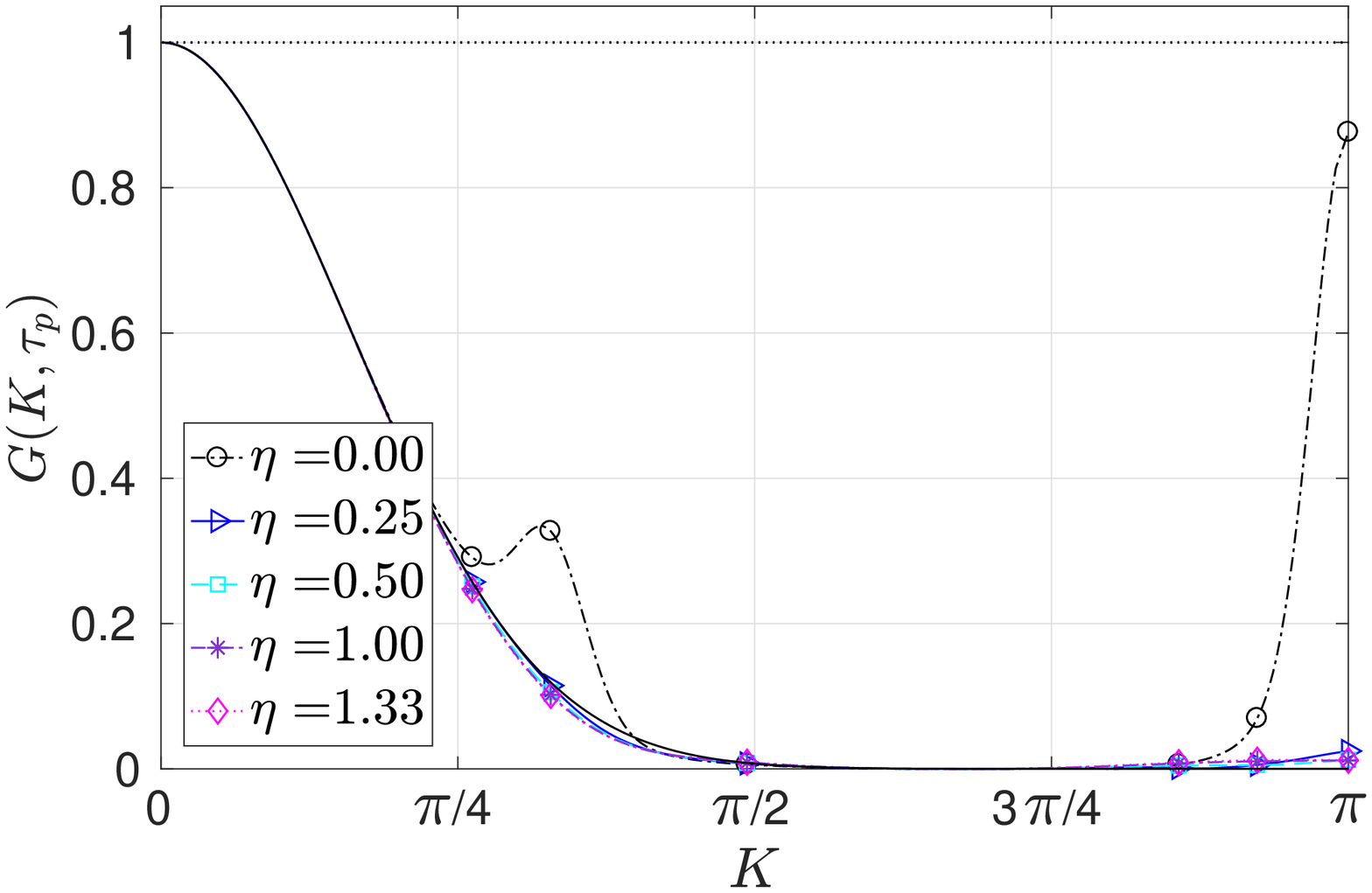} 
    \caption{Diffusion factor for BR$1$p$2$} 
    \label{fig:BR1p2_compareEta_sdisc_G_longTime}
    \end{subfigure}
    \, \,
   \begin{subfigure}[H]{0.5\textwidth}
    \includegraphics[width=0.965\textwidth]{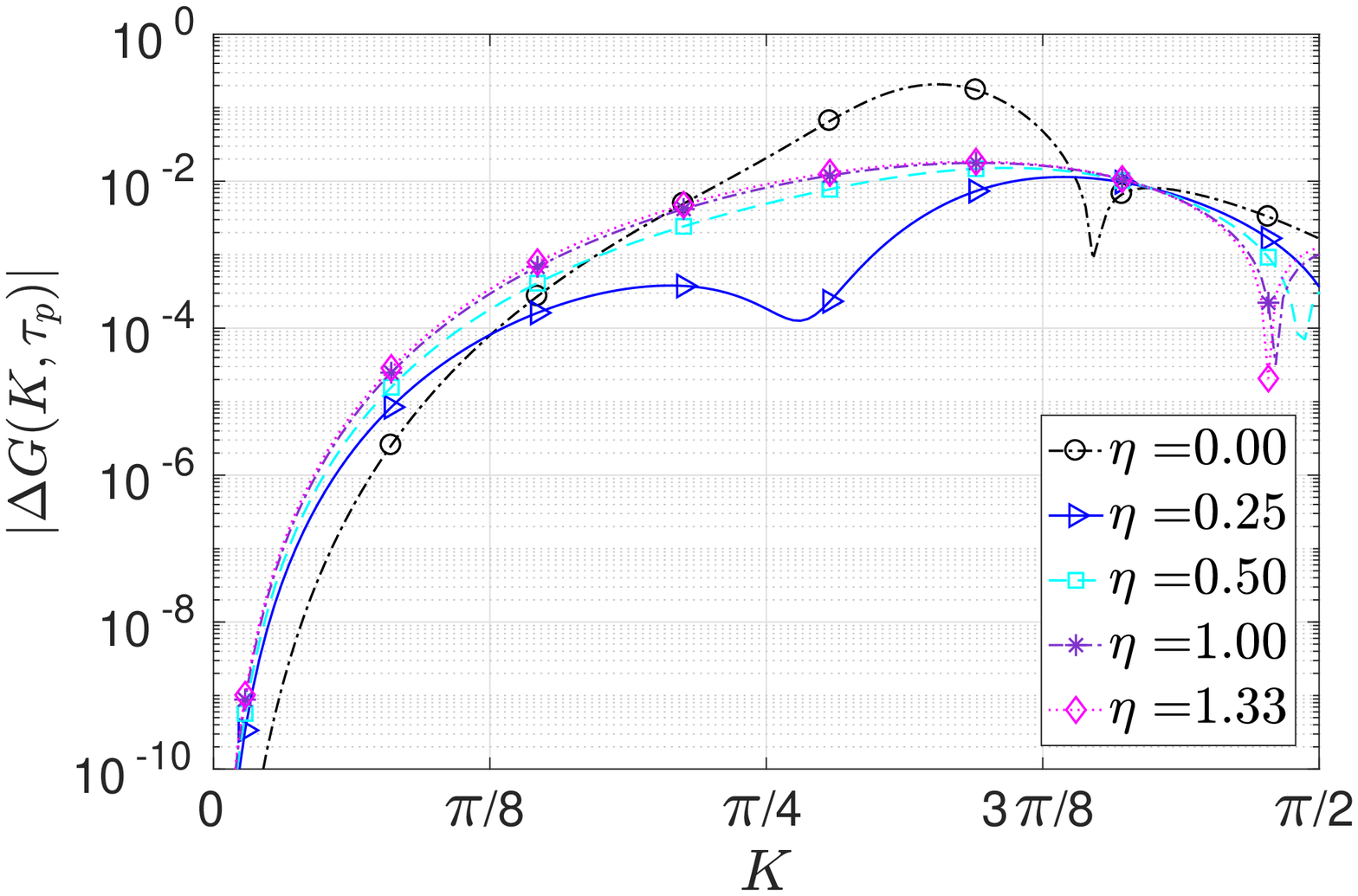}  
    \caption{Diffusion error for BR$1$p$2$} 
    \label{fig:BR1p2_compareEta_sdisc_Gerr_longTime}
    \end{subfigure}  \\ \\ \\
    \begin{subfigure}[h]{0.5\textwidth}
    \includegraphics[width=0.965\textwidth]{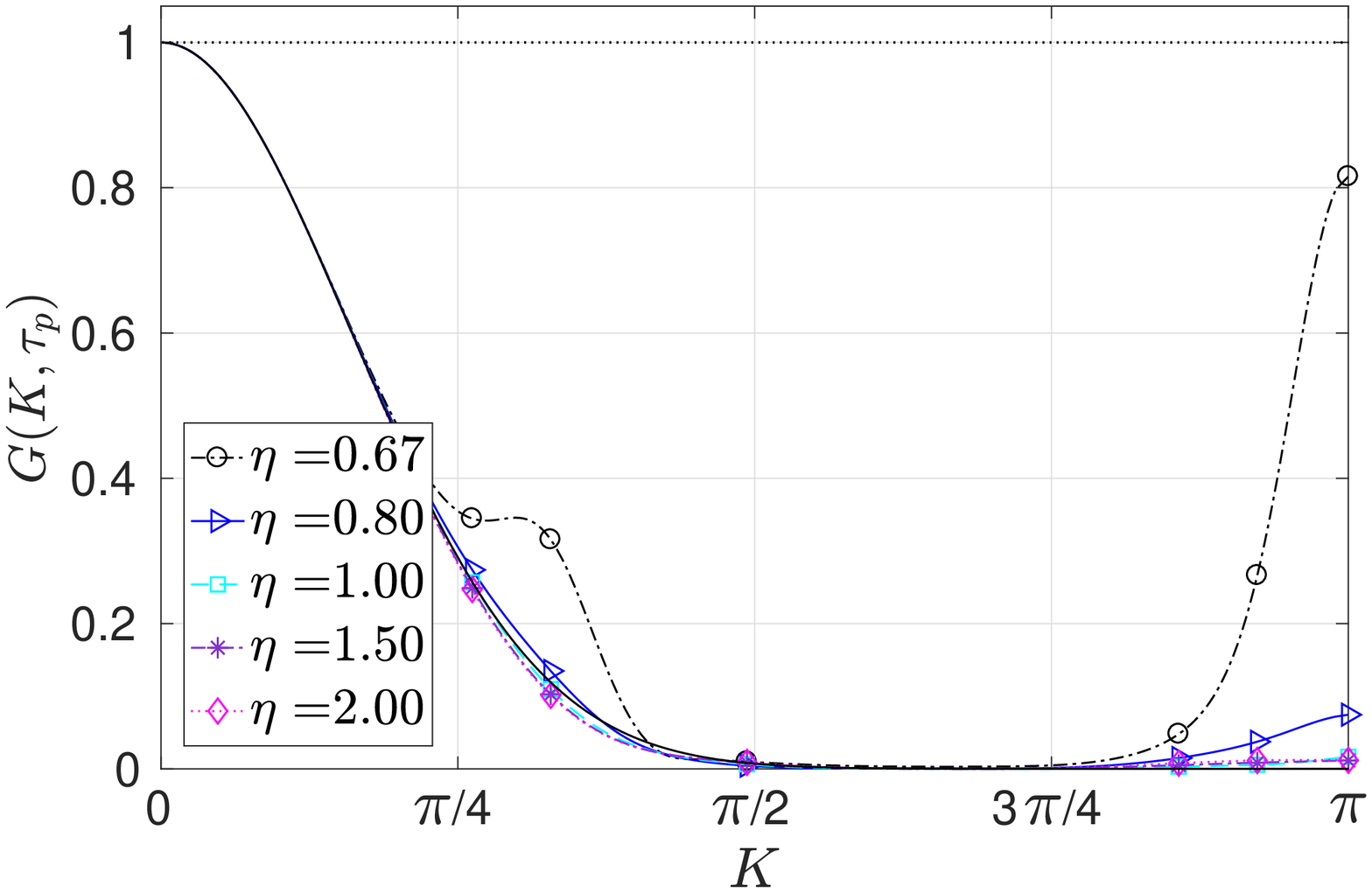} 
    \caption{Diffusion factor for BR$2$p$2$} 
    \label{fig:BR2p2_compareEta_sdisc_G_longTime}
    \end{subfigure}
    \, \,
   \begin{subfigure}[h]{0.5\textwidth}
    \includegraphics[width=0.965\textwidth]{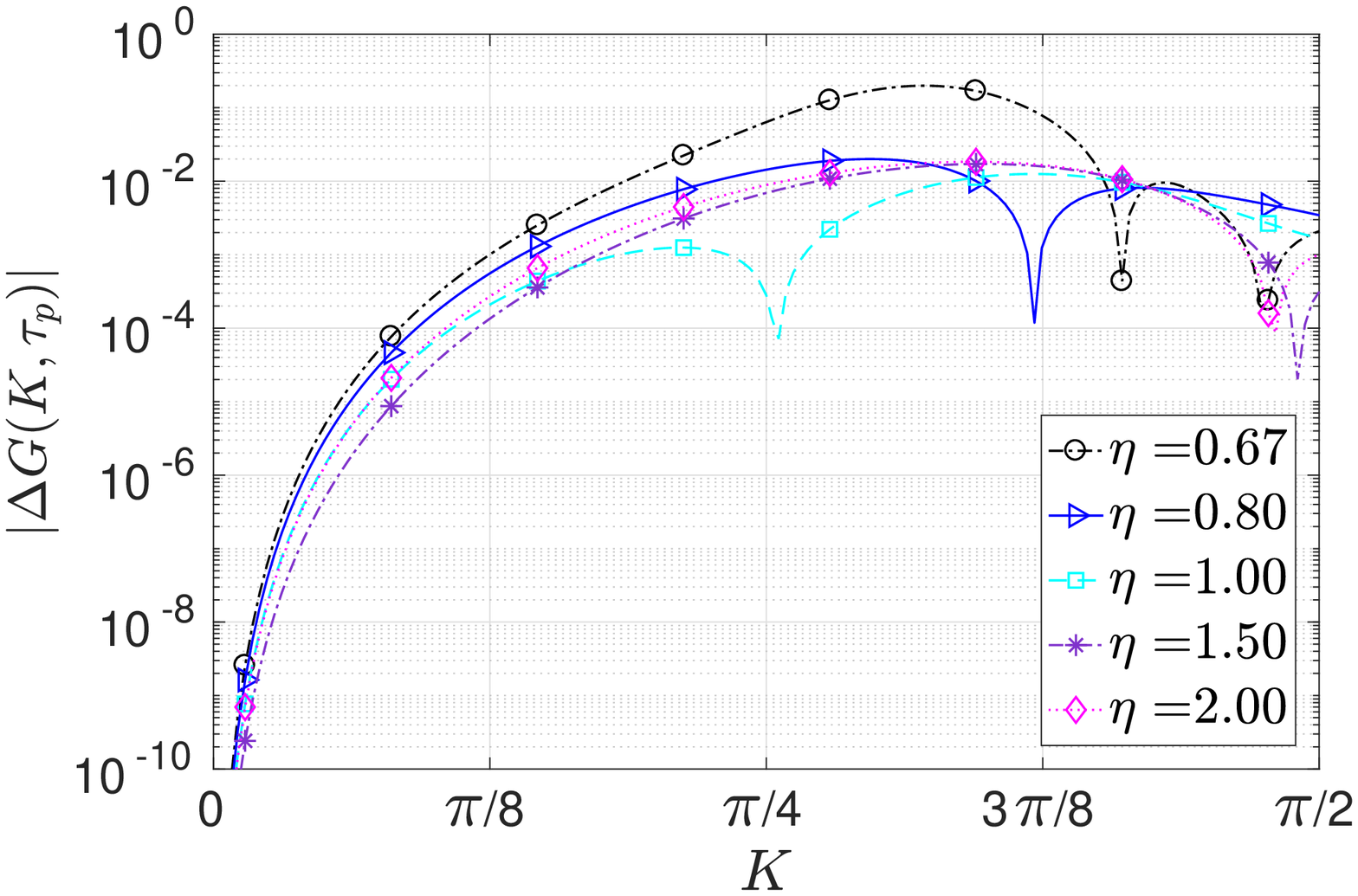}  
    \caption{Diffusion error for BR$2$p$2$} 
    \label{fig:BR2p2_compareEta_sdisc_Gerr_longTime}
    \end{subfigure} \\ \\ \\
    \begin{subfigure}[h]{0.5\textwidth}
    \includegraphics[width=0.965\textwidth]{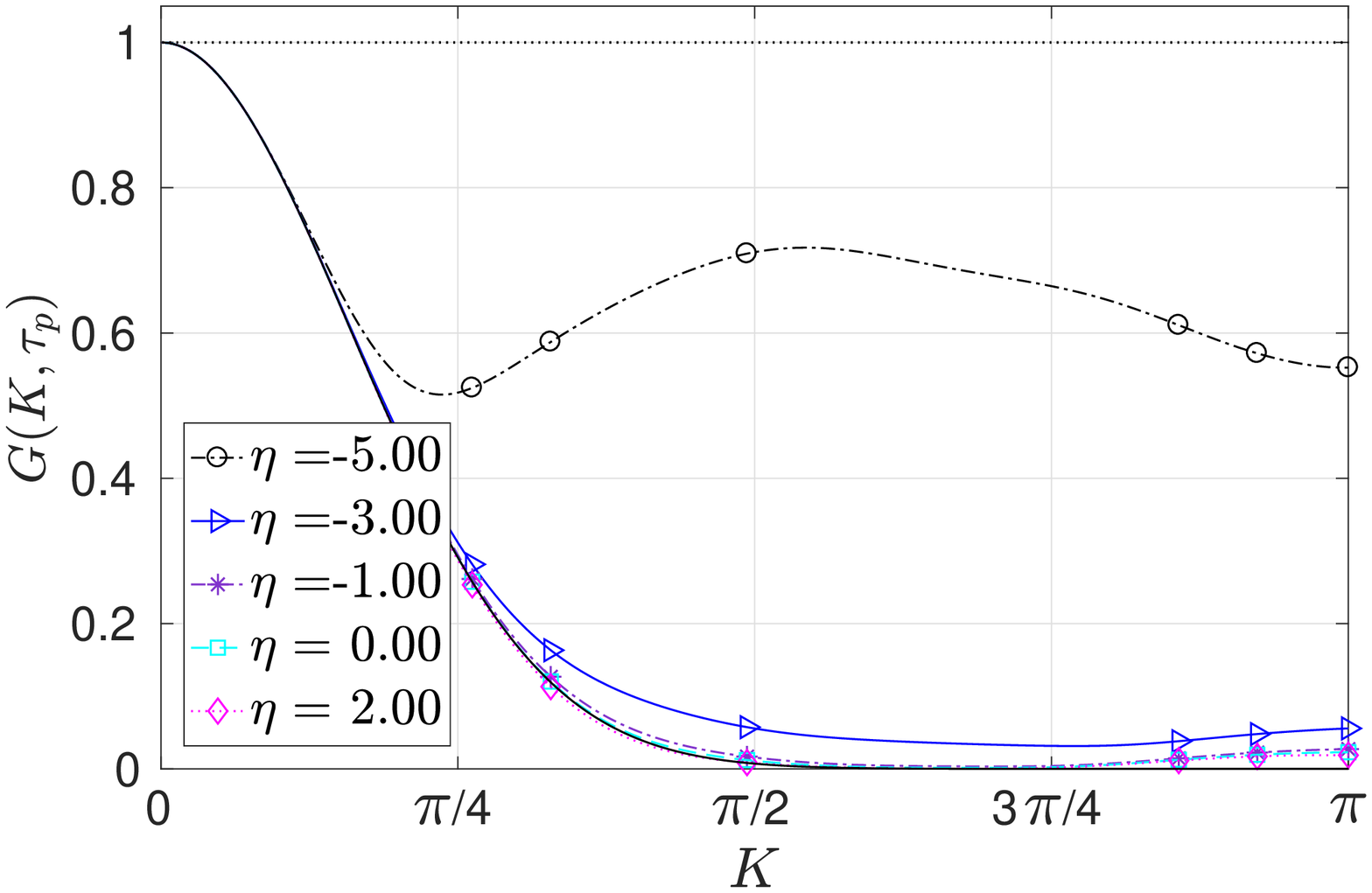} 
    \caption{Diffusion factor for LDGp$2$} 
    \label{fig:LDGp2_compareEta_sdisc_G_longTime}
    \end{subfigure}
    \, \,
   \begin{subfigure}[h]{0.5\textwidth}
    \includegraphics[width=0.965\textwidth]{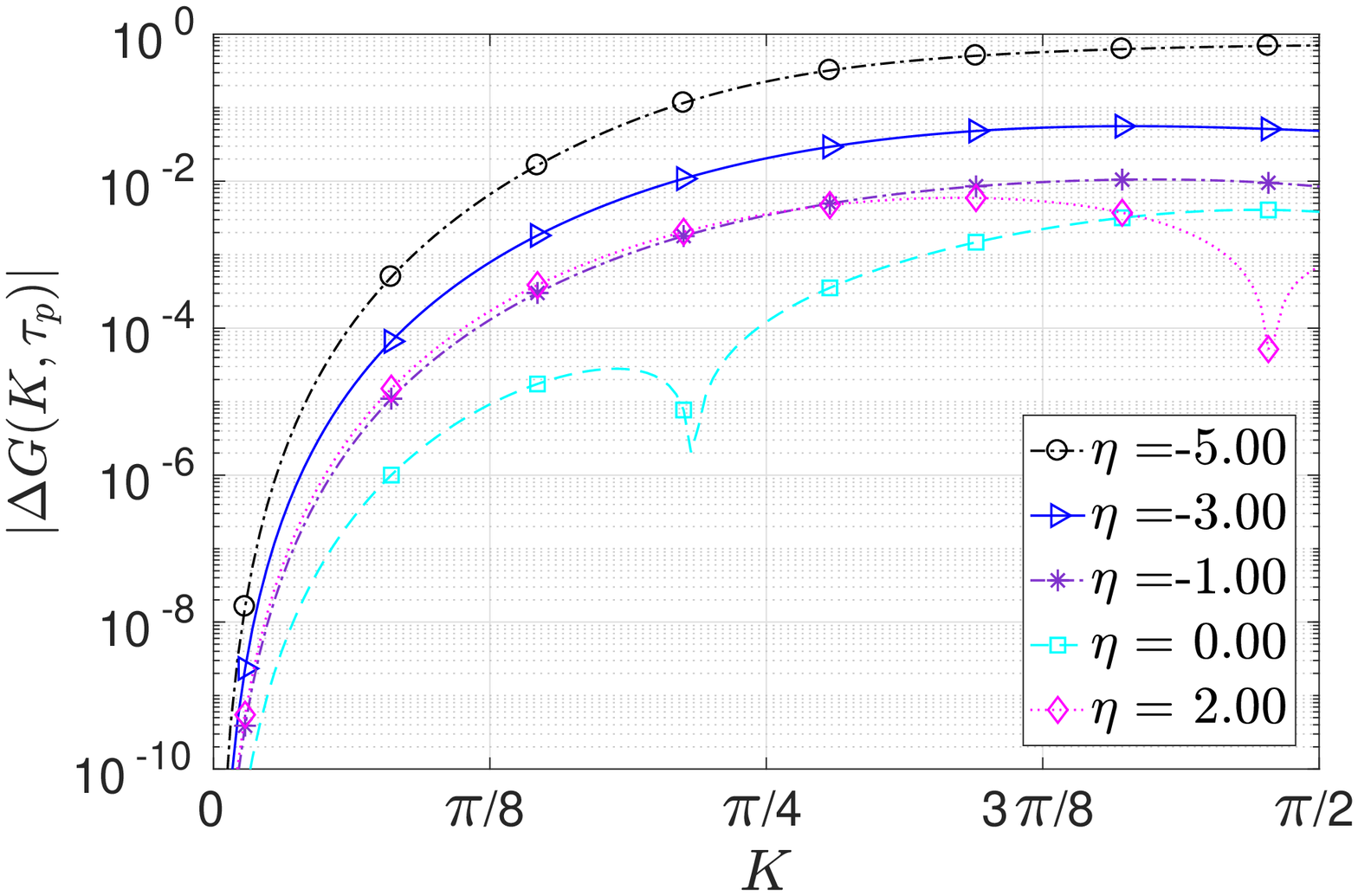}  
    \caption{Diffusion error for LDGp$2$} 
    \label{fig:LDGp2_compareEta_sdisc_Gerr_longTime}
    \end{subfigure}
	\caption{Effect of the penalty parameter on the semi-discrete~\mytrue behavior of diffusion schemes for long time $\tau_{p}=2.0$. In the left column of figures, the solid line without symbols represents the exact diffusion factor $G_{ex}=e^{-K^{2} \tau_{p}}$. }
    \label{fig:compareEta_sdisc_longTime}
\end{figure}

We start by studying the effects of the penalty parameter $\eta$ on the diffusion behavior of the three considered viscous flux formulations.~\hfigref{fig:compareEta_sdisc_shortTime} presents the short time diffusion factors and diffusion errors as a function of $\eta$ for p$2$ schemes with the three viscous flux formulations considered in this work. From this figure, as $\eta$ increases, all schemes have a general trend of increasing the dissipation for all wavenumbers. The long time diffusion behavior is presented in~\hfigref{fig:compareEta_sdisc_longTime} where we can see that for very low values of $\eta$ or minimum values for stability the scheme becomes non robust due to very low and less than exact dissipation for high wavenumbers ($K \gtrsim 3\pi/4$). %
%
For the \bff{BR2}p$2$ schemes, it can be seen from~\hfigsref{fig:BR2p2_compareEta_sdisc_G_shortTime}{fig:BR2p2_compareEta_sdisc_Gerr_shortTime} that an $\eta \geq 1.5$ provides larger than exact dissipation for $K \gtrsim 3\pi/4$ without sacrificing the low wavenumber accuracy and hence is a promising candidate for robust simulations. This can also be noticed for long time simulations in~\hfigsref{fig:BR2p2_compareEta_sdisc_G_longTime}{fig:BR2p2_compareEta_sdisc_Gerr_longTime}. In addition, for the \bff{BR1}p$2$-$\eta0$, the diffusion factor at the Nyquist frequency $K=\pi$ has a very slow decaying rate that almost saturates at some value for a long time of the simulation. This behavior is similar to the case of RKDG schemes for advection with central fluxes~\cite{AlhawwaryFourierAnalysisEvaluation2018}, where it was reported that the amplification factor of the scheme saturates at some level for a very large number of iterations at $K=\pi$ . It is mainly attributed to having a zero eigenvalue dissipation coefficient at $K=\pi$ for one of the eigenmodes as is shown in~\secref{sec:sdisc_eigenmode_analysis}. The \bff{LDG}p$2$-$\eta0$ scheme in its standard form is better than all the stabilized \bff{LDG}p$2$ schemes with different $\eta$ in terms of accuracy and robustness. For long time simulations, we can see in~\hfigref{fig:LDGp2_compareEta_sdisc_Gerr_longTime} that this scheme has the least diffusion errors for low wavenumbers ($K \leq 3 \pi/ 8$) among all other schemes studied in this work. Moreover, in~\hfigref{fig:LDGp2_compareEta_sdisc_G_shortTime} it is clear that, for short time behavior, this scheme is more dissipative than the exact solution for moderate to high wavenumbers $K > \pi/2$.

\begin{figure}[H]
 \begin{subfigure}[h]{0.5\textwidth}
    \includegraphics[width=0.975\textwidth]{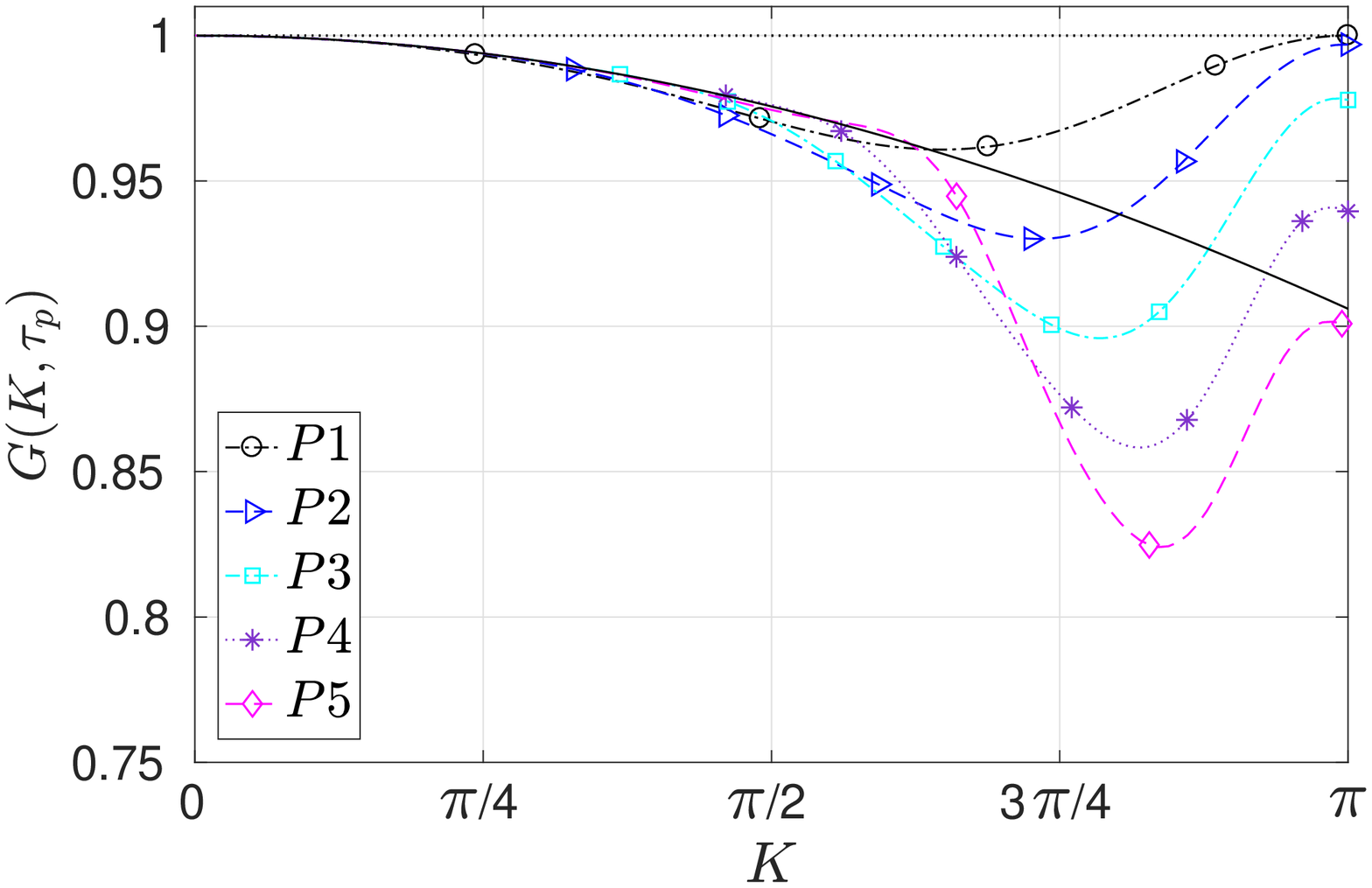} 
    \caption{Diffusion factor for \bff{BR1}-$\eta0.0$ at $\tau_{p}=0.01$} 
    \label{fig:BR1_compareP_G_tau0.01}
    \end{subfigure}
    \, \,
   \begin{subfigure}[h]{0.5\textwidth}
    \includegraphics[width=0.975\textwidth]{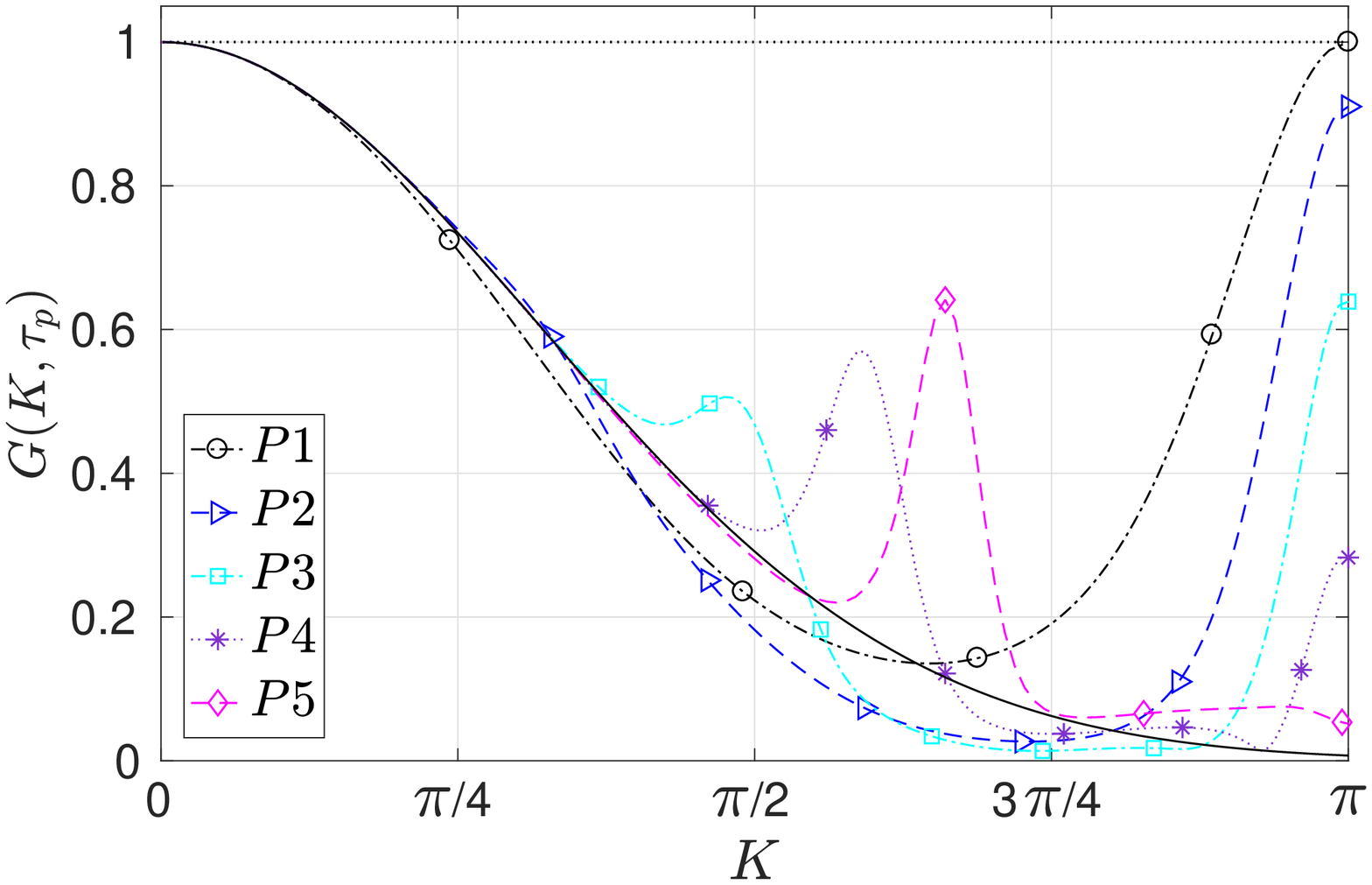}  
    \caption{Diffusion factor for \bff{BR1}-$\eta0.0$ at $\tau_{p}=0.5$} 
    \label{fig:BR1_compareP_G_tau0.5}
    \end{subfigure} \\ \\ \\
    \begin{subfigure}[h]{0.5\textwidth}
    \includegraphics[width=0.975\textwidth]{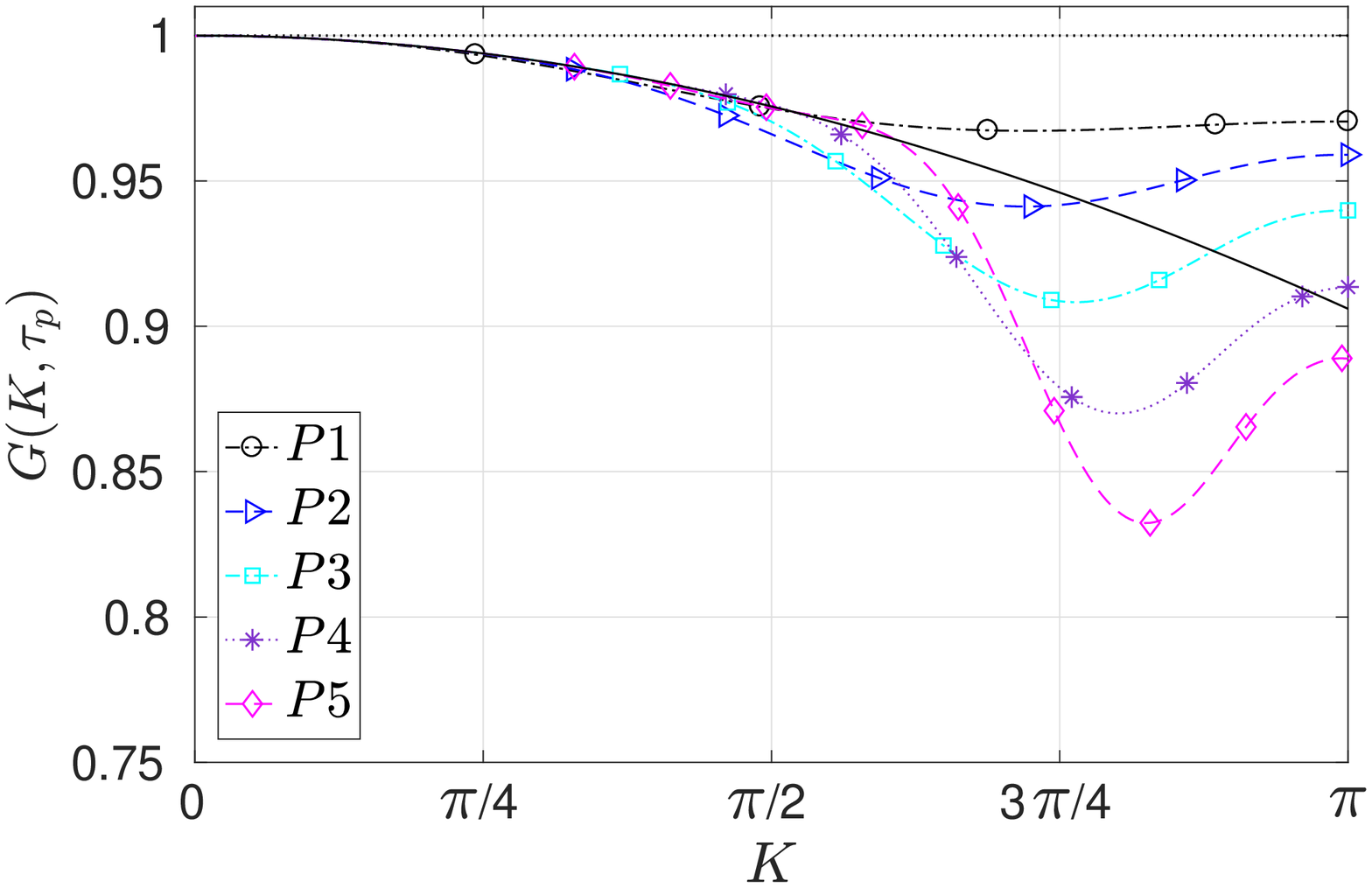} 
    \caption{Diffusion factor for \bff{BR2}-$\eta1.0$ at $\tau_{p}=0.01$} 
    \label{fig:BR2_compareP_G_tau0.01}
    \end{subfigure}
    \, \,
   \begin{subfigure}[h]{0.5\textwidth}
    \includegraphics[width=0.975\textwidth]{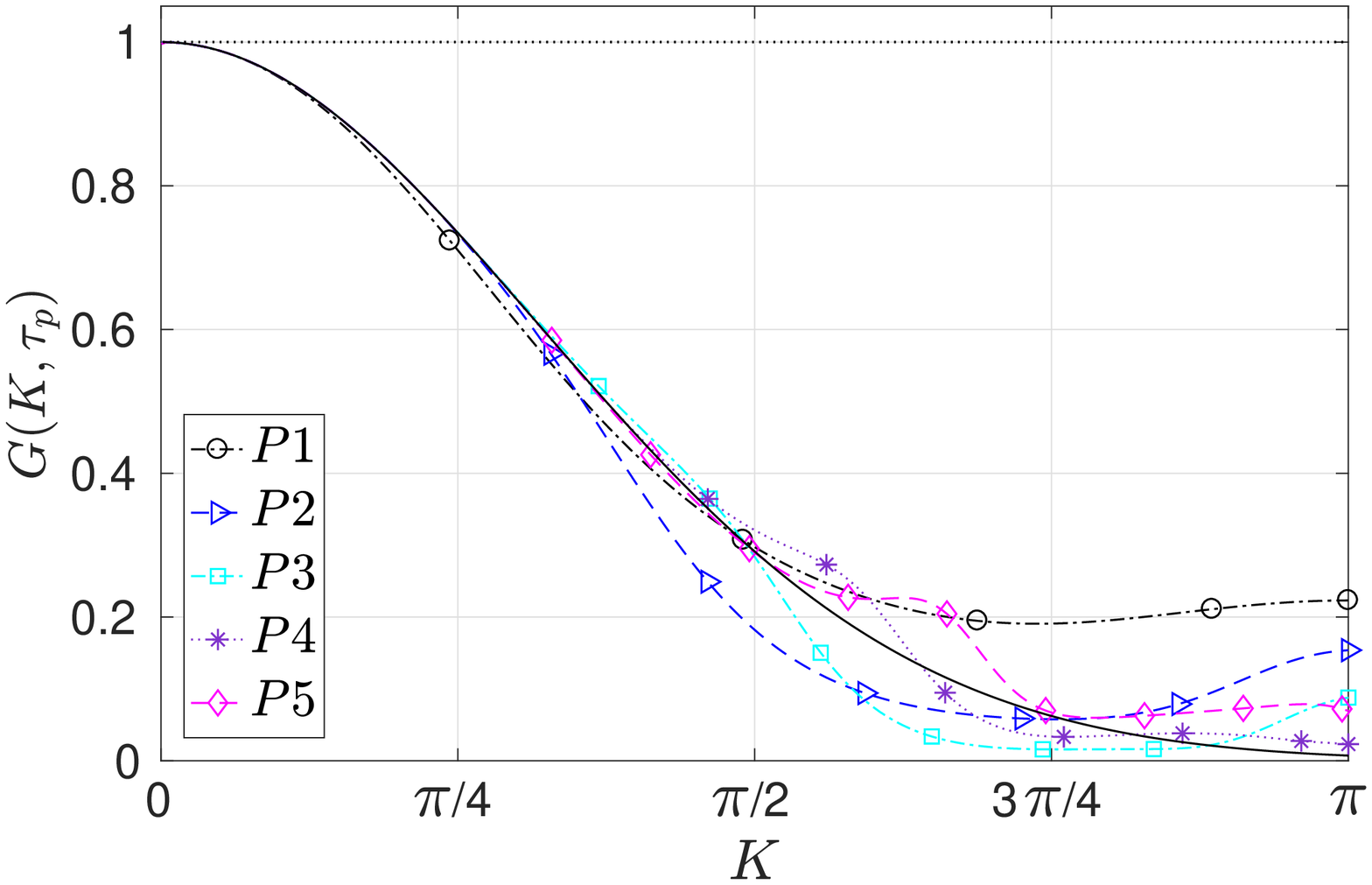}  
    \caption{Diffusion factor for \bff{BR2}-$\eta1.0$ at $\tau_{p}=0.5$} 
    \label{ffig:BR2_compareP_G_tau0.5}
    \end{subfigure} \\ \\ \\
    \begin{subfigure}[h]{0.5\textwidth}
    \includegraphics[width=0.975\textwidth]{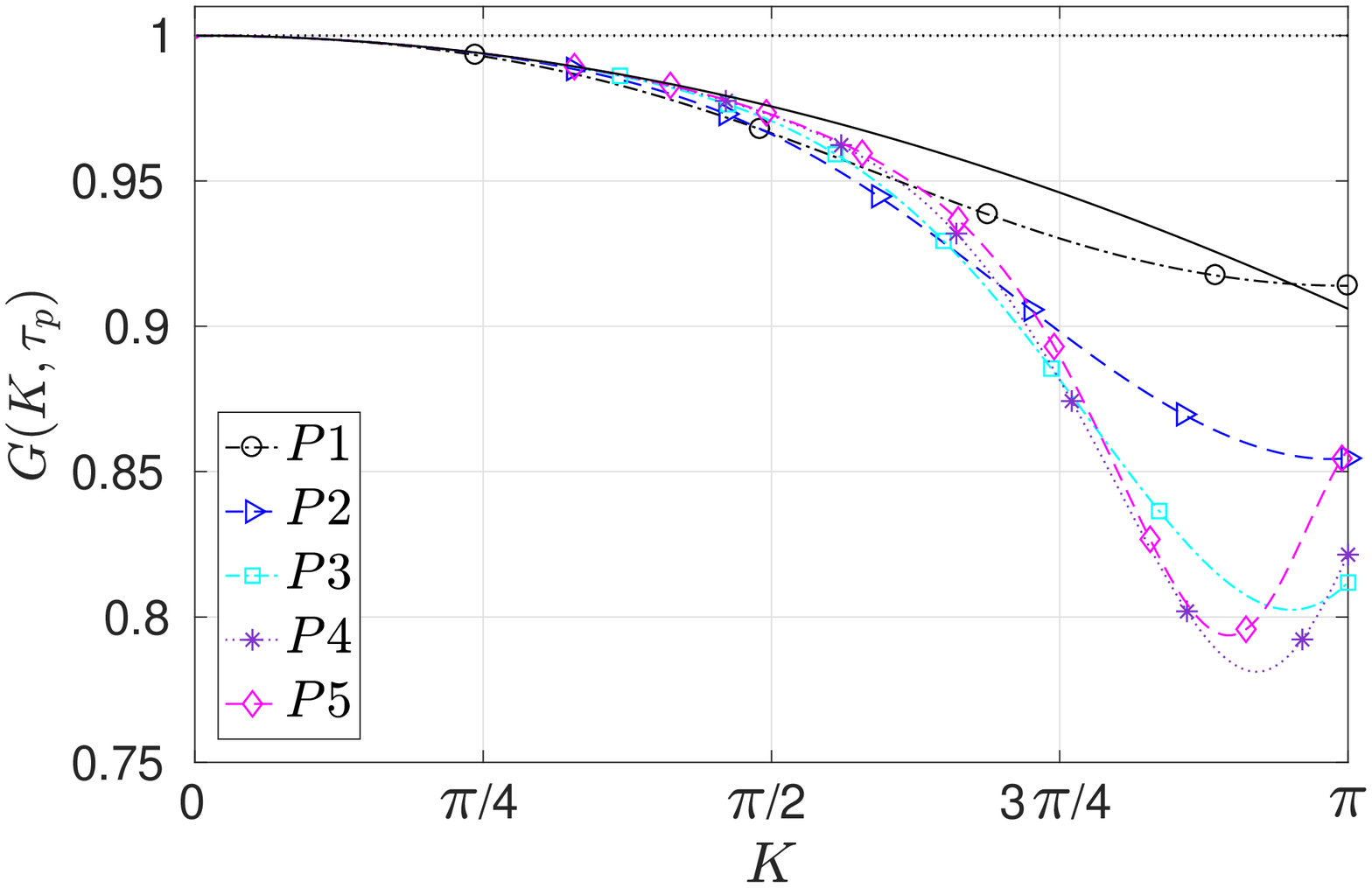} 
    \caption{Diffusion factor for \bff{LDG}-$\eta0.0$ at $\tau_{p}=0.01$} 
    \label{fig:LDG_compareP_G_tau0.01}
    \end{subfigure}
    \, \,
   \begin{subfigure}[h]{0.5\textwidth}
    \includegraphics[width=0.975\textwidth]{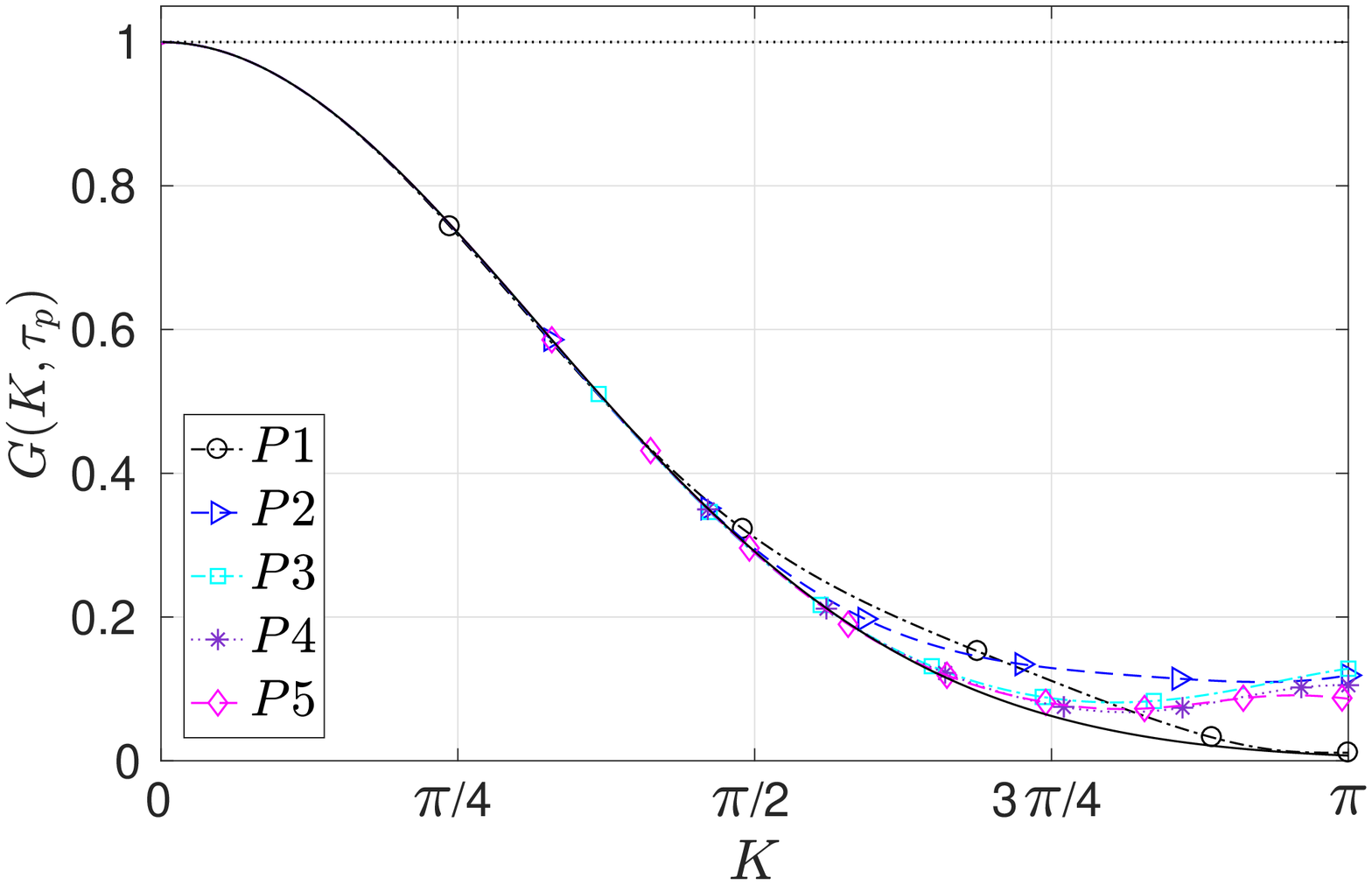}  
    \caption{Diffusion factor for \bff{LDG}-$\eta0.0$ at $\tau_{p}=0.5$}
    \label{fig:LDG_compareP_G_tau0.5}
    \end{subfigure} 
	\caption{Effect of the polynomial order on the semi-discrete~\mytrue diffusion factor.  In these figures, the solid line without symbols represents the exact diffusion factor $G_{ex}=e^{-K^{2} \tau_{p}}$.  }
    \label{fig:sdisc_compareP}
\end{figure}

Next, we show that the polynomial order has a significant effect on the \combined diffusion behavior of different schemes.~\hfigref{fig:sdisc_compareP} displays the diffusion factor of the standard \bff{BR1}-$\eta0$, \bff{BR2}-$\eta1$, and \bff{LDG}-$\eta0$ schemes for p$1$ to p$5$ polynomial orders. In general, as the order increases the scheme becomes more accurate in the low wavenumber range while more dissipative for high wavenumbers. That said, it is noticed that, as the order increases, both \bff{BR1}-$\eta0$ and \bff{BR2}-$\eta1$ schemes experience less than exact dissipation for moderate wavenumbers ($\pi/2 < K < 3\pi/4$) while they add more dissipation at the Nyquist frequency. The \bff{BR1}-$\eta0$ scheme is always away less robust than the \bff{BR2}-$\eta1$ scheme due areas of less than exact dissipation and also as we mentioned that the diffusion factors at the Nyquist frequency almost saturates at these levels in~\hfigref{fig:BR1_compareP_G_tau0.5}. Although \bff{LDG}-$\eta0$ schemes have higher dissipation than the exact one for short time simulations, they experience less dissipation for long time simulations at high wavenumbers $K > \pi/2$ especially p$2$ to p$5$ orders. This scheme indeed has the most accurate and monotonic behavior for long time simulations.

Finally, the relative efficiencies between the three viscous flux approaches are further discussed by comparing the resolution of all schemes at the same $\tau_{p}$. For the \bff{BR1} and \bff{BR2} approaches we choose a penalty parameter that gives almost the same diffusion factor at the Nyquist frequency as the \bff{LDG}-$\eta0$ when $\tau_{p}=0.01$ (short time). This diffusion factor provides higher than the exact dissipation for $K \gtrsim \pi/2$. This way we include a possibly robust versions of these schemes.

\begin{figure}[t]
    \begin{subfigure}[h]{0.5\textwidth}
    \includegraphics[width=0.965\textwidth]{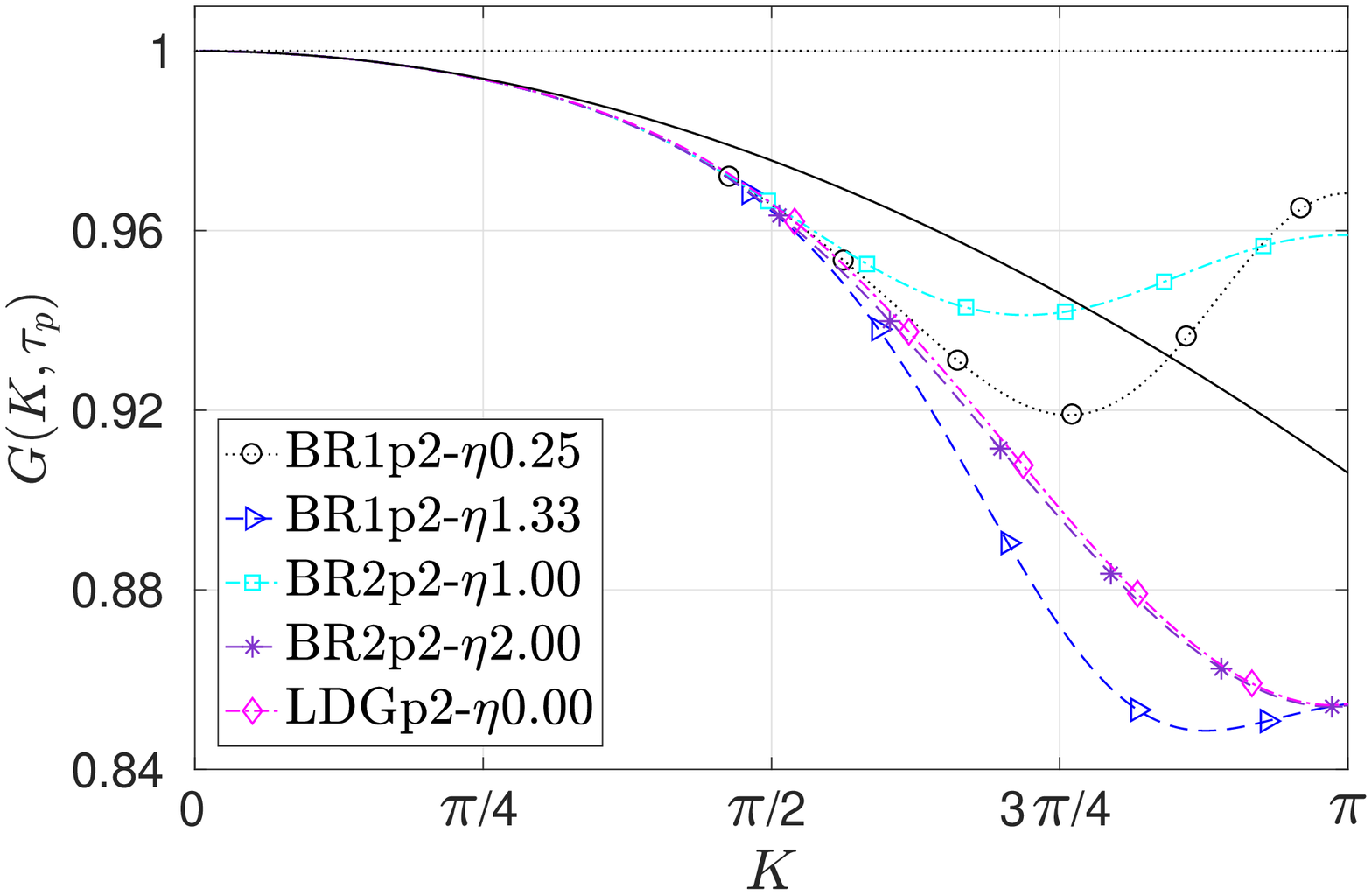} 
    \caption{Diffusion factor, $\tau_{p}=0.01$} 
    \label{fig:sdisc_compareP2schemes_tau0.01_G}
    \end{subfigure}
    \, \,
   \begin{subfigure}[h]{0.5\textwidth}
    \includegraphics[width=0.965\textwidth]{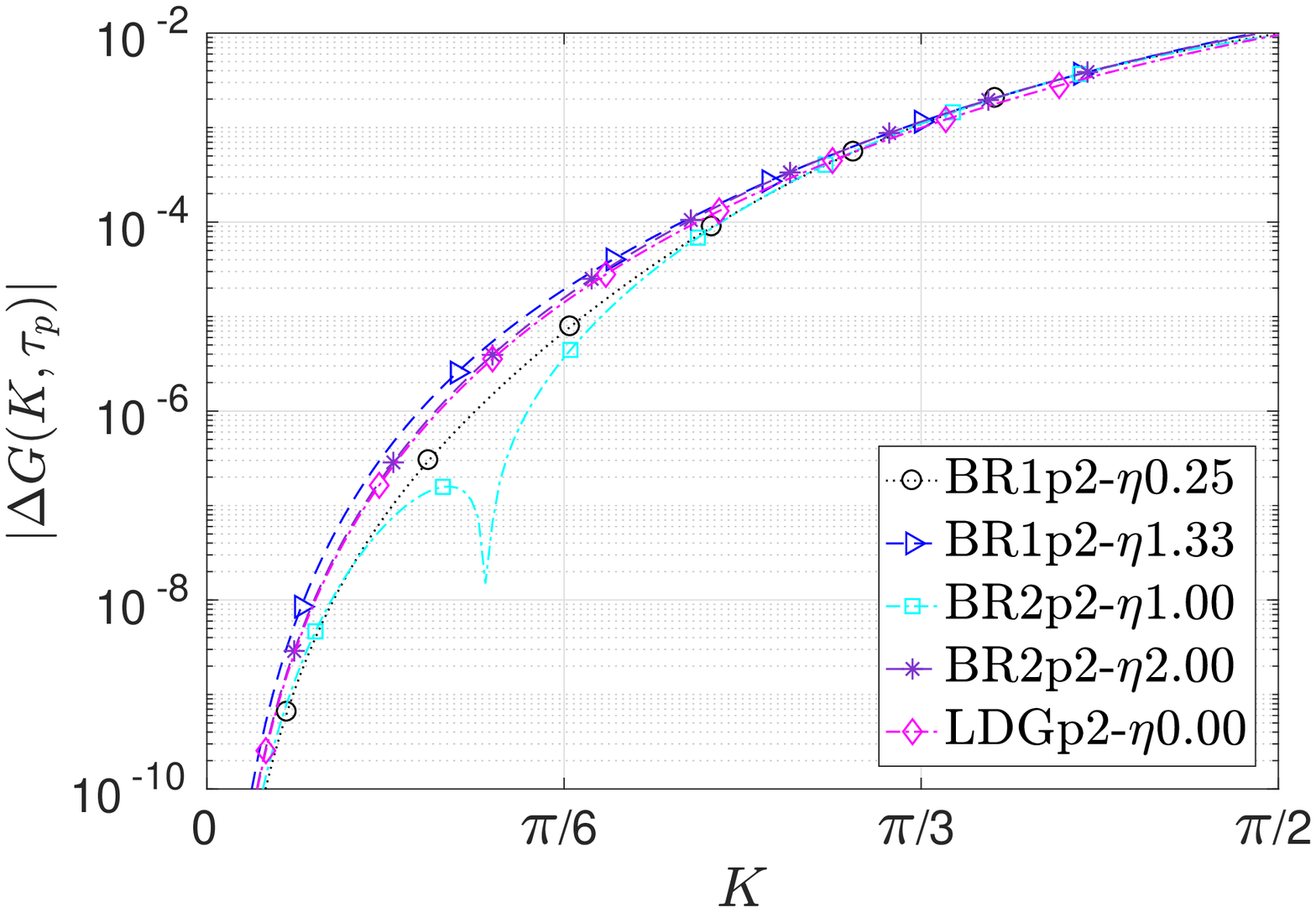}  
    \caption{Diffusion error, $\tau_{p}=0.01$} 
    \label{fig:sdisc_compareP2schemes_tau0.01_Gerr}
    \end{subfigure} 
    \\ \\ \\
    \begin{subfigure}[h]{0.5\textwidth}
    \includegraphics[width=0.965\textwidth]{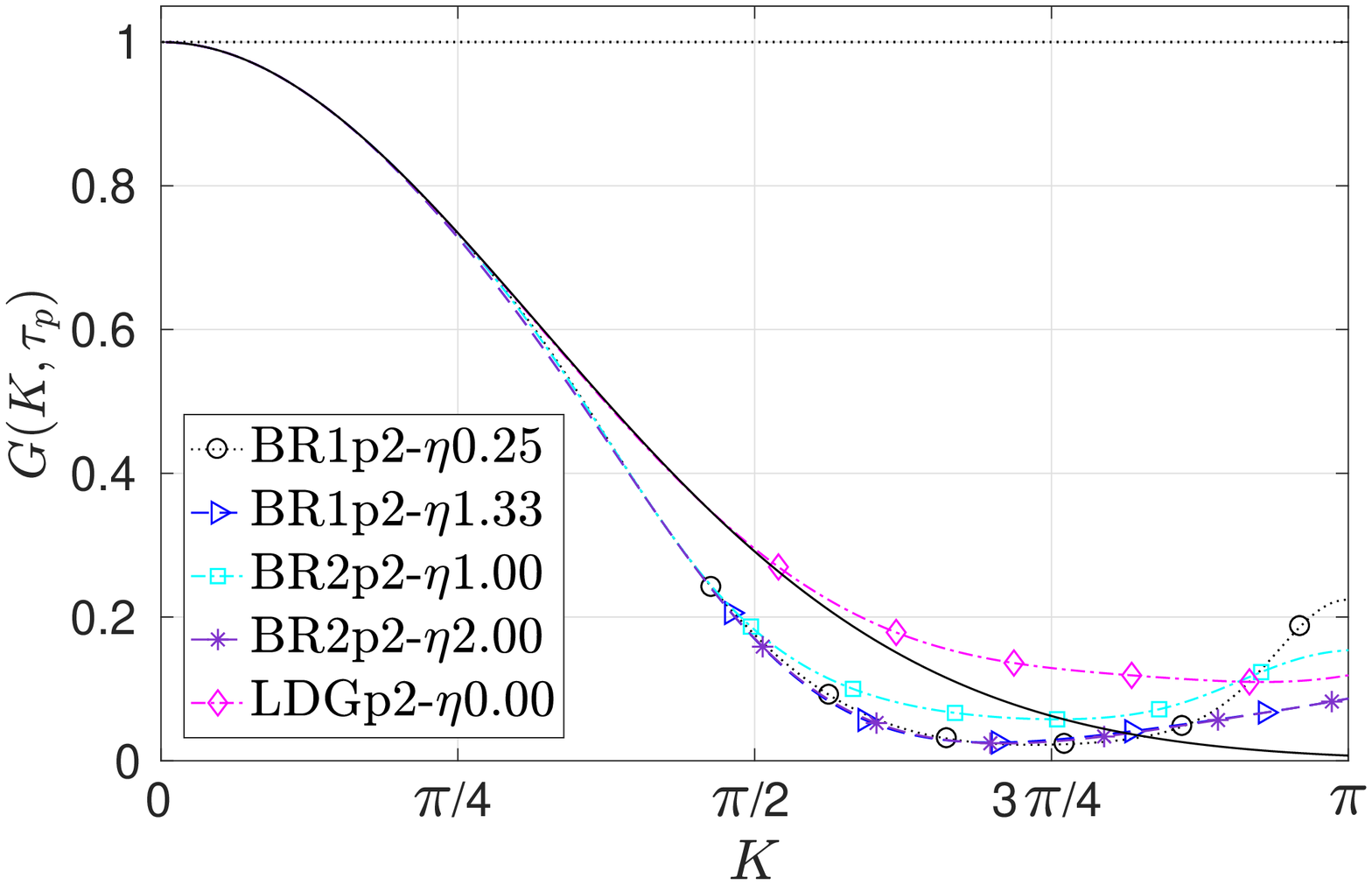} 
    \caption{Diffusion factor, $\tau_{p}=0.5$} 
    \label{fig:sdisc_compareP2schemes_tau0.5_G}
    \end{subfigure}
    \, \,
   \begin{subfigure}[h]{0.5\textwidth}
    \includegraphics[width=0.965\textwidth]{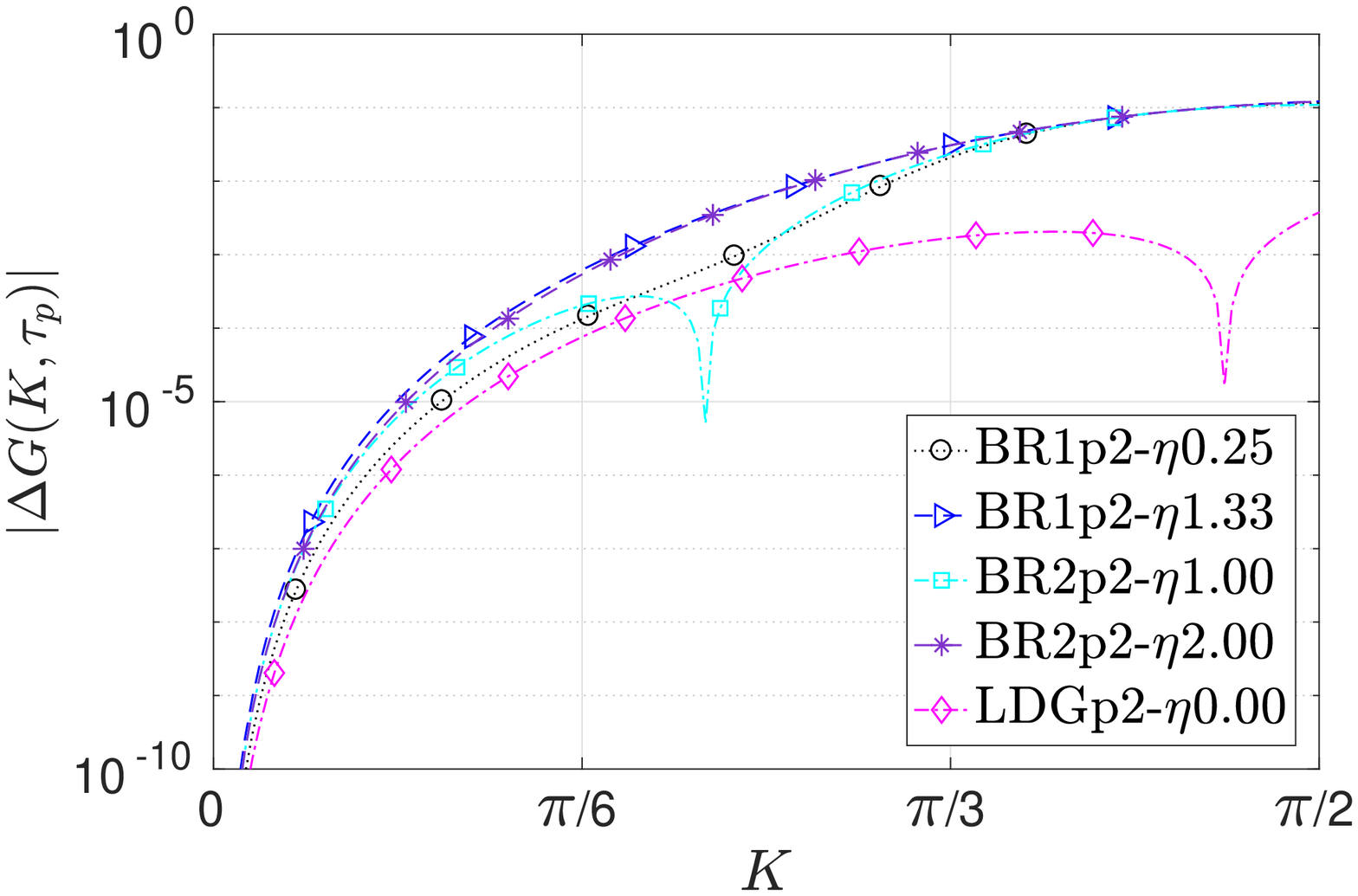}  
    \caption{Diffusion error, $\tau_{p}=0.5$} 
    \label{fig:sdisc_compareP2schemes_tau0.5_Gerr}
    \end{subfigure} 
    \caption{Comparison of the diffusion behavior of several p$2$ diffusion schemes. In the left figure, the solid line without symbols represents the exact diffusion factor $G_{ex}=e^{-K^{2} \tau_{p}}$. }
\label{fig:sdisc_compareP2schemes}
\end{figure}

In~\hfigref{fig:sdisc_compareP2schemes} we present the comparisons for short and long times, i.e., at $\tau_{p}=0.01, \, \tau_{p}=0.50$. From this figure it can be noticed that with careful tuning of $\eta$ we can achieve very similar behaviors between different schemes. According to~\hfigref{fig:sdisc_compareP2schemes_tau0.01_G}, the diffusion factor of \bff{BR2}p$2$-$\eta2$ is very close to the \bff{LDG}-$\eta0$ one. In addition, the diffusion factor of the \bff{BR1}p$2$-$\eta1.33$ scheme is very close to that of the \bff{BR2}p$2$-$\eta2$ scheme for both short and long times. This can be understood in the context of similarities between them through~\heqssref{eqn:SIP_expandedPrimal}{eqn:BR1_stabilized_Rg_simplified}{eqn:BR2_SIPG_equiv_Cp} where they involve the same form of the jump penalizing term. However, it is worth noting that their behavior for $K > \pi/2$ is not that similar, possibly due to the non-compactness of \bff{BR1} induced by additional terms in~\heqref{eqn:BR1_stabilized_Rg_simplified}. 

\begin{figure}[H]
    \begin{subfigure}[h]{0.5\textwidth}
    \includegraphics[width=0.965\textwidth]{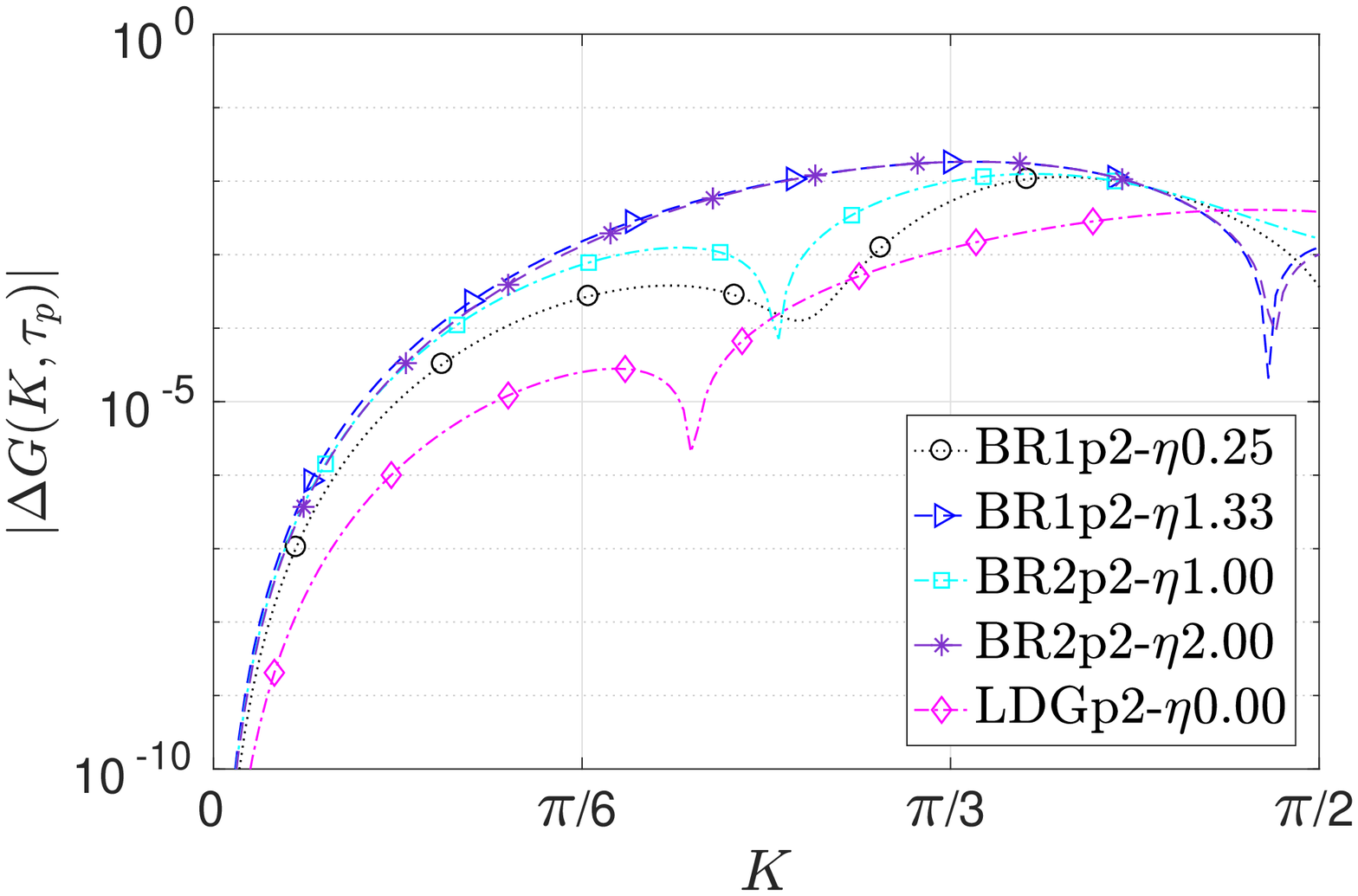}  
    \caption{p$2$ schemes} 
    \label{fig:sdisc_compareP2schemes_tau2_Gerr}
    \end{subfigure} 
    \, \,
   \begin{subfigure}[h]{0.5\textwidth}
    \includegraphics[width=0.965\textwidth]{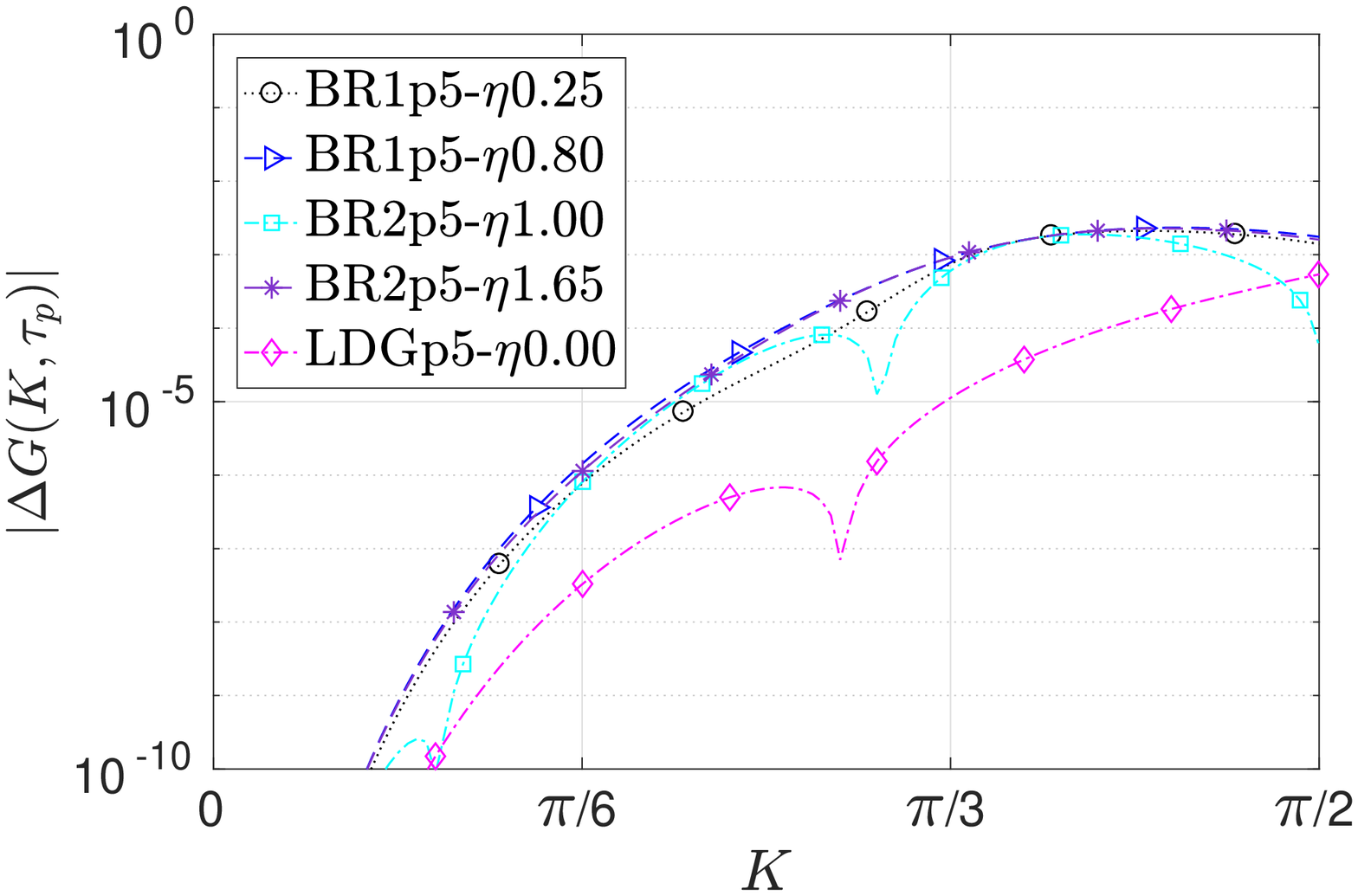}  
    \caption{p$5$ schemes} 
    \label{fig:sdisc_compareP5schemes_tau2_Gerr}
    \end{subfigure} 
    \caption{Comparison of the diffusion errors of p$2$ and p$5$ diffusion schemes at very long time, $\tau_{p}=2.0$.}
\label{fig:sdisc_compareP2P5schemes_tau2}
\end{figure}

 We also note that matching the same amount of dissipation at the Nyquist frequency requires different $\eta$ values for different orders, e.g., for p$5$ this can be achieved with $\eta^{br1}=0.80, \, \eta^{br2}=1.65$. This suggests that schemes with asymptotically similar diffusion factors can also be achieved for high-orders. One particularly interesting scheme is the \bff{BR1}p$2$-$\eta0.25$ which serves as a good alternative for the standard \bff{BR1}p$2$-$\eta0$ with even less low wavenumber errors than \bff{BR1}p$2$-$\eta1.33$,~\hfigref{fig:sdisc_compareP2P5schemes_tau2}. Despite the similarities between these schemes, the standard \bff{LDG} scheme maintained a better error bound than \bff{BR1} and \bff{BR2} schemes for low wavenumbers, $K \leq \pi/2$, see~\hfigref{fig:sdisc_compareP2P5schemes_tau2}.  %
%
%
 %
 \begin{figure}[H]
\centering
    \includegraphics[width=0.55\textwidth]{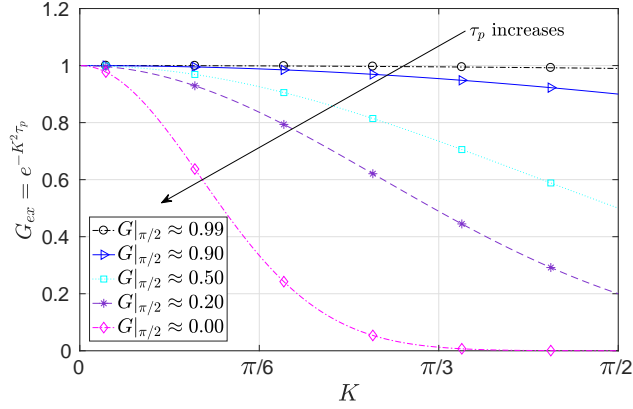} 
    \caption{Evolution of the exact semi-discrete diffusion factor with time. }
    \label{fig:Gex_pi_2_evolution}
\end{figure}%
%
\subsubsection{Evolution of the diffusion error with time} \label{sec:sdisc_time_err_evolution}
In this section we study the evolution of the diffusion error with time, for all schemes in the low wavenumber range, $K \leq \pi/2$. Based on the results of the previous section, we choose two values for $\eta$, namely, $\eta^{br1}=0.0,\, \eta^{br1}=0.25$, for the DG-\bff{BR1} scheme while using the standard form for the DG-\bff{BR2} and \bff{LDG} schemes with minimum stabilization. The error is computed at several instances in time where the diffusion factor $G$ of the initial Fourier mode at $K = \pi/2$ changes from $G|_{\pi/2} \approx 0.99$ to $G|_{\pi/2} \approx 0.00$,~\hfigref{fig:Gex_pi_2_evolution}. These values were chosen as to highlight the importance of the resolution in the low wavenumber range for short and long times.

\begin{figure}[h]
\begin{subfigure}[h]{0.5\textwidth}
    \includegraphics[width=0.975\textwidth]{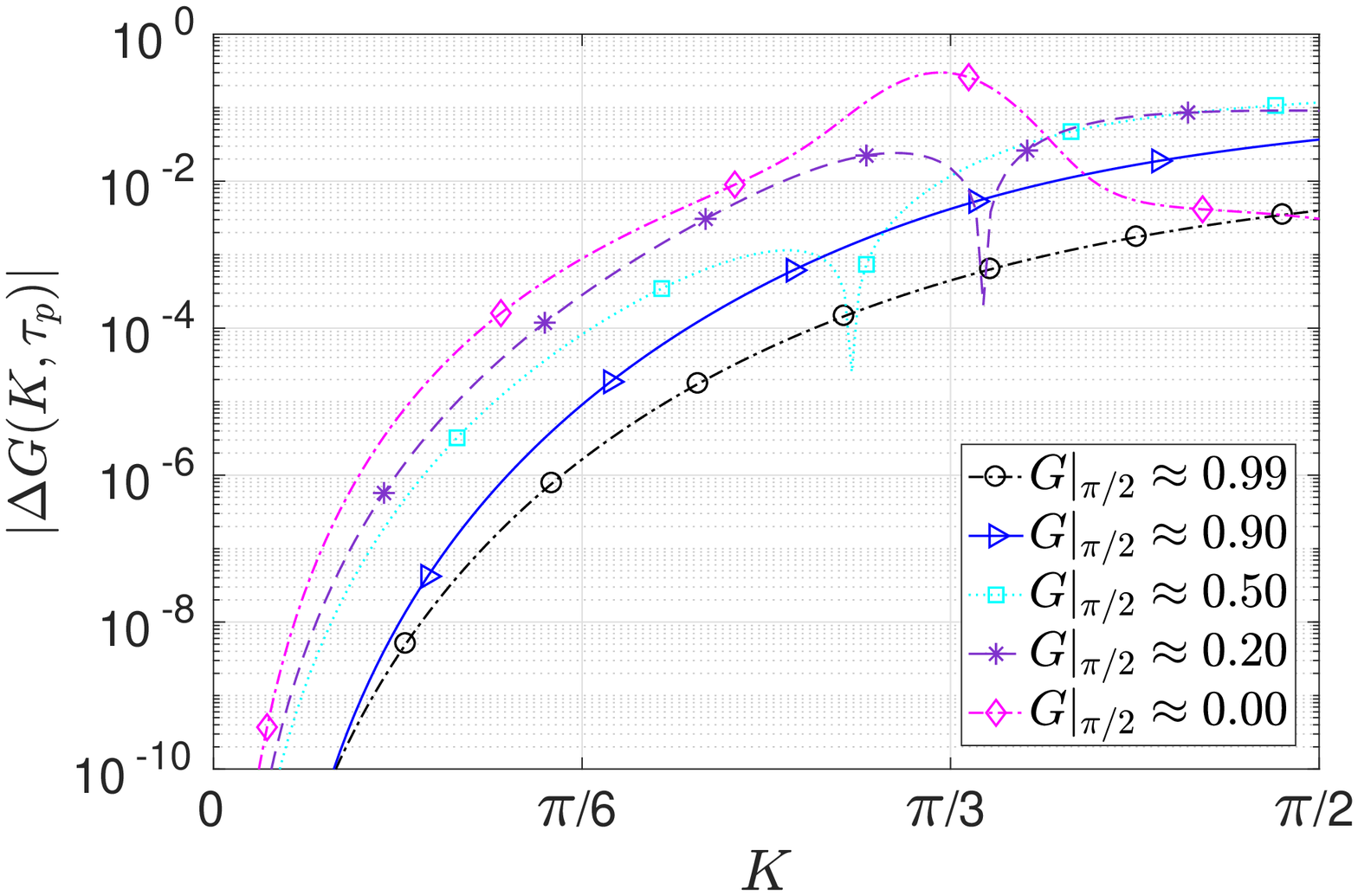} 
    \caption{Diffusion error for BR$1$p$2$-$\eta0.0$} 
    \label{fig:BR1p2_eta0.0_error_evolution}
    \end{subfigure}
    \, \,
   \begin{subfigure}[h]{0.5\textwidth}
    \includegraphics[width=0.975\textwidth]{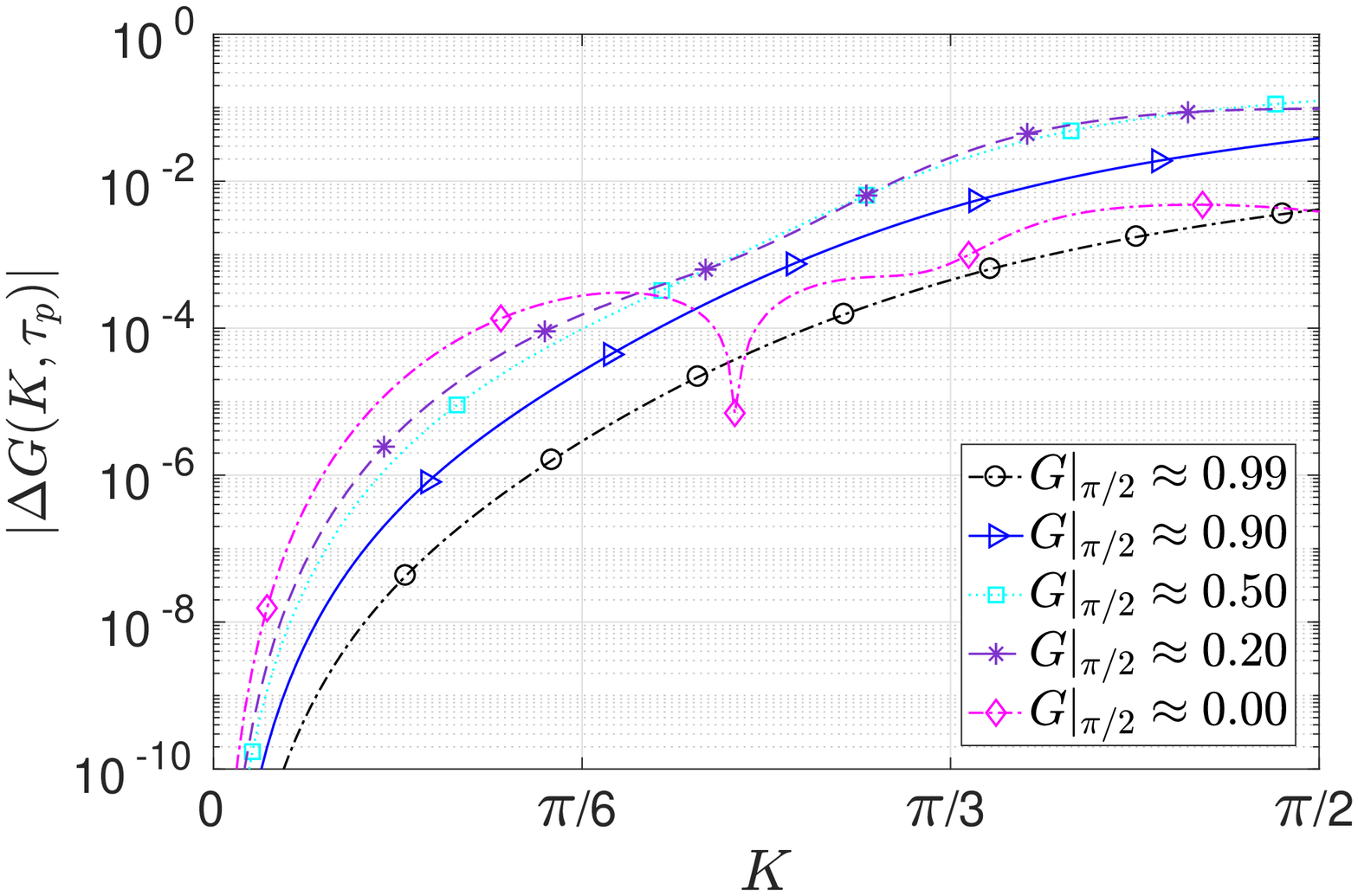}  
    \caption{Diffusion error for BR$1$p$2$-$\eta0.25$} 
    \label{fig:BR1p2_eta0.25_error_evolution}
    \end{subfigure}  \\ \\ \\
\, \,
    \begin{subfigure}[h]{0.5\textwidth}
    \includegraphics[width=0.975\textwidth]{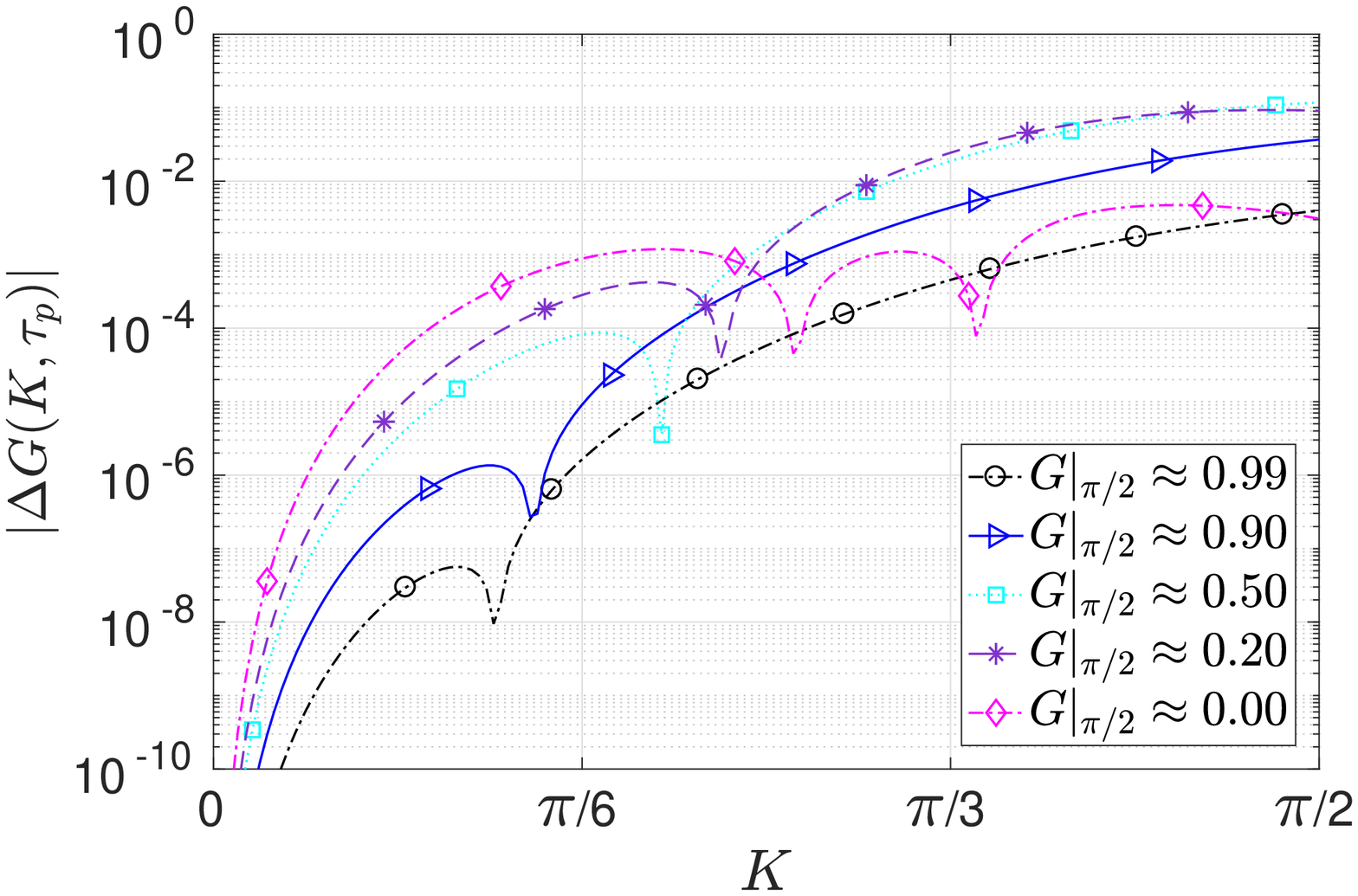} 
    \caption{Diffusion error for BR$2$p$2$-$\eta1.00$} 
    \label{fig:BR2p2_eta1_error_evolution}
    \end{subfigure}
    \, \,
   \begin{subfigure}[h]{0.5\textwidth}
    \includegraphics[width=0.975\textwidth]{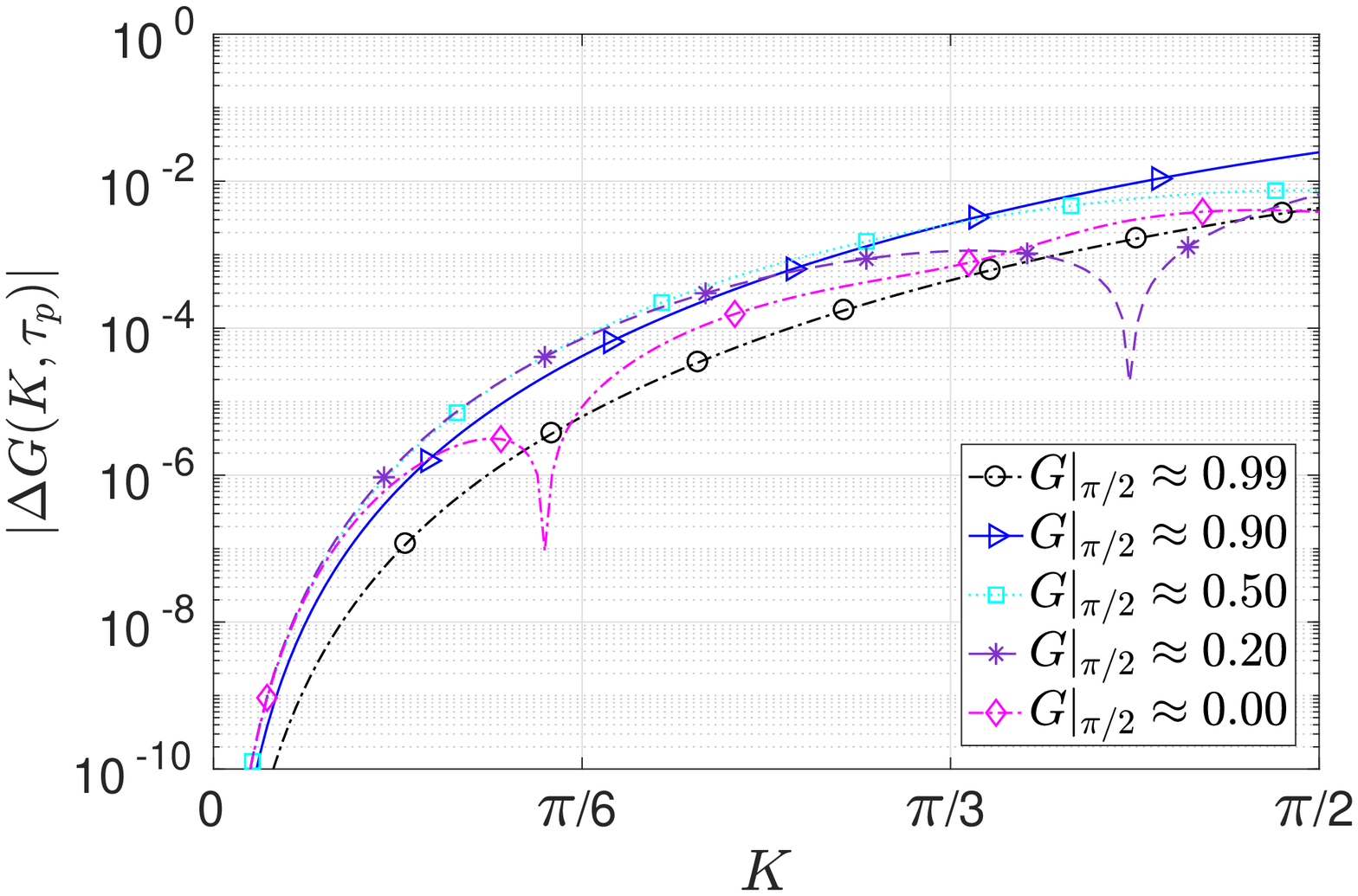}  
    \caption{Diffusion error for LDGp$2$-$\eta0.0$} 
    \label{fig:LDGp2_eta0_error_evolution}
    \end{subfigure}  
	\caption{Evolution of the semi-discrete diffusion error with time for p$2$ spatial schemes.}
    \label{fig:ErrorTimeEvolution_sdisc}
\end{figure}

The evolution of the error with time for a number of schemes is displayed in~\hfigref{fig:ErrorTimeEvolution_sdisc}. In this figure we can see that the \bff{LDG}p$2$-$\eta0$ scheme maintains an upper bound on the error that is not exceeded as time evolves. The maximum error of this scheme for $0 < K \lesssim \pi/4$ does not exceed $\approx 10^{-4}$ while for the \bff{BR1}p$2$-$\eta0$, \bff{BR1}p$2$-$\eta0.25$  and \bff{BR2}p$2$-$\eta1$ schemes the error clearly exceeds this level. The diffusion errors for the \bff{BR1}p$2$-$\eta0$, \bff{BR1}p$2$-$\eta0.25$  and \bff{BR2}p$2$-$\eta1$ schemes always increases as time evolves. It is also noticed the slightly stabilized \bff{BR1}p$2$-$\eta0.25$ scheme improves the properties of the standard \bff{BR1}p$2$-$\eta0$ scheme. %
%
\section{Fully-discrete Fourier Analysis}\label{sec:fdiscAnalysis}%
%
In order to analyze the wave propagation properties of a fully-discrete scheme, both time and space discretizations are applied to the linear heat equation~\heqref{eqn:heat_eqn}. Applying a RK time discretization scheme to one of the semi-discrete equations in~\secref{sec:num_methods} results in an update formula for the solution at $t^{n+1}= (n+1) \Delta t$ of the form%
\begin{equation}
\bs{u}^{n+1} =\mycalmat{G} \bs{u}^{n} = \mycalmat{G}^{n+1} \bs{u}^{0} , \quad \bs{u}^{0} = \bs{u}(x,0), \quad \mycalmat{G}= \mathcal{P} (\gamma \Delta t \mycalmat{A}) ,
\label{eqn:fdisc_temp_form} 
\end{equation}%
where $\mycalmat{G}$ is the fully-discrete operator, $\mathcal{P}$ is the time integration polynomial defined in~\heqref{eqn:RK_update_form_DG}, $\mycalmat{A} $ is the semi-discrete operator, and $\bs{u}^{n+1}$ is the vector of unknown DOFs. In addition, the eigenvalues of $\mymat{\mathcal{G}}$ can be expressed as%
\begin{equation}
\lambda_{\mymat{\mathcal{G}}} = \mathcal{P}(\gamma \Delta t \lambda_{\mycalmat{A}}),
\label{eqn:G_Lambda_Poly}
\end{equation}%
thus the behavior of the fully-discrete scheme depends on both the time-step $\Delta t$ and the form of the polynomial $\mathcal{P}$. The method of analysis used in this section is similar to the one used by Guo et al.~\cite{GuoSuperconvergencediscontinuousGalerkin2013} for analyzing the \bff{LDG} method. We proceed by seeking a wave solution in element $\Omega_{e}$ of the form %
\begin{equation}
\mathbf{U}^{e,n+1} = \bs{\mu} \: e^{i k x_{e} - \tilde{\omega} t^{n+1} } = e^{- i \tilde{\omega} \Delta t} \mathbf{U}^{e,n} ,
\label{eqn:blochwave_sol_fdisc}
\end{equation}%
and substituting in~\heqref{eqn:fdisc_temp_form} yields the following fully-discrete relation for element $\Omega_{e} \in \mathcal{D}$%
\begin{equation}
e^{- \tilde{\omega} \Delta t} \bs{\mu} = \mycalmat{G}\bs{\mu} ,
\label{eqn:DG_full_disc}
\end{equation}%
where $\tilde{\omega}$ is the numerical frequency which is a real number in this analysis. For a non-trivial solution the determinant of the following equation should be zero%
\begin{equation}
	\det \left(\mycalmat{G}- e^{- \tilde{\omega} \Delta t} \, \mycalmat{I} \right)  =  0.
\label{eqn:eigvalue_problem_fdisc}
\end{equation}%
This equation constitutes an eigenvalue problem similar to the case of semi-discrete analysis. The general representation of the vector of spatial solution coefficients can be expressed as a linear combination of the eigenmodes as follows%
\begin{equation}
\mathbf{U}^{e,n} = \sum_{j=0}^{p} \vartheta_{j} \lambda_{j}^{n} \bs{\mu}_{j} \: e^{i k x_{e} } =  \sum_{j=0}^{p} \vartheta_{j} \bs{\mu}_{j} \: e^{i  k x_{e} - \tilde{\omega}_{j} (n \Delta t) }  ,
\label{eqn:lin_combination_eigvec_fdisc_Un}
\end{equation}%
where the expansion coefficients $\vartheta_{j}$ are again given by~\heqref{eqn:eigenmode_weights}, and the general solution has the same form as~\heqsref{eqn:sdisc_fullsolform}{eqn:exact_sdisc_fullsolform}. To this end, it is clear that the eigenvalues, $\lambda_{j}$, are related to the numerical frequency $\tilde{\omega}$ through the following relation%
\begin{equation}
\lambda_{j} = e^{-\tilde{\omega}_{j} \Delta t} = e^{- \frac{\gamma \Delta t}{h^{2}} \frac{\tilde{\omega}_{j} }{\gamma} \frac{h^{2}}{(p+1)^{2}} (p+1)^{2} }  = e^{- (p+1)^{2} \Delta \tau \,  K_{m,j}^{2}  } , \quad j=0,...,p ,
\label{eqn:fdisc_lambda_K_rel1}
\end{equation}%
where $\Delta \tau=\gamma \Delta t/h^{2}$ is the non-dimensional time step, and the numerical wavenumber $K_{m}$ can be obtained from%
\begin{equation}
K^{2}_{m} = \frac{\: ln(\lambda)}{\Delta \tau_{p}}, \quad \Delta \tau_{p} = (p+1)^{2} \Delta \tau.
\label{eqn:fdisc_lambda_K_rel2}
\end{equation}%
Therefore, a numerical dissipation relation can be written as %
\begin{equation}
-(\operatorname{\mathcal{R}e}(K^{2}_{m}) ) \equiv - K_{m}^{2} \approx -K^{2} ,
\label{eqn:DG_num_full-disc_disper_relation}
\end{equation}%
while the numerical dispersion is  zero, i.e., $\operatorname{\mathcal{I}m} (K^{2}_{m})=0$.%
\begin{figure}[H]
\centering
    \includegraphics[width=0.5\textwidth]{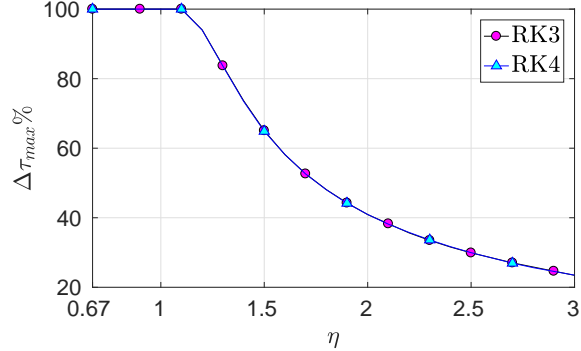} 
	\caption{Percentage ratios of $\Delta \tau_{max}(\eta)/\Delta \tau_{max}(\eta_{min})$ for the BR$2$p$2$ coupled with RK schemes. }
    \label{fig:BR2p2_maxdt_Eta}
\end{figure}%
In the previous sections we have identified minimum values of $\eta$ for the stability of DG diffusion schemes in the semi-discrete sense. The influence of these different $\eta$ on the performance was also studied. In this section, we investigate the impact $\eta$ has on the maximum time-step required for stability of a viscous DG scheme coupled with explicit RK for time integration. This is done in the same procedure as in~\secref{sec:sdisc_eigenmode_analysis}. Starting with an arbitrary large $\Delta \tau$, assume a prescribed wavenumber $K \in [0,\pi]$ and for a certain $\eta \geq \eta_{min}$ we check the sign of $-(K_{m})^{2}$ for all eigenmodes. If the sign is positive, i.e., $-(K_{m})^{2}>0$ this indicates a growing mode and the scheme is declared unstable for this $\Delta \tau$. Afterwards, we start to decrease the $\Delta \tau$ until we have all $-(K_{m})^{2}<0$ and this $\Delta \tau = \Delta \tau_{max}$, the stability limit for this particular scheme and $\eta$. As $\eta$ increases we noticed that $\Delta \tau_{max}$ always decreases which is counter intuitive with the effect of $\eta$ adding more dissipation. However, this is true because this dissipation is kind of wavenumber diffusion through the eigenvalues and not an even-order truncation error term. In other words, as the spectral radius of the $\mymat{\mathcal{G}}$ matrix increases, the convergence of the update equation becomes slower (can be interpreted as a point iteration). %

~\hfigref{fig:BR2p2_maxdt_Eta} presents the stability limits of the \bff{BR2}p$2$ scheme with a range of $\eta$ as a percentage of the $\max\limits_{\eta} \Delta \tau_{max}$.  From this figure it can be seen that for $\eta=\eta_{min} \approx 0.67$ to $\eta \approx 1.1$, the scheme has a constant maximum time step. Afterwards, the $\Delta \tau_{max}$ decreases up to $\approx 40\%$ of its maximum at $\eta=2$. This behavior is not dependent on the RK scheme in use and it only shifts to the right with the same form as the order $p$ increases due to increasing $\eta_{min}$. It is only p$1$ polynomial order where the plateau region of constant $\Delta \tau_{max}$ is from $\eta_{min}$ to $0.8$ instead of $1.1$. Maximum time steps for different orders, RK schemes, and with different $\eta$ for \bff{BR2} schemes are given in~\htablref{table:BR2_stability_limits_1D},~\happref{appx:A}{A}. In the semi-discrete~\combined analysis it was shown that far from $\eta_{min}$ is always recommended for a robust scheme with sufficient high frequency dissipation. However, a sever loss of $\Delta \tau_{max}$ of  about $50\%$ is expected, making explicit time integration stiffer. For the \bff{BR1} and \bff{LDG} approaches the dependence of $\Delta \tau_{max}$ on the $\eta$ parameter is nearly linear (not shown) with a negative (relatively small) slope, i.e., it has its maximum at $\eta_{min}$ for each scheme. This means that increasing $\eta$ reduces the time step but with a slower rate than in the case of \bff{BR2}. 

The von Neumann stability limits for all RKDG schemes considered in this work with minimum stabilization (standard form) for the linear heat equation is provided in~\htablref{table:BR1BR2LDG_stability_limits_1D},~\happref{appx:A}{A}. From this table it can be noticed that $\Delta \tau_{max,RK3} \approx 0.9 \times \Delta \tau_{max,RK4}$ for all schemes and orders. %
%
\subsection{Combined-mode analysis}\label{sec:fdisc_combinedMode}%
%
In the fully-discrete~\combined analysis of DG schemes for diffusion we utilize the same definitions introduced in~\secref{sec:sdisc_combinedMode} for the~\mytrue energy and diffusion factor. However, the solution DOFs are now given by~\heqref{eqn:lin_combination_eigvec_fdisc_Un}. The \physic diffusion factor of the fully-discrete scheme is given by%
\begin{equation}
G_{phy} (k, \tau_{p,n}) = e^{-n \Delta \tau_{p} \,  K_{m}^{2}  } , \quad \tau_{p,n} = \tau_{p,0} + n \Delta \tau_{p}, \quad \tau_{p,0} = 0.
\label{eqn:Gphys_fdisc_combined}
\end{equation}%
In this case the diffusion factor is dependent on both the selected $\Delta \tau$ and the number of iterations $n$. The objective of this section is to investigate the impact of $\Delta \tau$ on the performance of various schemes in the fully-discrete case. 

\begin{figure}[H]
 \begin{subfigure}[h]{0.5\textwidth}
    \includegraphics[width=0.975\textwidth]{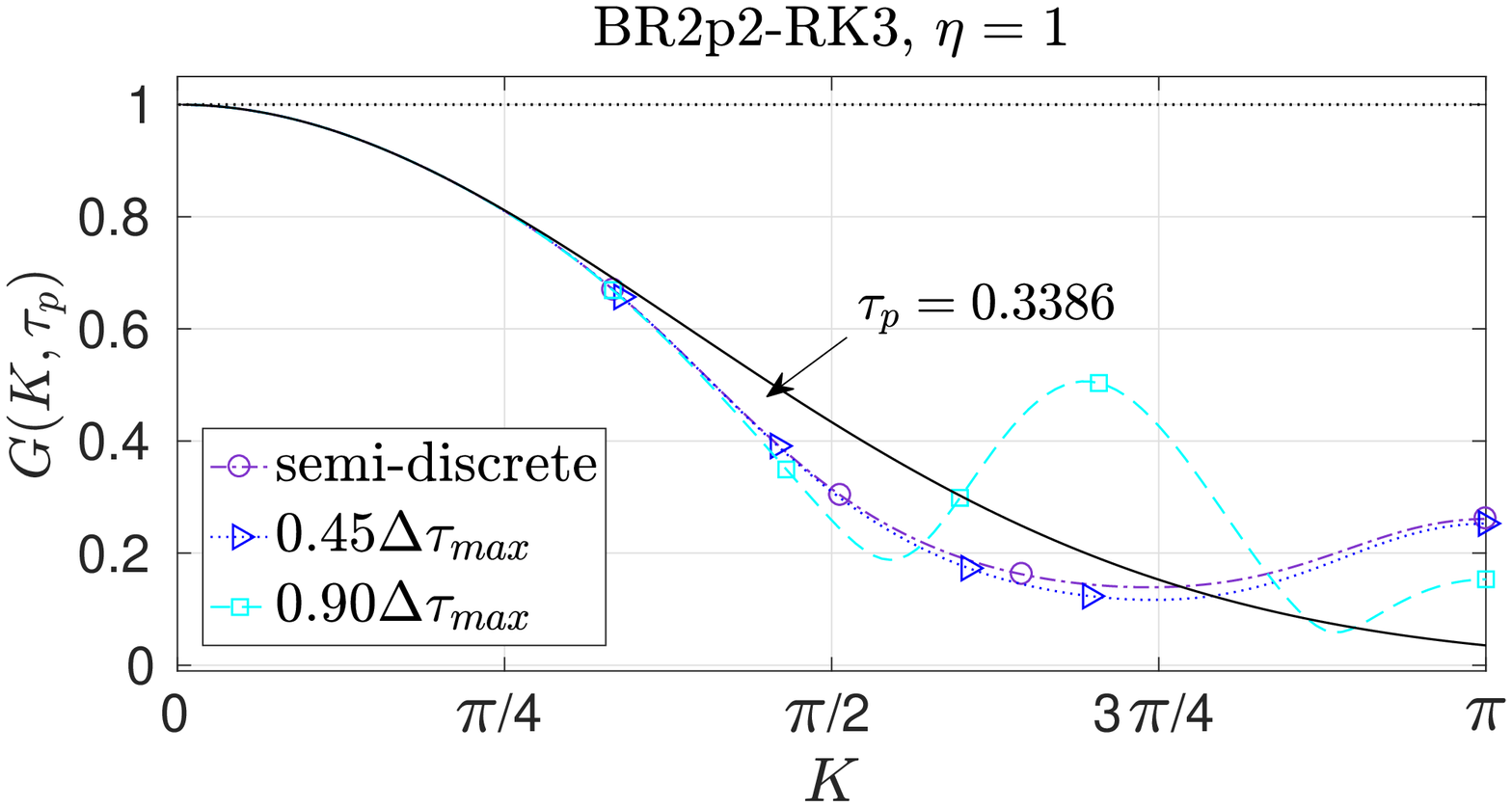} 
    \label{fig:BR2_fdisc_G_p2RK3}
    \end{subfigure}
    \, 
     \begin{subfigure}[h]{0.5\textwidth}
    \includegraphics[width=0.975\textwidth]{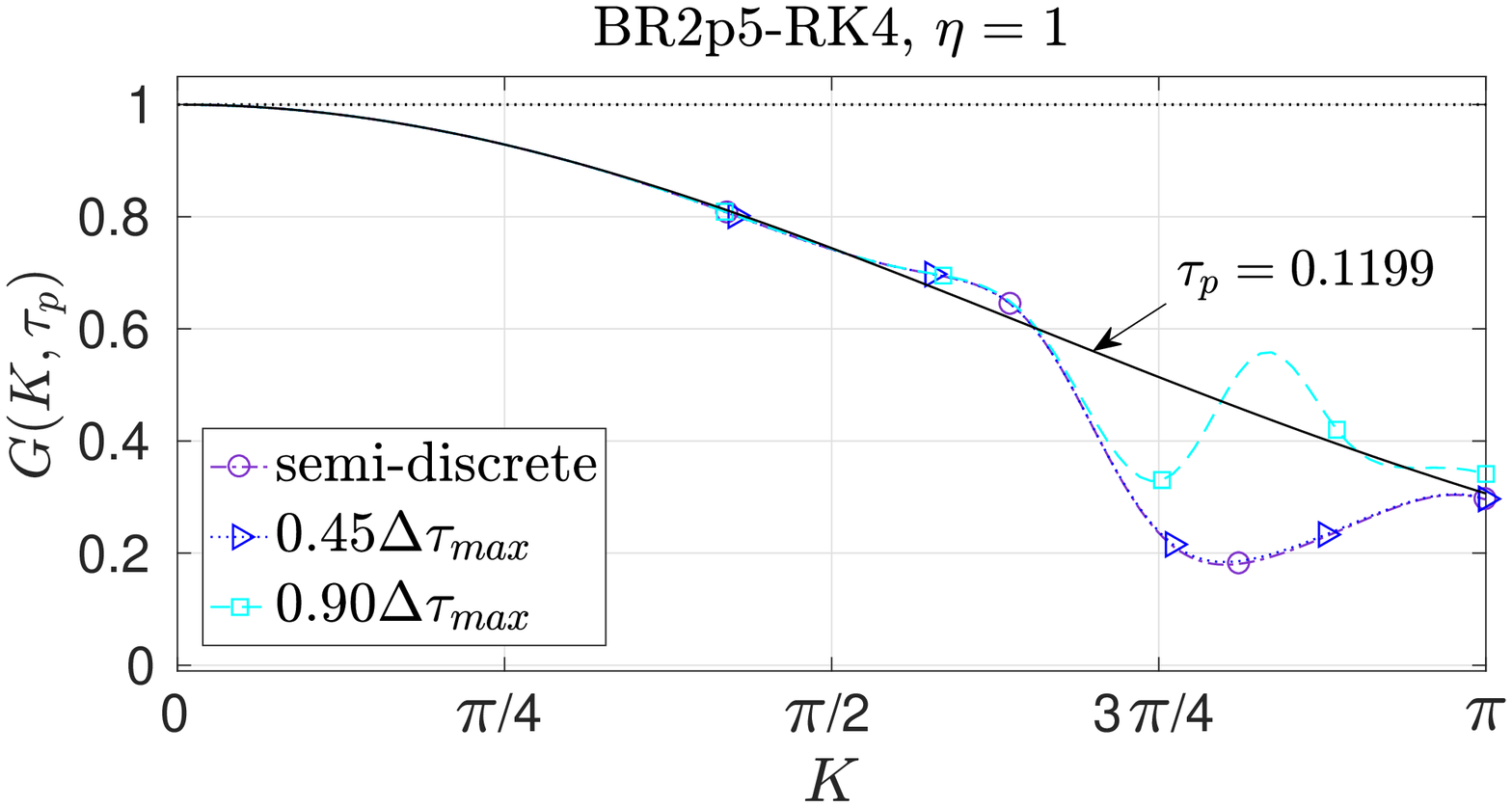} 
    \label{fig:BR2_fdisc_G_p5RK4}
    \end{subfigure} \\ \\ \\
    \begin{subfigure}[h]{0.5\textwidth}
    \includegraphics[width=0.975\textwidth]{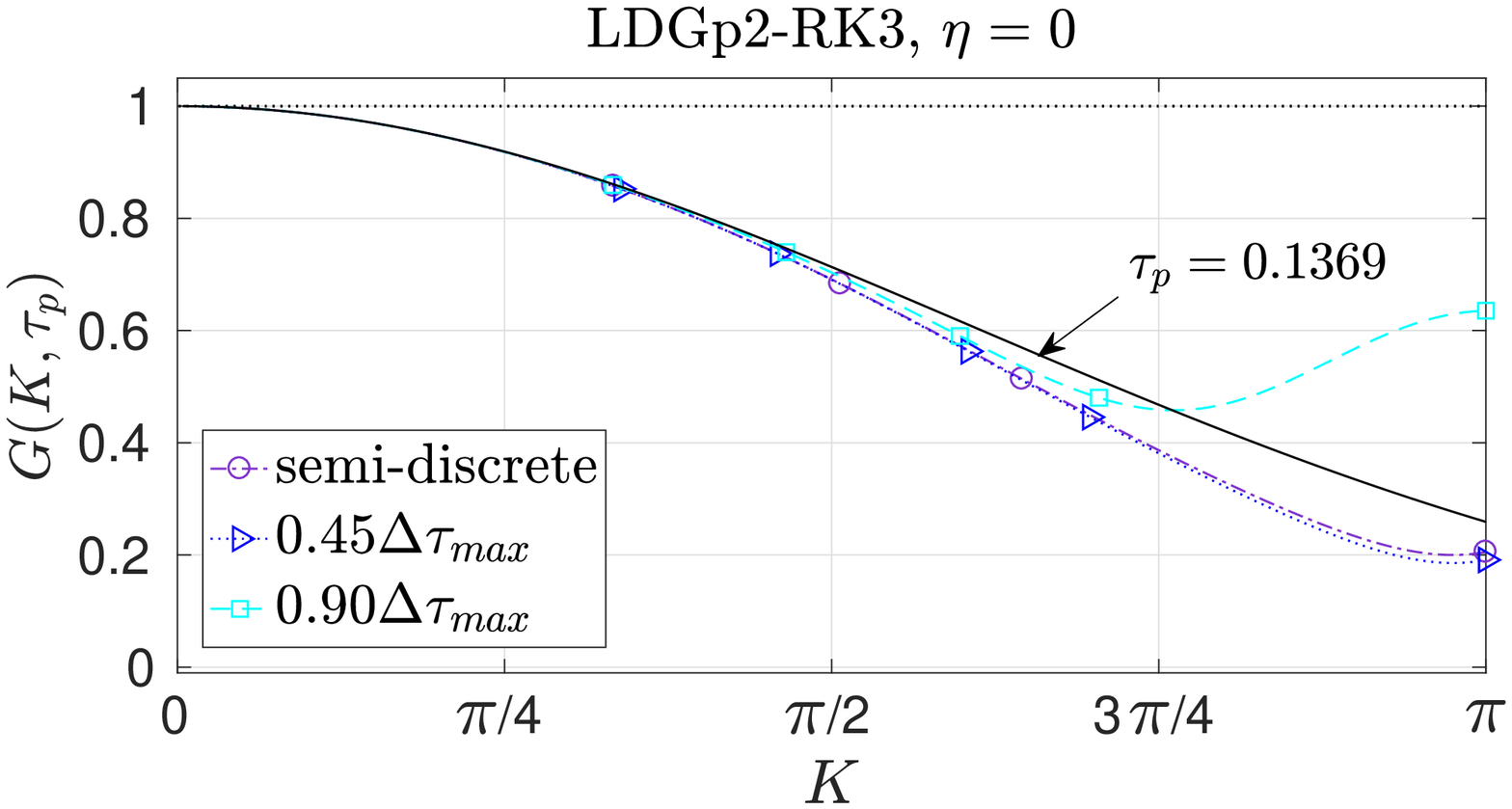} 
    \label{fig:LDG_fdisc_G_p2RK3}
    \end{subfigure}
    \, 
     \begin{subfigure}[h]{0.5\textwidth}
    \includegraphics[width=0.975\textwidth]{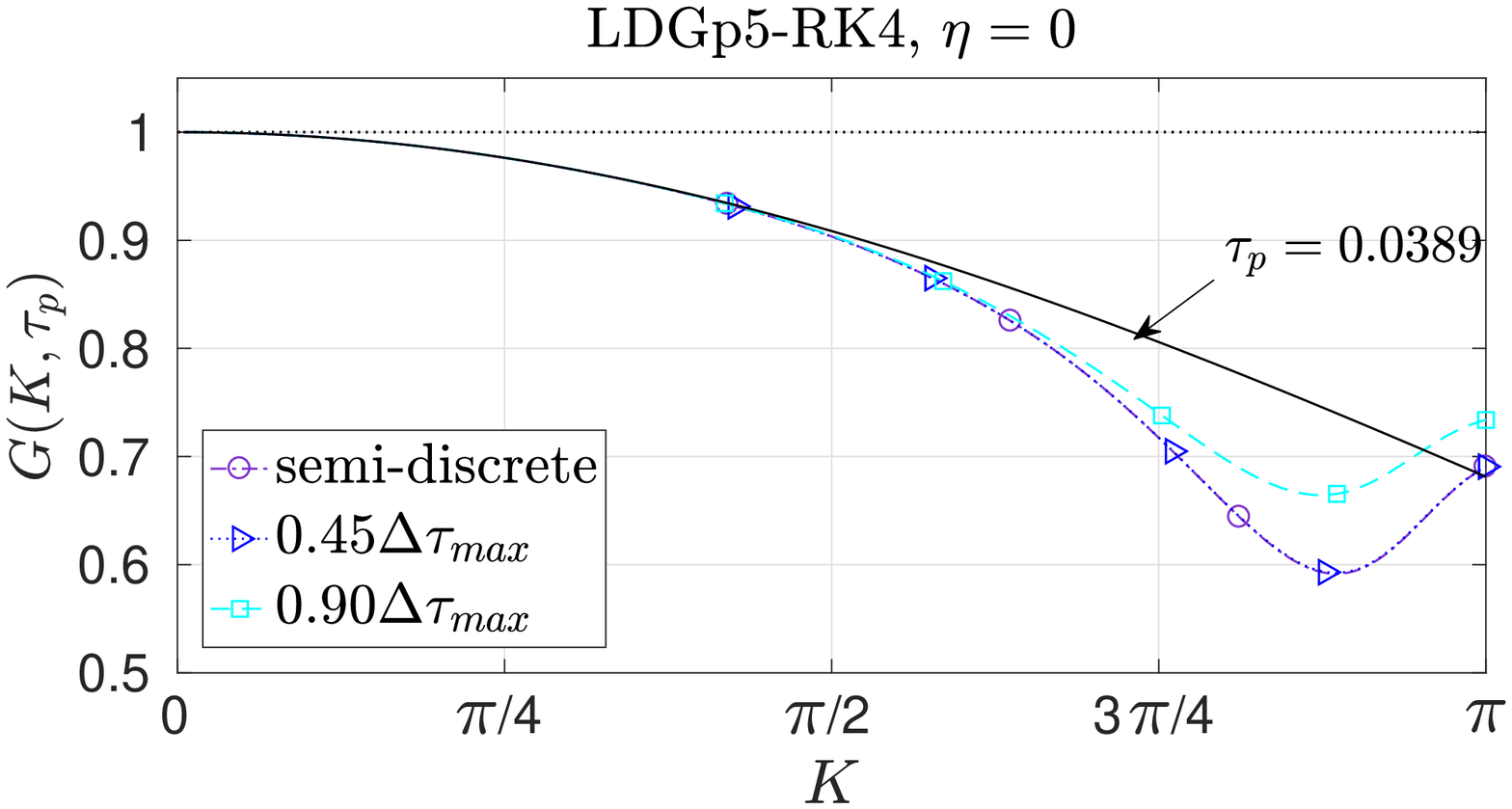} 
    \label{fig:LDG_fdisc_G_p5RK4}
    \end{subfigure}
    \caption{Diffusion factors for fully-discrete schemes with different $\Delta \tau$.  In these figures, the solid black line without symbols represents the exact diffusion factor $G_{ex}=e^{-K^{2} \tau_{p}}$.}
    \label{fig:BR2LDG_fdisc_G}
    \end{figure}

From~\hfigref{fig:BR2LDG_fdisc_G}, we can see that the effect of the $\Delta \tau$ is more pronounced at the high wavenumber range while there is almost no influence on the behavior for the low wavenumber range as can be seen by computing the diffusion error (not shown). This also indicates that near $\Delta \tau_{max}$ the behavior of the \bff{BR2}-$\eta1$ coupled with RK schemes is less robust than the semi-discrete case providing less dissipation at some moderate wavenumbers. However, in a non monotonic and oscillatory way. The influence of increasing $\Delta \tau$ is the same for \bff{LDG}-$\eta0$ coupled with RK schemes but without any oscillations. It is also worth noting that at $\approx 0.5 \, \Delta \tau_{max}$ all schemes become very close to their semi-discrete versions. We can also notice that indeed p$5$ schemes have a better performance than p$2$ schemes as have been seen in the semi-discrete case.%
\begin{figure}[H]
 \begin{subfigure}[h]{0.5\textwidth}
    \includegraphics[width=0.975\textwidth]{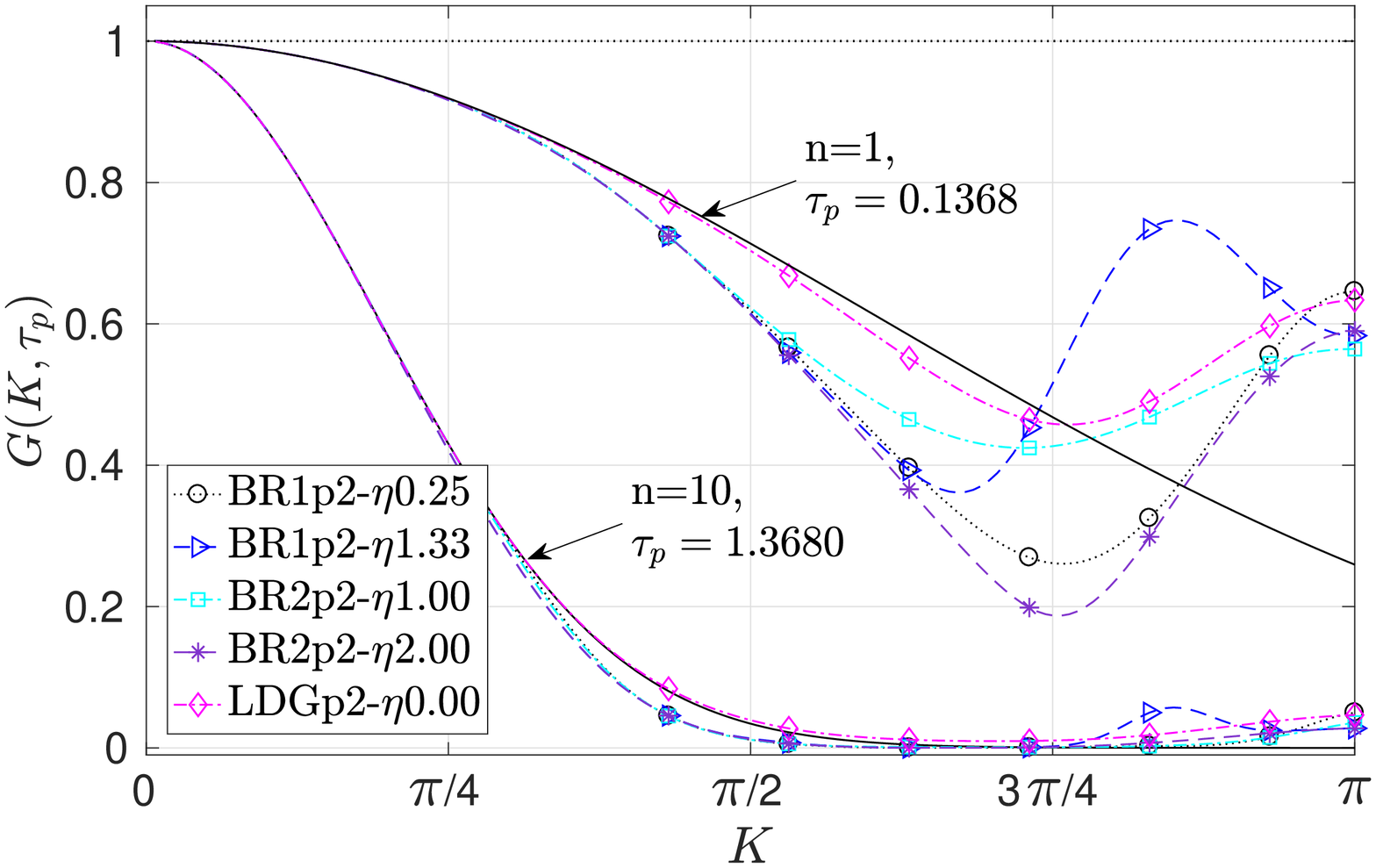} 
    \caption{p$2$-RK$3$, $\Delta \tau=0.0152$} 
    \label{fig:compareP2SchemesRK3_largedTau}
    \end{subfigure}
    \, \,
     \begin{subfigure}[h]{0.5\textwidth}
    \includegraphics[width=0.975\textwidth]{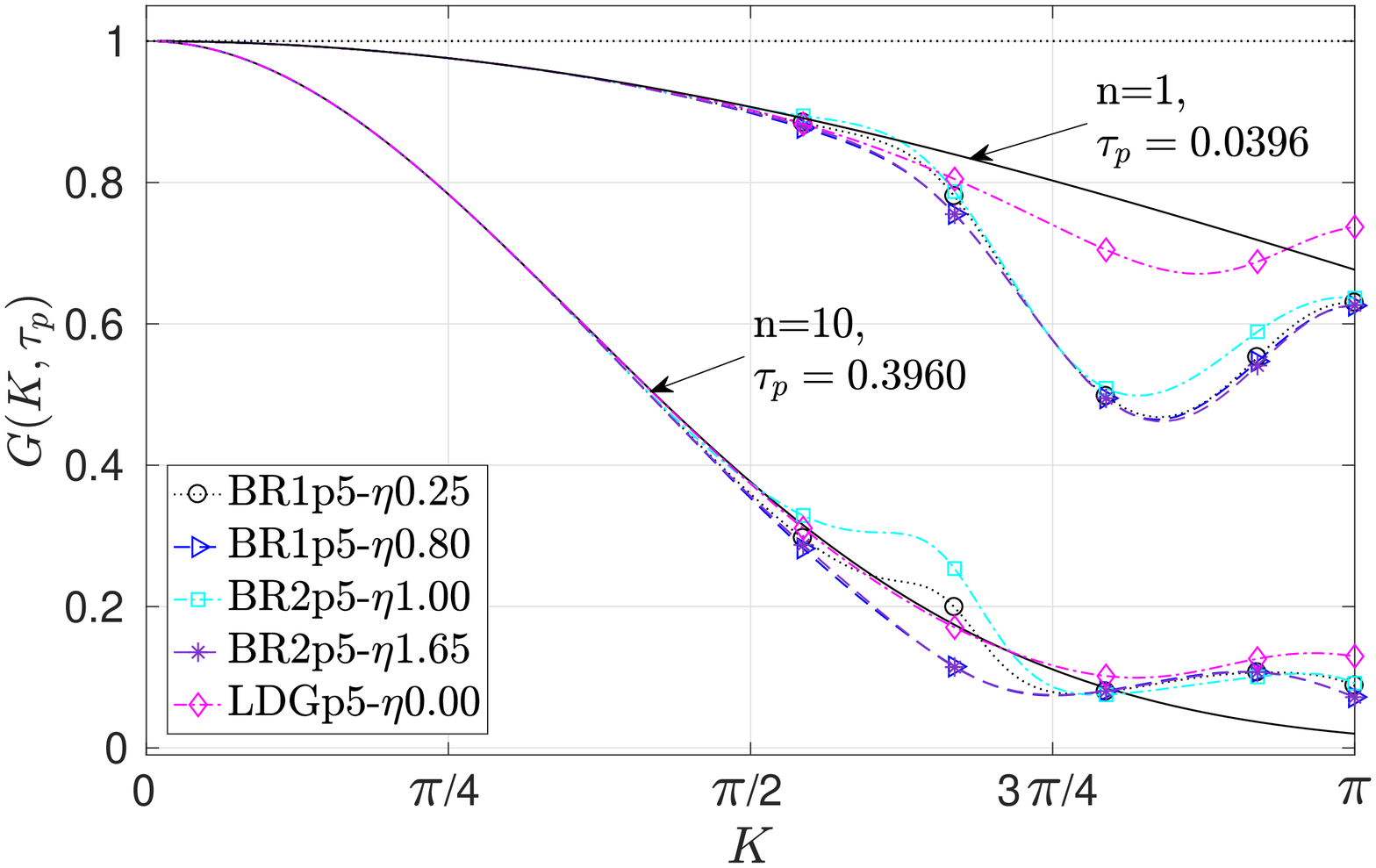} 
    \caption{p$5$-RK$4$, $\Delta \tau=0.0011$} 
    \label{fig:compareP5SchemesRK4_largedTau}
    \end{subfigure} 
    \caption{Comparison of the Diffusion factors of different fully discrete schemes for a fixed large $\Delta \tau \approx 0.9 \times \Delta \tau_{max}^{ldg}$.  In these figures, the solid black lines without symbols represents the exact diffusion factor $G_{ex}=e^{-K^{2} \tau_{p}}$.}
    \label{fig:CompareSchemes_RK4_largedTau}
    \end{figure}%
~\hfigref{fig:CompareSchemes_RK4_largedTau} compares the \mytrue diffusion factor of the three viscous flux formulations with a fixed but relatively large $\Delta \tau$. This $\Delta \tau$ is at about $90\% \times \Delta \tau^{ldg}_{max}$ of the \bff{LDG}-$\eta0$ scheme. Similar to the semi-discrete case, the \bff{LDG}-$\eta0$ is the most accurate scheme for low to moderate wavenumbers, $K\lesssim 3\pi/4$, while it has less dissipation at high wavenumbers for short time simulations. For long time simulations, it attains the best accuracy among all schemes with a similar high wavenumber damping as the other schemes. In addition, the \bff{BR1}-$\eta0.25$ scheme appears to have a good performance for both low and high wavenumbers, similar to the semi-discrete case. The p$2$ version of this scheme is even better than \bff{BR1}-$\eta1.33$, especially for high wavenumbers where the latter experience some oscillations.  Finally, \bff{BR2}p$2$-$\eta2$ scheme has more dissipation than the standard one \bff{BR2}p$2$-$\eta1$ for high wavenumbers $\pi/2 \leq K < \pi$. We can see this also for \bff{BR2}p$5$-$\eta1.65$ which is more dissipative than the \bff{BR2}p$5$-$\eta1$ for $\pi/2 \lesssim K \lesssim 3\pi/4$.

Note that these observations in~\hfigref{fig:CompareSchemes_RK4_largedTau} are at $90\% \times \Delta \tau^{ldg}_{max}$ for the \bff{LDG}-$\eta0.0$ scheme while for the \bff{BR2}-$\eta1.0$ scheme it is at about $\approx 30\%  \times \Delta \tau_{max}^{br2}$. This means that this \bff{BR2}-$\eta1.0$ scheme is almost identical to its semi-discrete version and hence should have better performance especially for high wavenumbers based on the results of the previous section. %
    \begin{figure}[H]
  \begin{subfigure}[h]{0.5\textwidth}
 \centering
    \includegraphics[width=0.965\textwidth]{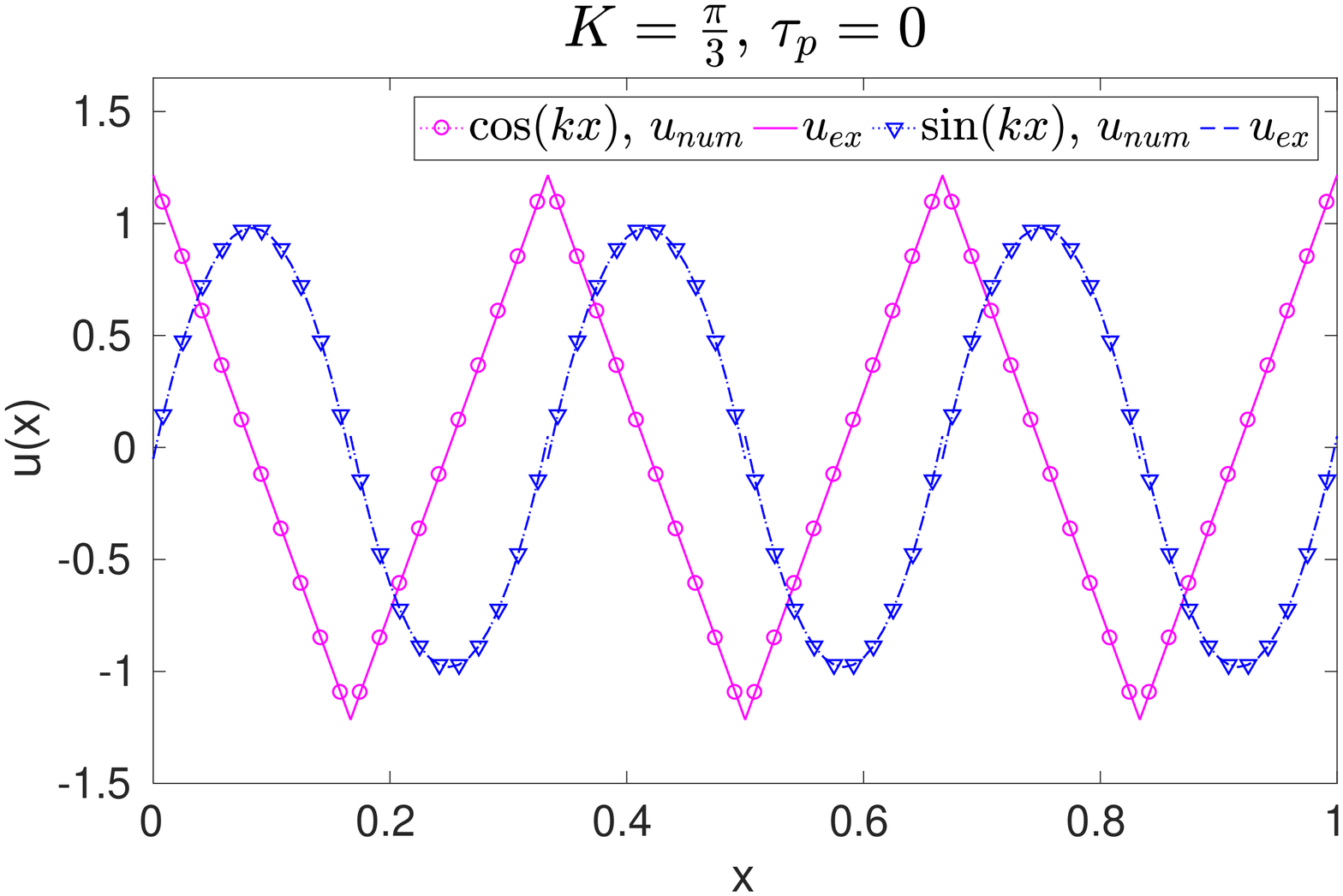}
    \caption{}
    \label{fig:BR2_numwave_sdisc_Kpi_3_t0}
    \end{subfigure} 
    \,
    \begin{subfigure}[h]{0.5\textwidth}
 \centering
    \includegraphics[width=0.965\textwidth]{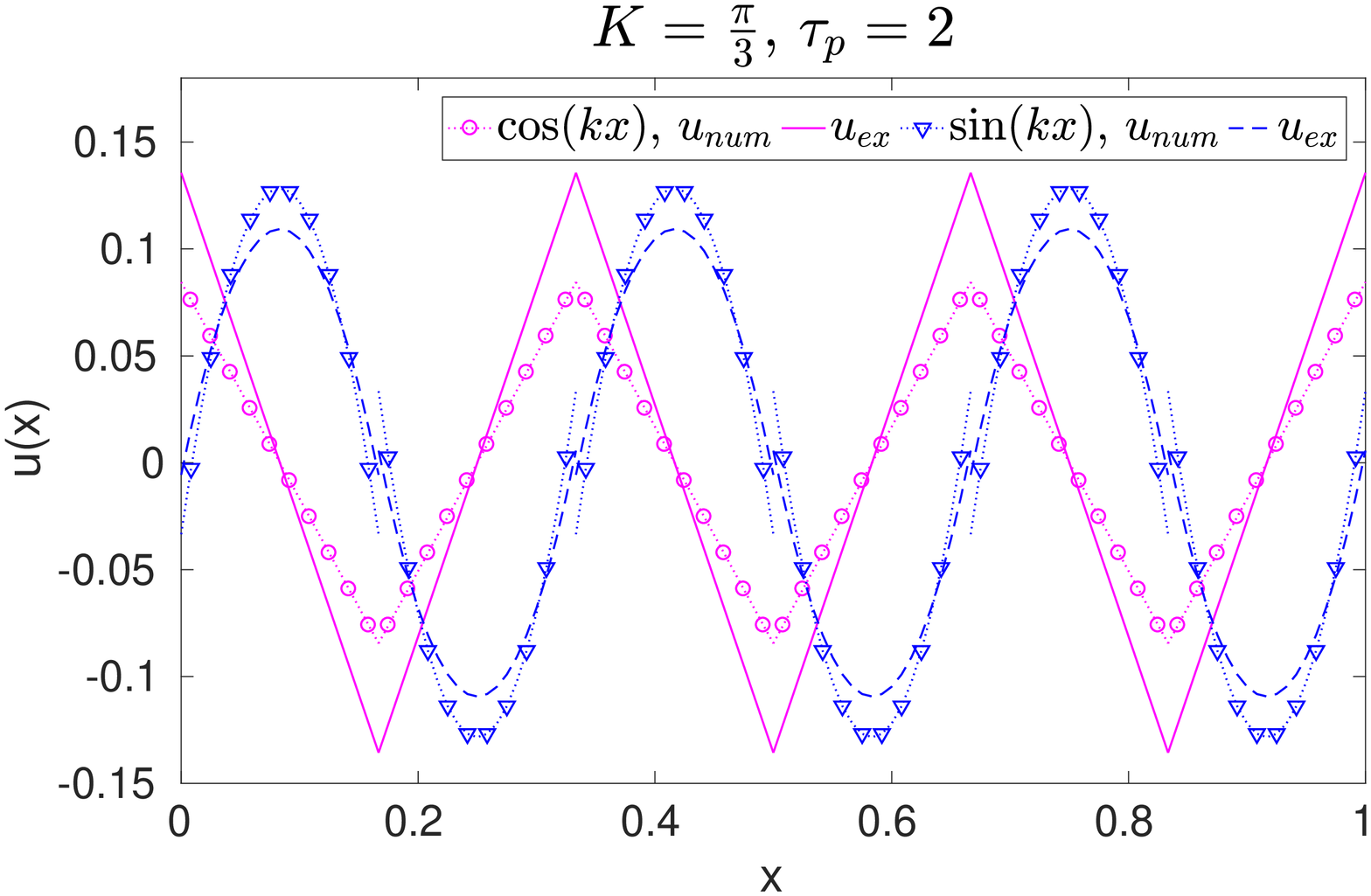} 
    \caption{}
    \label{fig:BR2_numwave_sdisc_Kpi_3_t2}
    \end{subfigure} \\ \\ \\
     \begin{subfigure}[h]{0.5\textwidth}
    \centering
    \includegraphics[width=0.965\textwidth]{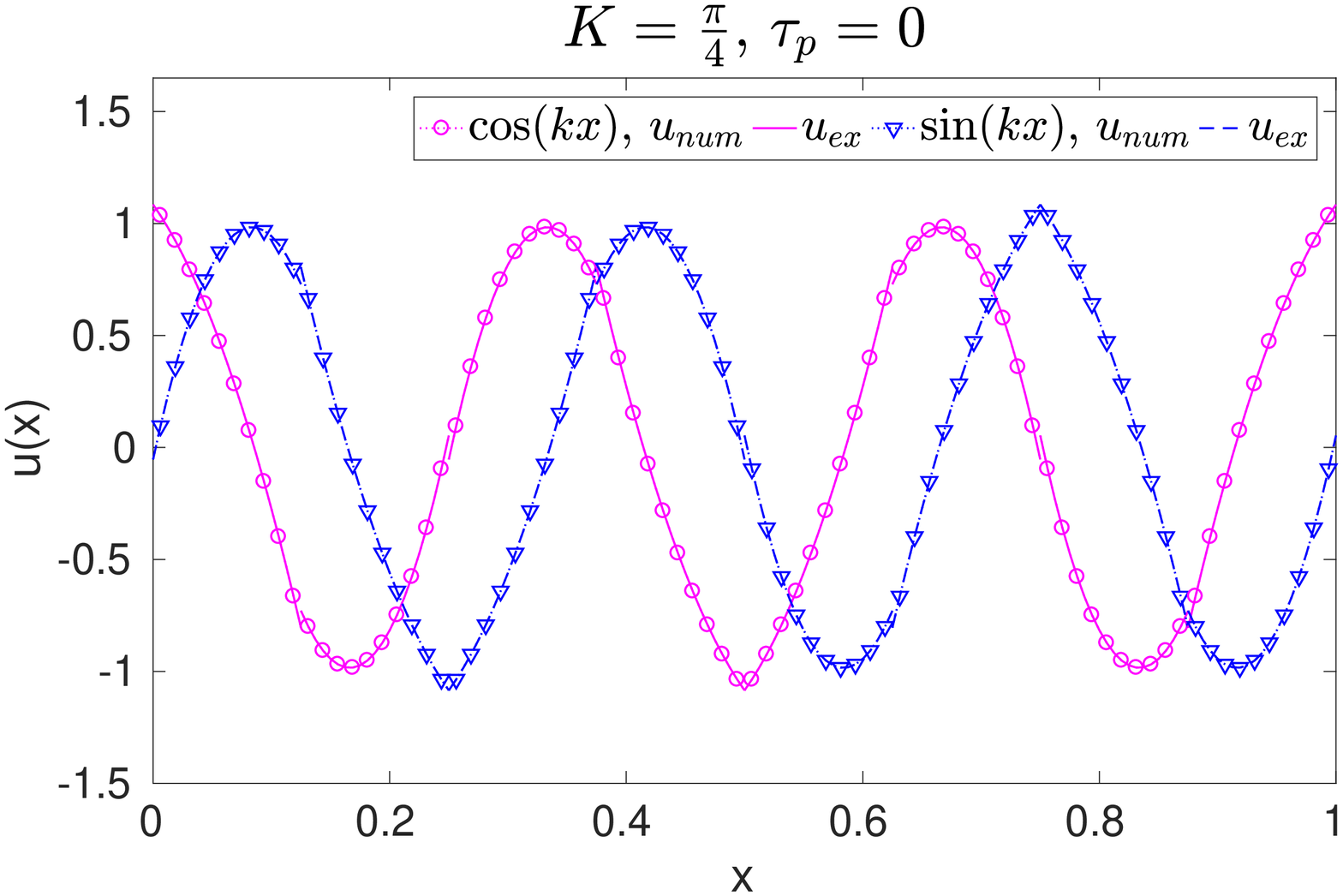} 
    \caption{}
    \label{fig:BR2_numwave_sdisc_Kpi_4_t0}
    \end{subfigure}
    \,
    \begin{subfigure}[h]{0.5\textwidth}
    \centering
    \includegraphics[width=0.965\textwidth]{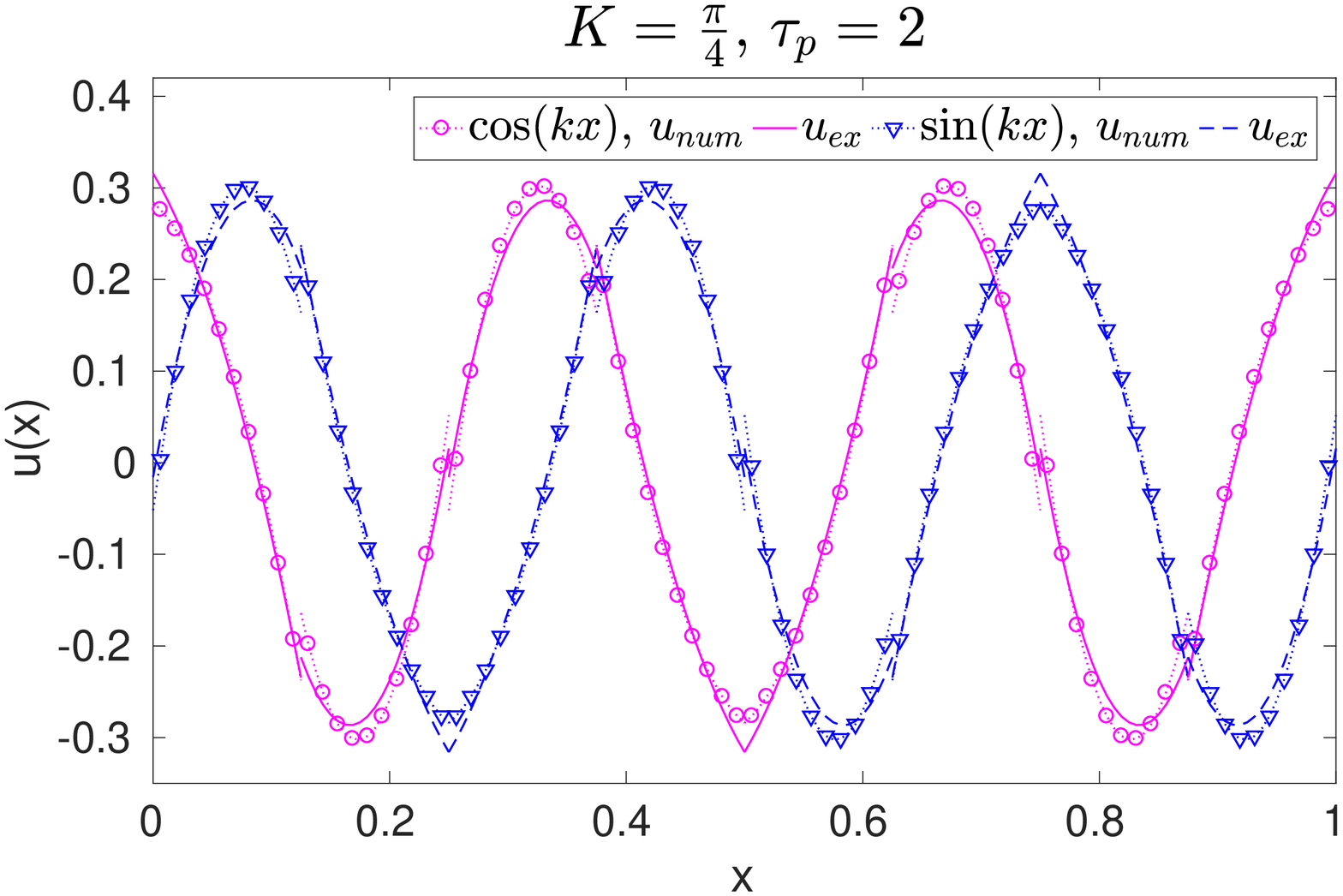} 
    \caption{}
    \label{fig:BR2_numwave_sdisc_Kpi_4_t2}
    \end{subfigure}
    \caption{Simulation of a single Fourier mode (either a sine or a cosine wave) using the \bff{BR2}p2-$\eta1$ scheme and RK$3$ for time integration, $\Delta t=1e-4$. At the first row $K=\pi/3, \: N_{e}=6$ while at the second $K=\pi/4, \: N_{e}=8$. The initial projected solutions ($\tau_{p}=0.0$) are on the left column while the right column is for $\tau_{p}=2.0$.  Note that the scales of each sub-figure are different.}
    \label{fig:BR2_numwave_sdisc}
    \end{figure}%
%
\section{Numerical Results}\label{sec:num_results}%
%
\subsection{Single Fourier mode}\label{subsec:num_Fourierwave}%
%
Consider the following initial condition for the linear-heat~\heqref{eqn:heat_eqn} %
\begin{alignat}{2}
\begin{split}
u(x,0) & = \sin \left( k x/L \right), \quad x \in [0,L]
\end{split}
\label{eqn:sine_wave_IC}
\end{alignat}%
where $k$ denotes the wavenumber, $L$ is the length of the domain, and $\eta=L$ is the wavelength. In all the test cases in this study we let $L=1, \: a=1$, and hence the period of the wave is $T=1$. For a DG scheme, this initial solution is projected using the $L_{2}$ projection onto the space of degree $p$ polynomials on each cell $\Omega_{e}$. 

In order to verify the semi-discrete analysis, a smooth sine/cosine wave (single Fourier-mode) of the form in~\heqref{eqn:sine_wave_IC} is simulated using the \bff{BR1}p$2$-$\eta0.25$, \bff{BR2}p$2$-$\eta1$, and \bff{LDG}p$2$-$\eta0$ schemes. For time integration, we use RK$3$ but with a relatively small $\Delta t=1e-4$ to guarantee being close to the semi-discrete case. For this verification to be consistent we compute the same average energy/amplitude quantity as in~\heqref{eqn:sdisc_E}, but here it is computed numerically for the whole domain%
\begin{alignat}{2}
\begin{split}
E_{num} & = \sqrt{\frac{1}{L} \, \int\limits_{0}^{L}  \left(u^{e}(x,t)\right)^{2} \, \dx } =\sqrt{\frac{1}{L} \sum\limits_{e=0}^{Ne} \, \int\limits_{x_{e-1/2}}^{x_{e+1/2}}  \left(u^{e}(x,t)\right)^{2} \, \dx },   \\ & = \sqrt{\frac{1}{2N_{e}}  \sum\limits_{e=0}^{Ne} \, \int\limits_{-1}^{1} (u^{e}(\xi,t))^{2} \, \dxi } , \quad \text{for uniform grids.}
\end{split}
\label{eqn:E_num}
\end{alignat}%
The numerical diffusion factor is defined in a similar way as%
\begin{equation}
G_{num}(\tau_{p}) = E_{num}(\tau_{p}) / E_{init}, \quad E_{init} = E_{num}(0).
\end{equation}%
In this section, we choose $k= 6 \pi$ and for direct comparison with the results of the~\combined analysis, the corresponding non-dimensional wavenumber of the wave is $K = k L / (N_{e} (p+1))$,  where $N_{e}$ is the number of elements used in the simulations. 

The definition of the average amplitude utilized in the~\combined analysis is based on the projection of a complete Fourier-mode which includes both a sine and a cosine wave. It is expected that both waves should retain the same amplitude after projection since a cosine is a $\pi$ shifted sine wave in principal. However, it turned out that for some under-resolved wavenumbers the $L_{2}$ projection splits the projected amplitude unevenly between the two waves, see for example~\hfigref{fig:BR2_numwave_sdisc_Kpi_3_t0} for a p$2$ projection of $K=\pi/3$ waves.  As a result, the numerical dissipation applied by the numerical scheme to each wave will be different as in~\hfigref{fig:BR2_numwave_sdisc_Kpi_3_t2}. In order to compare the results of the~\combined analysis with numerical simulations for such wavenumbers, one should compute the diffusion factor based on the energy of the two numerical simulations, i.e., one for a cosine wave and one for a sine wave. In other words, if we have a Fourier-mode of the form%
\begin{equation*}
e^{ikx} = \cos(kx) + i \sin(kx)  \rightarrow E_{num} = \sqrt{ E_{\cos}^{2} + E_{\sin}^{2}}
\end{equation*}%
which is numerically equivalent to the energy definitions in~\heqsref{eqn:sdisc_E}{eqn:sdisc_Eex}. As the order of the scheme gets higher and higher or as the mesh is refined (moving towards the resolved range of the wavenumbers), this difference in treating both wave forms is reduced and both waves are conceived by the numerical scheme in the same way.  Another wavenumber that is simulated for verification is $K=\pi/4$ (sufficiently resolved). The results are presented in~\hfigsref{fig:BR2_numwave_sdisc_Kpi_4_t0}{fig:BR2_numwave_sdisc_Kpi_4_t2}. For this particular case there is no difference in the projected energy/amplitude between the two wave forms and hence a direct comparison with either one of them can be made directly. %
\begin{table}[H]%
\caption{Numerical and~\combined analysis results for the simulation of a single Fourier-mode $e^{i k x}$, $K= \pi/3$ at $\tau_{p}=2.0$. The initial projected energy/average amplitude of this wave $E_{init}=0.9962$. This simulation is conducted using three p$2$ DG schemes with RK$3$ for time integration, $\Delta t=1e-4$, and $N_{e}=6$.} 
\centering 
\small
\begin{tabular}{|c|c c c c c|c|}
\hline \multirow{2}{*}{Scheme}  & \multicolumn{5}{c|}{numerical simulation} &   \multicolumn{1}{c|}{~\combined analysis} \\
\cline{2-7}
  & $E_{cos}$ & $E_{sin}$ & $E_{tot}$ & $G_{num}$ & $|\Delta G_{num}|$ & $|\Delta G_{true}|$  \\
 \hline   \bff{BR2}p$2$-$\eta1.00$ & $0.0488$ & $0.0887$ & $0.1012$ & $0.1016$ & $9.93e-3$ & $ 9.91e-2$  \\
  \hline   \bff{BR1}p$2$-$\eta0.25$ & $0.0488$ & $0.0939$ & $0.1058$ & $0.1062$ & $5.32e-3$ & $ 5.31e-3$  \\
 \hline   \bff{LDG}p$2$-$\eta0.00$ & $0.0786$ & $0.0803$ & $0.1124$ & $0.1128$ & $1.24e-3$ & $ 1.22e-3$  \\
\hline
\end{tabular}
\label{table:numwave_compareP2Schemes}
\end{table}%
In general, for both wavenumbers we were able to recover the same energy and diffusion factors of the semi-discrete~\combined analysis as well as diffusion errors  with very good accuracy. The results for the $K=\pi/3$ case are presented in~\htablref{table:numwave_compareP2Schemes}. These results can be compared with~\hfigref{fig:BR2_numwave_sdisc_Kpi_3_t2} for the \bff{BR2}p$2$-$\eta1$ numerical results as well as~\hfigref{fig:sdisc_compareP2schemes_tau2_Gerr} for the true diffusion errors comparison. The absolute errors $|\Delta G|$ in this table are defined in a similar way as in~\heqref{eqn:true_Gerr}.%
%
\subsection{Gaussian wave}\label{subsec:num_gaussian}%
\begin{figure}[H]
 \centering
    \includegraphics[width=0.975\textwidth]{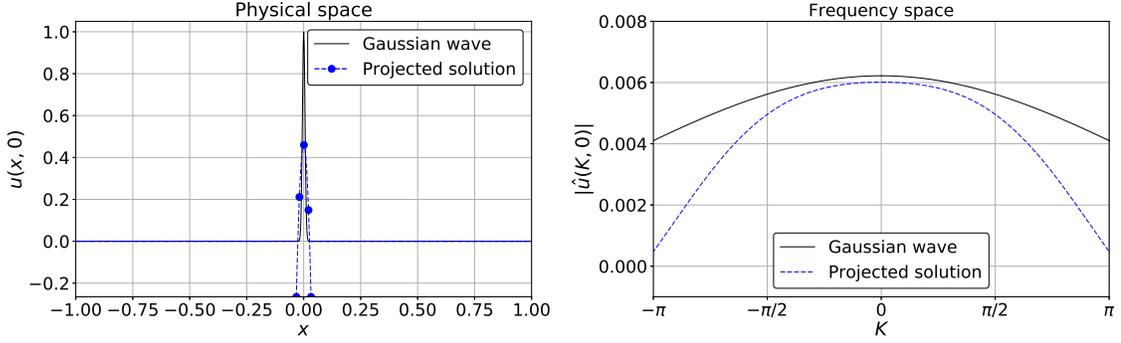}
    \caption{Representation of the initial Gaussian wave solution along with its p$2$ projection. In the left figure, the initial solution and its projection are plotted in the physical space $x$ while in the right figure, are their FFT results.  }
    \label{fig:gaussian_initial_sol}
    \end{figure}

\begin{figure}[H]
  \begin{subfigure}[h]{0.5\textwidth}
 \centering
    \includegraphics[width=0.975\textwidth]{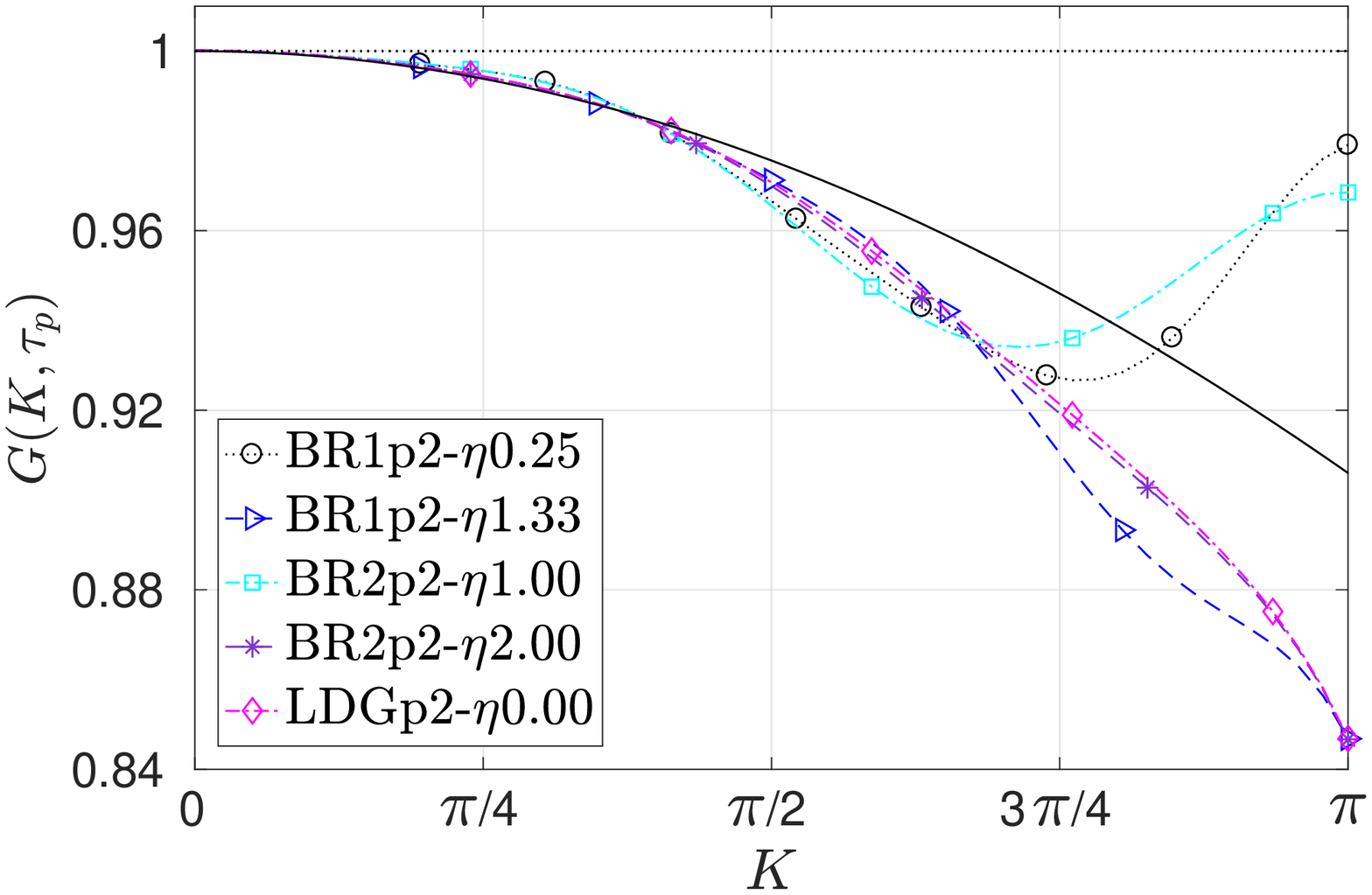}
    \caption{$\tau_{p}=0.01$}
    \label{fig:num_compareP2schemes_sharp_gaussian_t0.01}
    \end{subfigure} 
    \,
    \begin{subfigure}[h]{0.5\textwidth}
 \centering
    \includegraphics[width=0.975\textwidth]{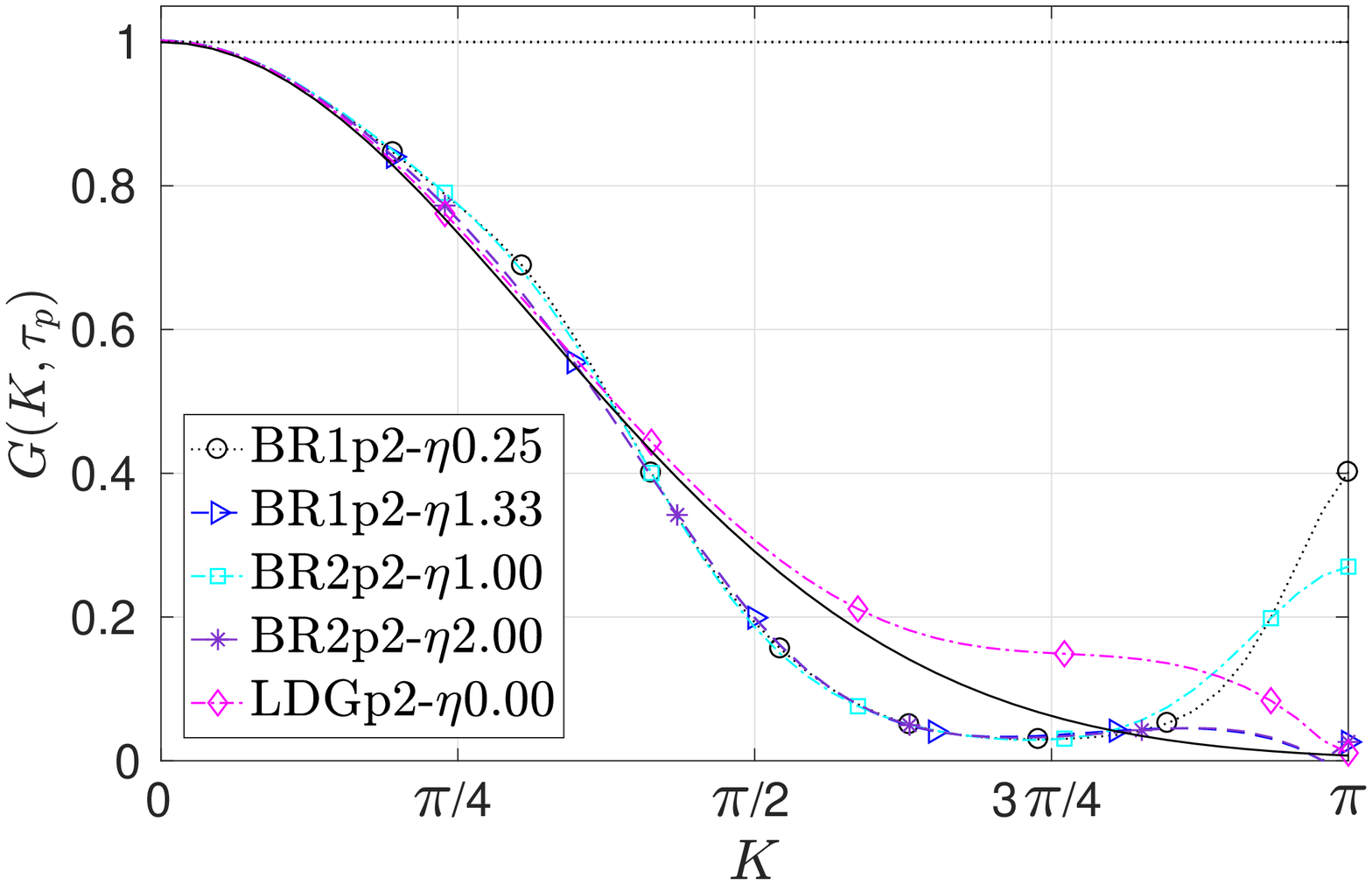} 
    \caption{$\tau_{p}=0.5$}
    \label{fig:num_compareP2schemes_sharp_gaussian_t0.5}
    \end{subfigure} 
    \caption{Numerical simulation of a Gaussian wave. This figure presents the comparison of the diffusion factors for several p$2$ schemes. These simulations were conducted using RK$3$ for time integration with small $\Delta t=1e-4$. Note that the vertical scale of each sub-figure is different.}
    \label{fig:num_compareP2schemes_sharp_gaussian}
    \end{figure}

%
In this section we compare the performance of several DGp$2$ diffusion schemes coupled with RK$3$ for a simulation of a broadband Gaussian wave. This helps in quantifying the numerical errors associated with a wide range of wavenumbers for a particular scheme. Consider an initial solution for the linear-heat~\heqref{eqn:heat_eqn} of the form%
\begin{equation}
u(x,0) = e^{-b \: x^{2}},\quad, \: x \: \epsilon \: [-L,L] ,
\label{eqn:Gaussian_function_IC}
\end{equation}%
where $2L$ is the length of the domain, with $L=1,\;a=1$. The exact solution can be computed analytically using the method of separation of variables for linear partial differential equations%
\begin{alignat}{2}
u_{ex}(x,t) = a_{0} + \sum\limits_{m=1}^{\infty} a_{m} \, \cos\left( m \pi / L \right) e^{-\gamma m^{2}  t}, 
\label{eqn:gaussian_exact_sol1}
\end{alignat}%
where $a_{0},\, a_{m}$ are constant coefficients provided in~\happref{appx:B}{B}. For this case, the parameter $b=1.5e4$ is chosen as to yield a sufficiently smooth case with a wide range of wavenumbers for the numerical simulations to be as close as possible to the~\combined analysis results. This results in a very sharp Gaussian wave in the spatial domain while its Fourier transform is very wide in the frequency space (under-resolved case). Similar to the previous test case, the initial condition is projected using the $L_{2}$ projection onto the space of degree $p$ polynomials on each cell $\Omega_{e}$.~\hfigref{fig:gaussian_initial_sol} displays both the Gaussian wave and the projected initial solution in addition to their Fast Fourier Transforms (FFT). We note the input to the FFT is a continuous solution on a uniform grid, and this continuous solution is obtained by averaging the interface solutions between two elements in addition to using a large number of points ($ > $nDOF) inside each element. This procedure is utilized as to converge the FFT results to a unique distribution with minimized errors. However, as we can see from~\hfigref{fig:gaussian_initial_sol}, there are still some differences between the Gaussian wave FFT and the projected solution FFT even for $K=0$. We attribute these differences and errors to both the projection error as well as the error introduced by the jumps at interfaces which can not be excluded by averaging. %
This test case is simulated using the same group of schemes used in previous sections with $p2$ polynomial orders. For marching in time we utilize RK$3$ with a sufficiently small $\Delta t=1e-04$ to be as close as possible to a semi-discrete case.~\hfigref{fig:num_compareP2schemes_sharp_gaussian} presents the~\mytrue diffusion factor results through numerical simulation for both short and long time simulations. These diffusion factors are computed based on the wave energy for each wavenumber determined using a Fast Fourier transform (FFT) algorithm. In this case the energy is simply defined as%
\begin{equation}
E\left( K;\tau_{p} \right)= |\hat{u} \left( K;\tau_{p} \right)|^{2}, 
\end{equation}%
where $|\hat{u} \left(K;\tau_{p} \right)|$ is the amplitude of the wave as determined by the FFT. The results in~\hfigref{fig:num_compareP2schemes_sharp_gaussian} show good agreement with the semi-discrete results in~\hfigref{fig:sdisc_compareP2schemes} for the diffusion factors at $\tau_{p}=0.01, \, \tau_{p}=0.5$. However, these two sets of figures are not exactly identical. The differences are attributed to the errors associated with the FFT results as pointed out in the previous discussion about the initial energy distribution. Nevertheless,  the numerical results provide the same general comparison results between different schemes and hence verify the applicability of the~\combined analysis for assessing the diffusion behavior of DG schemes. Indeed,  the \bff{LDG} method provides the most accurate scheme in approximating the exact diffusion behavior especially for long time simulations. %
    %
\subsection{The Burgers turbulence case}\label{subsec:num_burgers}%
%
Consider the viscous Burgers equation %
\begin{equation}
\frac{\partial u}{\partial t} + \frac{\partial F(u)}{\partial x} = \gamma \frac{\partial^{2} u}{\partial x^{2}}, \quad F = u^{2}/2 ,
\label{eqn:burgers_eqn}
\end{equation}%
where $F$ is the nonlinear inviscid flux function. For a randomly generated velocity distribution from an initial energy spectrum the solution of this equation undergoes a chaotic (turbulence like) behavior through time. In the literature, this is often referred to  as the decaying Burgers turbulence since it features a kinetic energy spectrum $E \propto k^{-2}$ with an energy decaying phases that reassembles the Navier-Stokes decaying turbulence behavior to some extent. 

In this work, we utilize this case in a slightly different manner than what is usually performed in the literature~\cite{DeStefanoSharpcutoffsmooth2002,AdamsSubgridScaleDeconvolutionApproach2002,SanAnalysislowpassfilters2016,AlhawwaryFourierAnalysisEvaluation2018}. That is by requiring  a small Pe\'clet ($\text{P\!e}$) number in addition to a higher viscosity value of $\gamma =0.015$ in order to make the diffusion part of the problem more dominant. Consequently, a direct comparison with the~\combined analysis for diffusion can be established. The Pe\'clet number is defined as $\text{P\!e} = E h / \gamma = u_{rms} h / \gamma$ and $E$ is defined as the total average energy or average amplitude in the domain~\heqref{eqn:E_num} which is mathematically equivalent to the $u_{rms}$ the root mean square of $u$. Based on this definition it is clear that the $\text{P\!e}$ will change with time as the velocity distribution changes. Therefore, we focus in this case on long time behavior where the $\text{P\!e}$ is ensured to be $<1$. 

Assuming a periodic domain $x\in[0,2\pi]$, we conduct the simulation of a decaying Burgers turbulence case similar to the case used in~\cite{AlhawwaryFourierAnalysisEvaluation2018}, with an initial energy spectrum%
\begin{equation}
E(k;0) = E_{0}(k) = A k^{4} \rho^{5} e^{-k^{2} \rho^{2} } ,
\label{eqn:energ_spect_omersan}
\end{equation}%
where $k$ is the prescribed wavenumber,  $(\rho, \, A)$ are constants to control the position of the maximum energy, and for $\rho =10, A = 2/(3\sqrt{ \pi}) $ the spectrum reaches its maximum at $k=13$. Assuming a Gaussian distribution in the frequency space, the initial velocity field reads%
\begin{equation}
\hat{v}(k) = \sqrt{2 E_{0}(k)} e^{i 2 \pi \Phi(k)} ,
\label{eqn:burger_velocity_fourierspace}
\end{equation}%
where $\Phi(k)$ is a random phase angle that is uniformly distributed in $ \left[ 0,1 \right] $ for each wavenumber $k$. In the physical space, an inverse Fourier transform can be employed to yield a real velocity field, provided that $\Phi(k) = -\Phi(-k)$,%
\begin{equation}
v(x) = \sum_{k=0}^{k_{max}} \sqrt{2 E_{0}(k)} \cos \left( k x + 2 \pi \Phi(k) \right),
\label{eqn:burger_velocity_physicalspace}
\end{equation}%
where $k$ is the prescribed integer wavenumber with $ k_{max}=2048$. 

This initial solution is projected onto the DG space of polynomials of degree at most $p$ using an $L_{2}$ projection methodology. For the convection terms we employ a standard DG discretization in space with the same polynomial degree as the diffusion terms and Rusanov~\cite{ToroRiemannSolversNumerical2009} numerical (upwind) fluxes are used at interfaces between cells. The diffusion terms are discretized by any scheme of the viscous diffusion schemes introduced in~\secref{sec:viscflux_formulations}. For both types of terms, exact integration is always performed to mitigate aliasing errors.

We perform simulations for the same group of schemes that were utilized in the Gaussian wave simulation. For all schemes, the Burgers equation~\heqref{eqn:burgers_eqn} is solved for a $64$ randomly generated samples of the initial velocity field~\heqref{eqn:burger_velocity_physicalspace}. After that an ensemble averaged FFT is performed for the $64$ cases to compute a sufficiently smooth kinetic energy spectrum. The discretization of the domain involves $N_{e}=50$ elements (nDOF$=150$) which results in an initial $\text{P\!e} \approx 2.1$ for all simulations. This $\text{P\!e}$ number quickly decreases with time allowing for a diffusion dominated behavior for long time simulations. For time integration we employ RK$3$ with a sufficiently small $\Delta t=1e-4$ for the time integration errors to be negligible. A DNS simulation was conducted using $N_{e}=500$ elements to provide a reference solution for comparisons. 

\begin{figure}[H]
 \centering
    \includegraphics[width=0.725\textwidth]{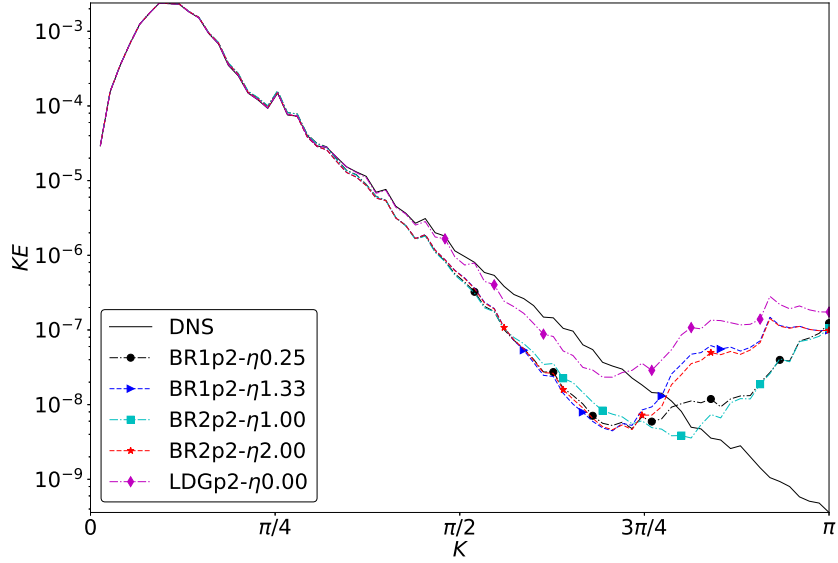} 
    \caption{Numerical results of the Kinetic energy ($K\!E$) spectra for the Burgers turbulence case at $t=0.5$, $\text{P\!e}\approx 0.15$. For time integration RK$3$ scheme is utilized and nDOF$=150$ for all schemes. }
    \label{fig:Burgers_KE_t0.5}
    \end{figure}

 The Kinetic energy spectra ($K\!E$) for all considered p$2$ schemes at $t=0.5$ with RK$3$ for time integration are presented in~\hfigref{fig:Burgers_KE_t0.5}. This point in time is chosen such that its $\text{P\!e} \lesssim 0.2$ for all the $64$ sampled simulations so as to ensure the domination of the diffusive effects.  From this figure, it can be inferred that indeed for long time simulations the \bff{LDG}-$\eta0$ scheme is the most accurate in approximating the exact diffusion behavior since it captures more $K\!E$ than all the other schemes. However, for high wavenumber it is less dissipative than all the other schemes but since this is a common behavior for all DG schemes~\cite{AlhawwaryFourierAnalysisEvaluation2018} we think that its accuracy outweigh this problem. All the other schemes are comparable in terms of low to moderate wavenumber accuracy in contrast to some differences at high wavenumbers. Indeed, the two equivalent schemes identified through the~\combined analysis are showing the same behavior for this case as well, namely, schemes \bff{BR1}p$2$-$\eta1.33$ and \bff{BR2}p$2$-$\eta2$. In addition, the \bff{BR1}p$2$-$\eta0.25$ scheme shows a good behavior that is close to the \bff{BR2}p$2$-$\eta1$and thus it is indeed enhancing the stability of the original \bff{BR1}p$2$ scheme with $\eta=0$. 
 
\begin{figure}[H]
 \centering
    \includegraphics[width=0.725\textwidth]{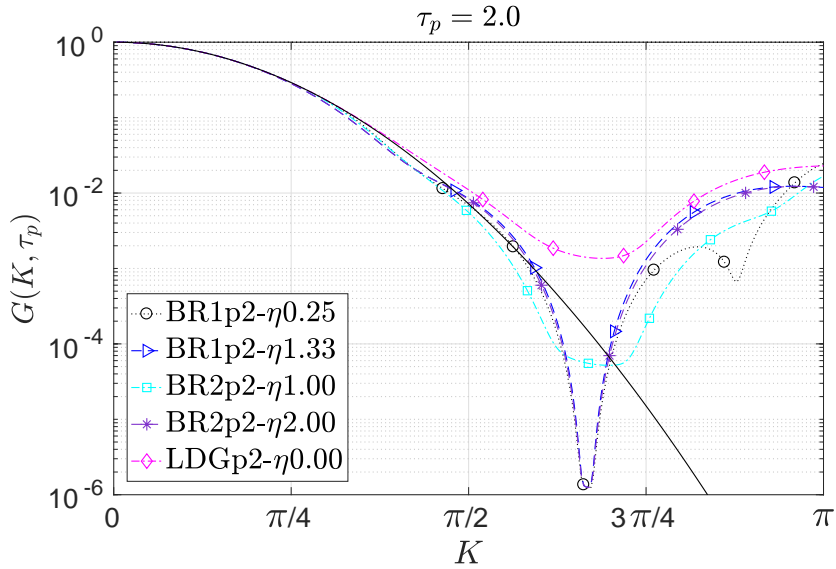}
    \caption{Diffusion factors of several p$2$ diffusion schemes through~\combined analysis. The solid black line without symbols represents the exact diffusion factor $G_{ex}=e^{-K^{2} \tau_{p}}$.  These results are to be compared with~\hfigref{fig:Burgers_KE_t0.5}. Note that the vertical axis is in log scale. }
    \label{fig:compareP2Schemes_tau2_burgers}
    \end{figure}

In order to compare this case numerical results with the~\combined analysis we have plotted the diffusion factors of all the considered schemes at a very long time $\tau_{p}=2.0$ in a similar log scale, see~\hfigref{fig:compareP2Schemes_tau2_burgers}. Note that this figure presents the comparisons for diffusion factors of a pure diffusion problem, while the $K\!E$ spectra of a nonlinear Burgers case are presented in~\hfigref{fig:Burgers_KE_t0.5}. Despite these differences, we can clearly see how the~\combined analysis gives accurate and useful predictions about the behavior of DG viscous schemes in a simple and informative manner. We also note that the short time behavior of such schemes is totally different from their long time behavior as was predicted through the~\combined analysis.%
    %
\section{Conclusions}\label{sec:conclusions}%
%
In this paper, we studied the semi-discrete and fully-discrete behavior of several viscous flux formulations in the context of DG methods in terms of stability and dissipation characteristics. Insights on the performance of each scheme were inferred through the \ita{combined-mode} analysis. It was shown that the classical definition of the~\physic cannot be utilized for studying the behavior of DG schemes for diffusion problems over the entire wavenumber range. The  \ita{combined-mode} analysis provides a more realistic picture of the diffusion characteristics for the entire wavenumber range in both semi-discrete and fully-discrete forms.

In studying the effects of the penalty parameters on the performance of several DG viscous schemes,  it was found that the minimum $\eta$ for stability should be avoided for all schemes. This is because these schemes with $\eta_{min}$ have a very slow dissipation rate in comparison to the exact one which may lead to unstable behavior for nonlinear cases. In addition, for short time simulations, schemes \bff{BR1}-$\eta0$, \bff{BR2}-$\eta1$ have a slower decaying rate than the exact one for high wavenumbers $K \gtrsim 3\pi/4$ in contrast to \bff{LDG}-$\eta0$. This short time behavior was demonstrated for all polynomial orders $p<5$. For long time simulations, all schemes including the \bff{LDG}-$\eta0$ have a slower than exact dissipation rate for $K \gtrsim 3\pi/4$ and up to $p=5$. In general, high order schemes are more robust than low order ones due to their high wavenumbers diffusion. 

Moreover, by introducing a stabilized version of the \bff{BR1} method we were able to identify \bff{BR1} schemes, with some $\eta$, that makes them close in their diffusion behavior to some \bff{BR2} schemes. For instance, the \bff{BR1}p$2$-$\eta1.33$ and \bff{BR2}p$2$-$\eta2$ schemes provide almost the same dissipation in both the semi-discrete and the fully-discrete cases, although they are closer in the semi-discrete case. This was also demonstrated for $p5$ schemes. Adjusting the penalty parameter $\eta$ was shown to have a strong influence on the scheme behavior for short time simulations. By comparing different schemes, we were able to show that the \bff{LDG} method has the best accuracy in approximating the exact dissipation behavior especially for long time simulations besides maintaining a lower bound on the diffusion error. 

For fully discrete schemes, maximum time steps for the stability were provided. It was also shown that at about $50\% \Delta \tau_{max}$ all schemes recovered their semi-discrete behavior while near $\Delta \tau_{max}$ they all experienced slower than exact decaying rates for high wavenumbers. 

The~\combined analysis results were verified numerically for the linear heat equation. In addition, a decaying Burgers turbulence case was utilized in order to assess the performance of the considered schemes for a more practical case. It was shown that indeed the \bff{LDG}-$\eta0$ method provides the best accuracy in capturing an energy spectrum among all other schemes if  the same small $\Delta \tau$ is used  for all schemes. This case also verified that long time behavior is always very different than short time behavior. %
\begin{acknowledgement}
The research outlined in the present paper has been supported by AFOSR under grant FA9550-16-1-0128, and US Army Research Office under grant W911NF-15-1-0505.
\end{acknowledgement}%
%
\section*{Appendix A. Stability limits for RKDG schemes}\label{appx:A}%
\addcontentsline{toc}{section}{Appendix A. Stability limits for DG schemes}%
%
\begin{table}[H]%
\caption{ Stability limits ($\Delta \tau_{max}$) for \bff{BR2} schemes with different stabilization parameter $\eta$ and RK time integration schemes for the $1$D linear heat equation. } 
\centering 
\small
\begin{tabular}{c c c c}
\hline
DG, $p$ & RK, $s$ & $\eta$ & $\Delta \tau_{max}$ \\ 
\hline
 \multirow{6}{*}{$1$} & \multirow{2}{*}{$2$} & $0.50-0.80$ & $0.1666$ \\ 
 & & $1.00$ & $0.1498$ \\
  \cline{2-4} 
 &  \multirow{2}{*}{$3$} & $0.50-0.80$ & $0.2093$ \\ 
 & & $1.00$ & $0.1882$ \\
   \cline{2-4} 
& \multirow{2}{*}{$4$} & $0.50-0.80$ & $0.2321$ \\ 
 & & $1.00$ & $0.2086$ \\
\hline
 \multirow{2}{*}{$2$} & $3$& $0.67-1.10$ & $0.0418$ \\ 
& $4$ & $0.67-1.10$ & $0.0464$ \\ 
\hline
 \multirow{2}{*}{$3$} & $3$& $0.75-1.10$ & $0.0147$ \\ 
& $4$ & $0.75-1.10$ & $0.0163$ \\ 
\hline
 \multirow{2}{*}{$4$} & $3$& $0.80-1.10$ & $0.0066$ \\ 
& $4$ & $0.80-1.10$ & $0.0073$ \\ 
\hline
 \multirow{2}{*}{$5$} & $3$& $0.84-1.10$ & $0.0034$ \\ 
& $4$ & $0.84-1.10$ & $0.0037$ \\ 
\hline
\end{tabular}
\label{table:BR2_stability_limits_1D}
\end{table}%
%
\begin{table}[H]%
\caption{ Stability limits ($\Delta \tau_{max}$) for all schemes in their standard forms with RK time integration schemes for the $1$D linear heat equation. } 
\centering 
\small
\begin{tabular}{c c| c c c|c}
\hline
DG, $p$ & RK, $s$  & \bff{BR1}-$\eta0$ & \bff{BR2}-$\eta1$& \bff{LDG}-$\eta0$ & $(\Delta \tau)_{br2} / (\Delta \tau)_{ldg}$\\ 
\hline
 \multirow{3}{*}{$1$} & $2$ & $0.1250$  &  $0.1498$ & $0.0555$ & $\approx 2.7$ \\ 
 & $3$ & $0.1570$ &  $0.1882$ &  $0.0697$ & \\
 & $4$  & $0.1740$ & $0.2086$ & $0.0773$ & \\ 
\hline
 \multirow{2}{*}{$2$}  & $3$ &  $0.0384$ & $0.0418$ & $0.0169 $ & $\approx 2.5$\\
 & $4$  & $0.0426$ & $0.0464$ & $0.0187$ &  \\ 
\hline
 \multirow{2}{*}{$3$}  & $3$ &  $0.0142$ & $0.0147$ & $0.0057$ & $\approx 2.6$\\
 & $4$  & $0.0158$ & $0.0163$ & $0.0063$ & \\ 
\hline
 \multirow{2}{*}{$4$}  & $3$ &  $0.0064$ & $0.0066$ &  $0.0024$ & $\approx 2.8$\\
 & $4$  & $0.0071$ & $0.0073$ & $0.0026$ & \\ 
\hline
 \multirow{2}{*}{$5$}  & $3$ &  $0.0033$ &  $0.0034$ & $0.0011$ & $\approx 3.1$\\
 & $4$  & $0.0037$ & $0.0037$ & $0.0012$ & \\ 
\hline
\end{tabular}
\label{table:BR1BR2LDG_stability_limits_1D}
\end{table}%
%
%
\section*{Appendix B. Exact solution for the linear heat equation with a Gaussian initial condition}\label{appx:B}%
\addcontentsline{toc}{section}{Appendix B. Exact solution for the linear heat equation with a Gaussian initial condition}%
%
For the following Gaussian initial condition%
\begin{equation*}
u(x,0) = e^{-b \: x^{2}},\quad, \: x \: \epsilon \: [-L,L] ,
\end{equation*}
the exact solution for the linear heat equation~\heqref{eqn:heat_eqn} with periodic boundary conditions, takes the form%
\begin{equation*}
u_{ex}(x,t) = a_{0} + \sum\limits_{m=1}^{\infty} a_{m} \, \cos\left( m \pi / L \right) e^{-\gamma m^{2}  t}, 
\end{equation*}where ,
\begin{alignat*}{2}
\begin{split}
a_{0} &= \frac{\sqrt{\pi}}{2c}, \quad a_{m} = \frac{\sqrt{\pi}}{c} \, \mathcal{R}e\left( \text{erf}(z_{m}) \right) e^{-m^{2}/(4b)} , \\
z_{m} &= c + i \frac{m \pi}{2c} , \quad c= L\sqrt{b}.
\end{split}
\label{eqn:gaussian_exact_sol_param1}
\end{alignat*}%
In the above equation $\text{erf}(.)$ is the error function operating on a complex number $z$ and $\mathcal{R}e(.)$ is the real part of a complex number. This form of the exact solution is drived through the method of the separation of variables for linear partial differential equations. 

\bibliographystyle{spmpsci}      
\bibliography{references}   
\addcontentsline{toc}{section}{References}%

\end{document}